\documentclass[pdflatex,sn-mathphys-num]{sn-jnl}

\usepackage{graphicx}%
\usepackage{multirow}%
\usepackage{amsmath,amssymb,amsfonts}%
\usepackage{amsthm}%
\usepackage{mathrsfs}%
\usepackage[title]{appendix}%
\usepackage{xcolor}%
\usepackage{textcomp}%
\usepackage{manyfoot}%
\usepackage{booktabs}%
\usepackage{algorithm}%
\usepackage{algorithmicx}%
\usepackage{algpseudocode}%
\usepackage{comment}
\usepackage{listings}%
\usepackage{bbm}

\theoremstyle{thmstyleone}%
\newtheorem{theorem}{Theorem}[section]
\newtheorem{lemma}{Lemma}[section]
\newtheorem{proposition}{Proposition}[section]%

\theoremstyle{thmstyletwo}%
\newtheorem{example}{Example}[section]%
\newtheorem{remark}{Remark}[section]%
\newtheorem{corollary}{Corollary}[section]%
\newtheorem{assumption}{Assumption}[section]%
\newtheorem{definition}{Definition}[section]%
\usepackage{url}
\usepackage{color,soul}

\begin{document}

\title[Article Title]{Centralized Reduction of Decentralized Stochastic Control with General Spaces: Existence, Weak-Feller Regularity, and Approximate Optimality}


\author*[1]{\fnm{Omar} \sur{Mrani-Zentar}}\email{o.mranizentar@queensu.ca}

\author[1]{\fnm{Serdar} \sur{Y\"uksel}}\email{yuksel@queensu.ca}
\equalcont{These authors contributed equally to this work.}

\affil*[1]{\orgdiv{Department of Mathematics and Statistics}, \orgname{Queen's University}, \orgaddress{\street{48 University Avenue}, \city{Kingston}, \postcode{K7L 3N6}, \state{Ontario}, \country{Canada}}}


\abstract{Stochastic control problems involving agents with decentralized information structures and general state, measurement, and action spaces are intrinsically difficult to study because of the inapplicability of standard tools from centralized (single-agent) stochastic control. In this paper, we address some of these challenges for decentralized stochastic control with standard Borel spaces under two different but tightly related information structures: the one-step delayed information sharing pattern, and the $K$-step periodic information sharing pattern. We will show that  problems with these information patterns can be reduced to a centralized Markov Decision Process (MDP), generalizing prior results which considered finite, linear, or static models, by addressing several measurability and topological questions. We then provide sufficient conditions for the transition kernels of both centralized reductions to be weak-Feller. The existence and separated nature of optimal policies under both information structures are then established. The weak Feller regularity also facilitates rigorous approximation and learning theoretic results, as shown in the paper.}




\maketitle

\section{Problem Description} \label{problem description}

Let there be $N$ agents in a decentralized stochastic control model. Let $\mathbb{X};\mathbb{U}^{1},\cdots,\mathbb{U}^N;\mathbb{Y}^{1},\cdots,\mathbb{Y}^N$ be standard Borel spaces, i.e., Borel subsets of metric Polish spaces (complete, and separable metric spaces); here $\mathbb{X}$ denotes the state space, the action spaces for agents $1,...,N$ are $\mathbb{U}^{1},...,\mathbb{U}^N$, and their measurement spaces are given by $\mathbb{Y}^{1},...,\mathbb{Y}^N$. We denote by $\mathcal{B}(\mathbb{X})$ the Borel $\sigma$-algebra on $\mathbb{X}$ and the same notation holds for all spaces. At every time $t \in \mathbb{Z}_+$ the system dynamics evolve as follows:
\begin{eqnarray} 
   \label{Prob 1} x_{t+1}&=&\Phi(x_{t},u^{1}_{t},...,u^{N}_{t},v_{t})    \\
   \label{obs 1} \forall i \in \{1,..,N\} \text{: } y^{i}_{t}&=&\phi^i(x_{t},w^{i}_{t}) 
\end{eqnarray}
where $\{ v_{t}\}_{t}$ and $\{ w^{i}_{t}\}_{t}$ are i.i.d noise processes for all $i \in \{1,..,N\}$, and $\Phi$ and $\phi^i$ are Borel measurable functions. We assume that all random variables are defined on a common probability space $(\Omega, {\cal F}, P)$. The dynamics and measurement equations above induce the following kernels for any Agent $i$: let $A\in {\cal B}(\mathbb{X})$ and $B \in {\cal B}(\mathbb{Y}^{i})$: 
\begin{align}
    \tau(A|x,u^{1},...,u^{N})&:=P(x_{t+1} \in A|x_t=x,u^1_t=u^{1},...,u^N_t=u^{N})\\
    Q^{i}(B|x)&:=P(y_{t} \in B|x_t=x)
\end{align}

At every time $t$, each Agent $i$ has access to their local information 
\begin{align}\label{localInf}
I_{t}^{i}=(y_{[0,t]}^{i}, u_{[0,t-1]}^{i})
\end{align}
as well as some common information $I_{t}^{C}$. We denote the spaces in which the random variables $I_{t}^{i}$ and $I_{t}^{C}$ take values by $\mathbb{I}_{t}^{i}$ and $\mathbb{I}_{t}^{C}$ respectively. Each Agent applies a policy $\gamma^{i}=\{\gamma^{i}_{t} \}_{t\in \mathbb{Z}^{+}}$, such that $\gamma_{t}^{i}:\mathbb{I}_{t}^{C}\times \mathbb{I}_{t}^{i}\rightarrow \mathbb{U}^{i}$ and is $\sigma(I_{t}^{i},I_{t}^{C})$-measurable, with the objective of minimizing 
\begin{equation} \label{cost}
   J(\gamma)=\sum_{t=0}^{\infty} \beta^{t} E^{\gamma}[c(x_{t},\mathbf{u_{t}})] 
\end{equation}
where $\beta \in (0,1)$ is the discount factor, $\gamma=(\gamma^{1},...,\gamma^{N})$ is the team policy, and $c: \mathbb{X} \times \mathbb{U} \rightarrow [0,\infty)$ is a cost functional which is assumed to be continuous and bounded. We refer to the set of all admissible team policies as $\Gamma$. We also assume that $x_{0}\sim \mu$ and that the initial distribution $\mu$ is known to all agents at time $t=0$. Next, we will describe the information structures that will be considered in this paper. 

\subsection{Various information structures} \label{information structures}
We define the various information structures that we will consider in this paper. In all cases, each agent has access to their own measurements and previous actions (\ref{localInf}). However, the information shared among agents, i.e., the common information will vary.
\begin{itemize}
    \item[(i)] \textit{One-step delayed information sharing pattern (OSDISP)}: \[I_{t}^{C}=\{y_{[0,t-1]}^{[1,N]},u_{[0,t-1]}^{[1,N]}\}.\]
    \item[(ii)] $K$\textit{-step periodic information sharing pattern (KSPISP)}: Let $t=qK+r$ where $q$ and $r$ are non-negative integers such that $r\leq K-1$ and $K\in \mathbb{N}$ represents the period between successive information sharing between the Agents. Here the common information at time $t$ is given by $I^{C}_{t}=\{y^{[1,N]}_{[0,qK-1]},u^{[1,N]}_{[0,qK-1]}\}$.
    \item[(iii)] \textit{Completely decentralized information structure (CDIS)}: $I_{t}^{C}=\{\emptyset\}$ 
\end{itemize}
\begin{remark}
    Note that when $K=1$, the KSPISP problem reduces to the OSDISP problem. We note that in this paper we will not address the CDIS directly, however, the results which will be presented for the KSPISP can be used to obtain near optimal solutions to CDIS problems. 
\end{remark}

\section{Literature Review and Contributions}

\subsection{Literature review}
Contributions of Witsenhausen to decentralized stochastic control theory have been foundational, as reviewed in \cite{wit71,YukselBasarBook24}, and it is by now well known that decentralized stochastic control problems are intrinsically difficult \cite{Witsenhausen1968counterexample} even in simple linear Gaussian models despite the abundance of application areas. In particular, the inapplicability of standard methods of single-agent stochastic control makes arriving at numerical solutions untenable even when the dynamics of a given problem are relatively simple \cite{bernstein2022comp-maor,wit71,Witsenhausen1968counterexample,YukselBasarBook24}. 

Dynamic programming has nonetheless been a useful method for decentralized stochastic control since the beginning of the 1970s. A possible approach toward establishing a tractable dynamic program is through the construction of a controlled Markov chain where the controlled Markov state may now live in a larger state space and the actions are elements in possibly function spaces. Such a {\it dynamic programming approach} has been adopted extensively in the literature for linear models and models with finite spaces (see for example, \cite{Athans,VaraiyaWalrand,WitsenhausenSeparation}, \cite{yos75}, \cite{ChongAthans}, \cite{AicardiDavoli}, \cite{YukTAC09}, \cite{lamperski2015optimal}). This was significantly generalized and termed as the {\it common information approach} in \cite{NayyarMahajanTeneketzis} and \cite{NayyarBookChapter}) through the use of a team-policy which uses common information to generate prescriptions (local policies) for each DM to generate its actions using local information. For finite models, another related method utilizing common knowledge among agents, is presented in \cite{nayyar2019common}. In \cite{tavafoghi2018unified}, it is shown that that when a certain compression function exists for the private information, optimal policies which rely on such a compression are near optimal. Similarly, this result is extended in \cite{Subramanian-2021-common} to the compression of the common information. Recently \cite{Saldi2023Commoninformationapproachstaticteamproblems} studied the standard Borel setup for the case of static team problems. 


The construction above requires a common knowledge among decision makers, which is a unifying assumption in the aforementioned contributions in the literature. Alternatively, references \cite{WitsenStandard} and \cite{YukselWitsenStandardArXiv} developed general dynamic programming algorithms for all sequential teams in standard Borel spaces (teams where agents act in a predetermined order), which are mathematically consequential for existence and optimality analysis. 

Also for finite spaces, in the computer science literature, for completely decentralized information structure problems a related approach was presented in \cite{dibangoye2016} where it is shown that such problems are equivalent to a deterministic control problem where local policies are viewed as actions and the state consists of the window of all local policies used up to a given time stage. Moreover, it is shown \cite{dibangoye2016} that an equivalent MDP can be constructed whereby the new state, called \textit{occupancy state}, consists of the conditional probability on the state as well as the history process of each agent given the window of local policies. Additionally, representations of the aforementioned MDP which lead to dimension reduction are discussed \cite{dibangoye2016}. This discussion also motivates a study on robustness and continuity to approximations and learning, which is executed in our paper.

A further primary contribution is on existence of optimal policies in decentralized stochastic control. Present existence results, in their most general forms to our knowledge as reviewed in \cite[Theorems 4.5-4.7]{saldiyukselGeoInfoStructure} require either applicability of change of measure arguments (under absolute continuity conditions of the measurements given past actions) or nested information structures, which are not applicable to the setups we present in this paper. Additionally, the studies of learning for decentralized stochastic control problems with limited or no information sharing have typically focused on the development of algorithms that converge but with no performance guarantees (see \cite{mao2023decentralized,gupta2017cooperative,hu2019simplified,wu2013monte,liu2015stick,kraemer2016multi} for example) in part because the necessary results on structural and rigorous approximations had not been available. By developing structural and regularity results for decentralized control problems, we are able to provide conditions under which dynamic programming can be used to arrive at solutions analytically. Further, we will illustrate how those structural and regularity results can be used in unison with the results of \cite{SaLiYuSpringer} to obtain near optimal numerical solutions. To arrive at our results, we show that both the one-step delayed information sharing pattern and $K$-step periodic information sharing pattern problems can be reduced to a centralized problem. We then use our results to find conditions under which the centralized reduction of the OSDISP problem and KSPISP problem is weak-Feller. The weak-Feller property is important as it enables one to establish existence of optimal solutions and use approximations \cite{SaYuLi15c,SaLiYuSpringer}, and Q-learning \cite{KSYContQLearning,kara2021convergence}.


Our analysis involves discrete-time decentralized stochastic control models. Research in continuous-time under decentralized information structures has been relatively recent beyond the perfect-recall setting \cite{rishel1970necessary,jacka2015informational,elliott1977optimal} or the mean-field setting \cite{caines2018peter}. See \cite{charalambous2016decentralized} and \cite{charalambous2017centralizedI} for a class of convex problems, and \cite{huang2022general,jackson2023approximately} for recent studies on continuous-time stochastic team problems; \cite{jackson2023approximately} studies near optimality of decentralized policies for large models. Under absolute continuity conditions, \cite{pradhanyuksel2023DTApprx} obtains existence results and develops a discrete-time approximation theory under general information structures. 

\subsection{Main contributions} \label{main contributions}
\begin{itemize}
\item[(1)] {\bf On the one-step delayed information sharing pattern.} (i)  We show that under the one-step delayed information sharing pattern (see Theorem \ref{centralized reduction one step part one} and \ref{them2young}) optimal policies can be expressed in a separated form and the problem reduces to an equivalent centralized Markov Decision Process (MDP); this generalizes the prior results which considered finite models as well as a recent result  \cite{Saldi2023Commoninformationapproachstaticteamproblems} which considers static Borel models. A key measurability argument that we prove is critical for our result (see Theorem \ref{centralized reduction one step part one}) which shows that infima over control policies under the original information structure and the reduced form  are equivalent. (ii) We then use this result to show that under mild conditions on the original kernel and measurement channels, the new centralized MDP is weak-Feller (Theorem \ref{theorem5}).  (iii) We use the weak-Feller property to establish the existence of optimal solutions as well as stationarity (for the infinite horizon discounted setup) of optimal solutions (see Corollary \ref{existence OSDISP}).

\item[(2)]{\bf On the $K$-step periodic information sharing pattern.}
    (i) By expanding on the results for the one-step delayed information sharing pattern, we show that under the $K$-step periodic information sharing pattern with standard Borel spaces (see Theorem \ref{themkstepMDP} and \ref{centralized reduction kstep step part two}), policies can be taken to be of separated form without any loss, and the problem reduces to an equivalent centralized MDP. (ii) We use this result to show that under some conditions on the original kernel and the observation channels the centralized MDP is weak-Feller (Theorem \ref{theorem6}). (iii) Moreover, we use the weak-Feller property to establish the existence of optimal solutions as well as stationarity of optimal solutions (see Corollary \ref{existence KSDISP}).  

    We note that the one-step delayed sharing pattern problem is a special case of the $K$-step periodic sharing pattern problem with $K=1$. However, establishing many of the important properties of the $K$-step periodic information sharing pattern problem, such as the weak-Feller regularity, requires that one first establishes those properties for the one-step delayed information sharing pattern problem. Moreover, one can establish the weak-Feller property under weaker assumptions for the OSDISP problem than those needed for the KSPISP problem. 
    
\item[(3)]{\bf On existence and approximations in decentralized stochastic control.} The problem of existence of optimal policies in decentralized stochastic control has been studied in several publications \cite{WuVer11}, \cite{wit68} \cite{gupta2014existence}, \cite{YukselWitsenStandardArXiv} \cite{SaldiArXiv2017}\cite{charalambous2017centralizedI}, however these do not apply to the setup presented in our paper, due to a lack of static reducibility, non-nested nature of information patterns, or generality of the models. We thus provide further existence and approximation results in decentralized stochastic control.   
\end{itemize}

An important motivation for our work is that CDIS problem can be approximated with KSPISP problems as $K$ becomes sufficiently large. Thus the results of this paper which allow us to arrive at analytical and numerical solutions to KSPISP problems will have significant implications for CDIS problems.

\subsection{Notation and Preliminaries}

\subsubsection{Notation}

   (i) For any standard Borel set $\mathbb{X}$ (Borel measurable subset of a complete and separable metric space) we will use the notation $\mathcal{B}(\mathbb{X})$ to denote the Borel sigma-algebra on $\mathbb{X}$. (ii) For any standard Borel set $\mathbb{X}$ we denote the space of probability measure on $\mathbb{X}$ under the weak convergence topology by $\mathcal{P}(\mathbb{X})$. (iii) We will denote a finite sequence $\{x_{t}\}_{t=l}^{t=k}$, where $k,l\in \mathbb{Z}^{+}$ by $x_{[l,k]}$ or $x^{[l,k]}$. (iv)  Similarly, for a double-indexed sequence $\{x_{t}^{i}\}_{\{t\in\{l,...,k\},i\in\{l',...,k'\}\}}$ we will use the notation $x^{[l',k']}_{[l,k]}$. (v) Whenever the superscript of a double-indexed sequence ranges over all agents ($i\in \{1,...,N\}$) and when there is no ambiguity, we will use the notation $\mathbf{x}_{[l,k]}$ to refer to the double-indexed sequence $\{x_{t}^{i}\}_{\{t\in\{l,...,k\},i\in\{1,...,N\}\}}$. (vi) Moreover, we will use the notation $\mathbb{U}=\mathbb{U}^{1}\times...\times\mathbb{U}^{N}$ and $\mathbb{Y}=\mathbb{Y}^{1}\times...\times\mathbb{Y}^{N}$. (vii) We will use $C_{b}(\mathbb{X})$ to denote the set of continuous and bounded functions on the space $\mathbb{X}$. (viii)  For a measurable function $f:(\mathbb{X},\mathcal{F}_{X})\rightarrow (\mathbb{U},\mathcal{B}(\mathbb{U}))$ where $\mathcal{F}_{X}$ is a $\sigma$-algebra on $\mathbb{X}$ and $\mathbb{U}$ is standard Borel, we will use the notation 
   \[ f(du|x): x \mapsto f(du|x) \in \mathcal{P}(\mathbb{U}),\]
   to refer to the conditional probability induced by $f$ on $\mathbb{U}$ so that $f(.|x)=\delta_{f(x)}(.)$ where $\delta_{f(x)}(.)$ is the Dirac measure concentrated at $f(x)$.\\
   
   Throughout this paper all probability valued spaces are endowed with the weak topology.

\subsubsection{Convergence of probability measures}
In this section, we will introduce three notions of convergence for sequences of probability measures that will be useful later on: weak convergence, convergence under total variation, and convergence under the Wasserstein metric of order one. For a complete, separable and metric space $\mathbb{X}$, a sequence $\left\{\mu_n \right\}_{n \in \mathbb{N}}$ $\subset$ $\mathcal{P}(\mathbb{X})$ is said to converge to $\mu \in \mathcal{P}(\mathbb{X})$ weakly if and only if $\int_{\mathbb{X}} f(x) \mu_n(d x) \rightarrow \int_{\mathbb{X}} f(x) \mu(d x)$ for every continuous and bounded $f: \mathbb{X} \rightarrow \mathbb{R}$. One important property of weak convergence is that the space of probability measures on a complete, separable, metric (Polish) space endowed with the topology of weak convergence is itself complete, separable, and metric. An example of such a metric is the bounded Lipschitz metric, which is defined for $\mu, \nu \in \mathcal{P}(\mathbb{X})$ as
$$
\rho_{B L}(\mu, \nu):=\sup _{\|f \|_{BL\leq 1}}\left|\int f d \mu-\int f d \nu\right|
$$
where
$$
\|f\|_{B L}:=\|f\|_{\infty}+\sup _{x \neq y} \frac{|f(x)-f(y)|}{d(x, y)}
$$
and $\displaystyle \|f\|_{\infty}=\sup _{x \in \mathbb{X}}|f(x)|$. Here, $d(.,.)$ denotes the metric on $\mathbb{X}$. For probability measures $\mu, \nu \in \mathcal{P}(\mathbb{X})$, the total variation metric is given by
\begin{eqnarray*}
    &&\|\mu-\nu\|_{T V}=2 \sup _{B \in B(\mathbb{X})}|\mu(B)-\nu(B)|=\sup _{\|f\|_{\infty \leq 1}}\left|\int f(x) \mu(\mathrm{d} x)-\int f(x) \nu(\mathrm{d} x)\right|,
\end{eqnarray*}

\noindent where the supremum is taken over all measurable real-valued functions $f$ such that $\displaystyle \|f\|_{\infty}=\sup_{x \in \mathbb{X}}|f(x)| \leq$ 1. A sequence $\mu_{\mathrm{n}}$ is said to converge in total variation to $\mu \in \mathcal{P}(\mathbb{X})$ iff $\| \mu_n-\mu\|_{T V} \rightarrow 0$.\\
The Wasserstein order one metric for measures with bounded support is given by
\begin{equation*}
    W_{1}(\mu,\nu)=\sup_{\|f\|_{\text{Lip}}\leq 1}\bigg| \int f d\mu - \int f d\nu   \bigg| \text{, where } \|f\|_{\text{Lip}}:=\sup _{x \neq y} \frac{|f(x)-f(y)|}{d(x, y)}
\end{equation*}

\section {Equivalent Formulations via Centralized MDP Reductions} \label{equivalent formulation}

\subsection{One-step delayed information sharing pattern as a centralized MDP}
It was shown in \cite{yos75} and \cite{VaraiyaWalrand} that optimal policies for the one step-delayed information sharing pattern, under finite horizon cost criteria, are of separated form in the case where the measurement spaces, action spaces, and state spaces are all finite. In both \cite{yos75} and \cite{VaraiyaWalrand}, separability of policies was shown via dynamic programming. This approach requires a continuity, and in this case a weak-Feller continuity, assumption in order for measurable selection to be applicable-thus allowing for the use of dynamic programming [Chapter 3 \cite{HernandezLermaMCP}]. This assumption always holds in the case of finite state and action spaces. In this section, we will generalize the separability result by presenting an  approach that is suitable for problems with standard Borel spaces and infinite horizon discounted cost criteria.

First, we will show that the OSDISP problem admits a centralized MDP reduction where the one-step predictor acts as the state. 

\noindent{\bf Young topology on Maps as Actions.} Here, we will introduce the  Young topology on control policies (see, e.g. \cite{BorkarRealization,yuksel2023borkar}). Suppose there exists a reference measure $\psi\in \mathcal{P}(\mathbb{Y}^{i})$ and a measurable function $\displaystyle h^{i}: \mathbb{X} \times  \mathbb{Y}^{i}\rightarrow [0,\infty)$ such that the measurement channel, of each individual agent $i$, $Q^{i}$ satisfies:
\begin{equation} \label{Young top}
    Q^{i}(B|x_{t})=\int_{B} \psi(dy^{i}_{t})h^{i}(x_{t},y^{i}_{t}) \text{ for all } B\in\mathcal{B}(\mathbb{Y}^{i})
\end{equation}
For the sake of readability, throughout the rest of the paper, we will assume that $h^{i}\equiv h$ for all $i \in \{1,...,N\} $ where $h$ is a measurable function such that $\displaystyle h: \mathbb{X} \times  \mathbb{Y}^{i}\rightarrow \mathbb{R}$. All the proofs for the more general case proceed in the same manner as the ones included in this paper.

\begin{definition} \label{topology for actions onestep} For all $i \in \{1,...,N\}$ let $\{f^{i}_{n,t}\}$ be a sequence of measurable functions such that $f^{i}_{n,t}:\mathbb{Y}^{i} \rightarrow \mathcal{P}(\mathbb{U}^{i})$ and let $f^{i}_{t}:\mathbb{Y}^{i} \rightarrow \mathcal{P}(\mathbb{U}^{i})$ be a measurable function. We say
$f_{n,t}=(f^{1}_{n,t},...,f^{N}_{n,t})\rightarrow f^{i}_{t}=(f^{1}_{t},...,f^{N}_{t})$ under the Young topology if and only if for all $i \in \{1,...,N\}$ and for all $g \in C_{b}(\mathbb{Y}^{i},\mathbb{U}^{i})$ : 
\begin{eqnarray*}
  \int g(y^{i},u^{i})f^{i}_{n,t}(du^{i}|y^{i})\psi(dy^{i})\rightarrow \int g(y^{i},u^{i})f^{i}_{t}(du^{i}|y^{i})\psi(dy^{i}).  
\end{eqnarray*}

\end{definition}
\begin{remark}
    It follows that, since the marginal on the measurement space is fixed, the convergence in Definition \ref{topology for actions onestep} is equivalent to the requirement that $\int g(y^{i},u^{i})f^{i}_{n,t}(du^{i}|y^{i})\psi(dy^{i})\rightarrow \int g(y^{i},u^{i})f^{i}_{t}(du^{i}|y^{i})\psi(dy^{i})$ for any measurable function $g$ which is bounded and continuous in $y^{i}$ (\cite[Theorem 3.10]{Schal},  \cite[Theorem 2.5]{balder2001}).
\end{remark}

This topology leads to a convex and compact formulation (see Section 2 in\cite{BorkarRealization} as well as Section 3 in \cite{Bor91}). In particular, consider the set of probability measures on $\mathbb{Y}^{i} \times \mathbb{U}^{i}$ with fixed marginal $\psi\in \mathcal{P}(\mathbb{Y}^{i})$: 
\begin{eqnarray*}
    E^{i}&=&\big\{P\in \mathcal{P}(\mathbb{Y}^{i} \times \mathbb{U}^{i}) \mid P(B)=\int_{B} f(du^{i}|y^{i})\psi(dy^{i}),\text{ for some } f\in \mathcal{F}^{i} \text{ and for all } B\in\mathcal{B}(\mathbb{Y}^{i} \times \mathbb{U}^{i})  \big \}
\end{eqnarray*}
where $\mathcal{F}^{i}$ is the space of regular conditional probability measures from $\mathbb{Y}^{i}$ to $\mathbb{U}^{i}$ (see \cite[section 10.2]{Dud02} or \cite[Appendix 4]{dynkin1979controlled1}): 
\begin{eqnarray*}
    \mathcal{F}^{i}=\big\{f:(\mathbb{Y}^{i},\mathcal{B}(\mathbb{Y}^{i})) \rightarrow (\mathcal{P}(\mathbb{U}^{i}),\mathcal{B}(\mathcal{P}(\mathbb{U}^{i})))| \text{ }f \text{ is measurable}  \big \}
\end{eqnarray*}
Recall that $\mathcal{P}(\mathbb{U}^{i})$ is endowed with the weak topology. Then, every deterministic map $f^{i}: \mathbb{Y}^{i} \rightarrow \mathbb{U}^{i}$ can be identified as an element of the set of extreme points of $E^{i}$ \cite[Lemma 2.2]{BorkarRealization}:
\begin{eqnarray*}
   \nonumber \Theta^{i}&=&\big\{P\in \mathcal{P}(\mathbb{Y}^{i} \times \mathbb{U}^{i}) \mid P(B)=\int \mathbbm{1}_{\{(y^{i},g(y^{i}))\in B\}}\psi(dy^{i}),\text{ for some measurable}\\
   &&\text{function } g:\mathbb{Y}^{i}\rightarrow \mathbb{U}^{i}\text{ and for all } B\in\mathcal{B}(\mathbb{Y}^{i} \times \mathbb{U}^{i})  \big \}
\end{eqnarray*}
Thus, $\Theta^{i}$ inherits the Borel measurability and topological properties
of the Borel measurable set $E^{i}$ \cite[Section 2]{BorkarRealization}.

The following is a key result.

\begin{theorem} \label{centralized reduction one step part one}
    The problem (\ref{Prob 1}),(\ref{obs 1}),(\ref{cost}) is equivalent to one where the team policy at time $t$ is given by $(\Tilde{\gamma}_{t}^{1},...,\Tilde{\gamma}_{t}^{N})$ such that for all $i$, $\Tilde{\gamma}^{i}_{t}:I_{t}^{C} \mapsto f^{i}_{t}$ is $\sigma(I_{t}^{C})$-measurable and $f^{i}_{t}:y^{i}_{t} \mapsto u_{t}^{i}$ is $\sigma(y^{i}_{t})$-measurable and represents the action of Agent $i$ at time $t$. Here, by equivalent, we mean that $\displaystyle\inf_{\gamma\in\Gamma}J(\gamma)=\inf_{\Tilde{\gamma}\in\Tilde{\Gamma}}J(\Tilde{\gamma})$ where $\Gamma$ is the space of all team policies such that at time $t$ Agent $i$ applies a measurable policy $\gamma^{i}_{t}:I^{C}_{t}\times I^{i}_{t}\mapsto u^{i}_{t}$ and $\Tilde{\Gamma}$ is the space of all policies of the form $\Tilde{\gamma}=\{(\tilde{\gamma}^{1}_{t},...,\tilde{\gamma}^{N}_{t})\}_{t\in \mathbb{Z}^{+}}$ where $\tilde{\gamma}^{i}_{t}$ is as defined above.
\end{theorem}

\noindent\textbf{Proof.} Before, we proceed with the proof of Theorem \ref{centralized reduction one step part one} we recall a result of \cite{DubinsFreedman} that will be useful.
\begin{theorem}  \cite[Theorem 2.1]{DubinsFreedman} \label{DubinsFreedman- theorem}
    Consider two measurable spaces $(\Omega,\mathcal{F})$, and $(\mathbb{X}, \Sigma)$ where $\Sigma$ is countably generated. Let $\mathbb{M}$ denote the set of all signed measures on $\mathbb{X}$. Define $\Sigma^{*}$ as the smallest sigma algebra of subsets of $\mathbb{M}$ such that for all $A\in \Sigma$ the function $\zeta$ such that $\zeta: \begin{matrix}
        (\mathbb{M},\Sigma^{*})\rightarrow (\mathbb{R},\mathcal{B}(\mathbb{R}))\\
        \mu \mapsto \mu(A)
    \end{matrix}$ 
     is measurable. Then, a function $\phi$ from $\Omega$ to $\mathbb{M}$ is measurable if and only if the function $\omega \mapsto \phi(\omega)(A)$ is measurable from $(\Omega,\mathcal{F})$ to $(\mathbb{R},\mathcal{B}(\mathbb{R}))$ for all $A \in \Sigma$. 
\end{theorem}
Let $\{(\gamma_{t}^{1},...,\gamma_{t}^{N})\}_{t=0}^{\infty}$ be a team policy as defined in Section \ref{problem description}.

\noindent{\bf Step 1}: We have, by iterated expectations, that

\begin{align*}
       &\sum_{t=0}^{\infty} \beta^{t} E\big[c(x_{t},\mathbf{u_{t}})\big]= \sum_{t=0}^{\infty} \beta^{t} E\big[E[c(x_{t},\mathbf{{u}_{t}})| I^{C}_{t}]\big] \\
   &= \sum_{t=0}^{\infty} \beta^{t} E\big[\int P(dx_{t}|I^{C}_{t})(\prod_{i=1}^{N}Q^{i}(dy^{i}_{t}|x_{t}))c(x_{t},\gamma^{1}_{t}(I^{C}_{t},y^{1}_{t}),...,\gamma^{N}_{t}(I^{C}_{t},y^{N}_{t}))\big]
\end{align*}


\noindent{\bf Step 2}: Next, to complete the proof, we will show that for any policy $\gamma^{i}_{t}:\mathbb{I}^{C}_{t}\times \mathbb{Y}^{i}\rightarrow \mathbb{U}^{i}$ there exists a policy $\Tilde{\gamma}^{i}_{t}:\mathbb{I}^{C}_{t}\rightarrow \Theta^{i}$ such that $\Tilde{\gamma}^{i}_{t}:I^{C}_{t}\mapsto f\in \Theta^{i}$ and leads to the same cost. Conversely, for every policy $\Tilde{\gamma}^{i}_{t}$ there exists a policy $\gamma^{i}_{t}$ that achieves the same cost. To do this, it is sufficient to check that the induced {\it strategic measures} (that is, probability measures on the sequence spaces of actions and measurements) are equivalent almost surely: $P^{\gamma^{i}_{t}}(du^{i}_{t},dy^{i},dI^{C}_{t})=P^{\Tilde{\gamma}^{i}_{t}}(du^{i}_{t},dy^{i},dI^{C}_{t})$ a.s. To prove this, it is sufficient to check that for a countable family of measure-determining continuous and bounded functions $g\in \mathbf{C}_{b}(\mathbb{U}^{i}\times\mathbb{Y}^{i}\times\mathbb{I}^{C}_{t})$ \cite[Theorem 3.4.5]{ethier2009markov} (note that without loss of generality these functions can be taken so that $\|g\|_{\infty}\leq 1$) the following holds
\begin{align} \label{check1}
    \int g(u^{i}_{t},y^{i}_{t},I^{C}_{t})P^{\gamma^{i}_{t}}(du^{i}_{t},dy^{i},dI^{C}_{t})=\int g(u^{i}_{t},y^{i}_{t},I^{C}_{t})P^{\Tilde{\gamma}^{i}_{t}}(du^{i}_{t},dy^{i},dI^{C}_{t})
\end{align}
Note that whenever $\gamma^{i}_{t}$ is a kernel, it can be realized \cite{Bor93} as $\gamma^{i}_{t}(A\mid y^{i},I^{C}_{t})=\int \mathbbm{1}_{\{\Delta(y^{i},I^{C}_{t},\Bar{w})\in A\}}P(d\Bar{w})$ where $\Delta$ is a measurable $\mathbb{U}^{i}$-valued function and $\Bar{w}$ is some $[0,1]-$valued random variable which is independent of $y^{i}$ and $I^{C}_{t}$. \\
Consider $\Delta(.,I^{C}_{t},\Bar{w}):(I^{C}_{t},\Bar{w})\mapsto \theta \in \Theta^{i}$, by \cite[Theorem 8.5]{Rud87} $\Delta(.,I^{C}_{t},\Bar{w})$ is measurable. Define $\Tilde{\gamma}^{i}_{t}$ as a map from $\mathbb{I}^{C}_{t}$ to $\mathcal{P}(\mathbb{Y}^{i}\times \mathbb{U}^{i})$ which is realized by the measurable function $\Delta(.,I^{C}_{t},\Bar{w})$. First, observe that for every $A\in\mathcal{B}(\mathbb{Y}^{i}\times \mathbb{U}^{i})$ $\Tilde{\gamma}^{i}_{t}(I^{C}_{t})(A)=\int \mathbbm{1}_{\{(y,\Delta(y,I^{C}_{t},w))\in A\}}\psi(dy^{i})P(dw)$ which is $\mathcal{B}(\mathbbm{I}^{C}_{t})/\mathcal{B}(\mathbb{R})$ measurable. Hence, by Theorem \ref{DubinsFreedman- theorem} $\Tilde{\gamma}^{i}_{t}$ is measurable. It then follows that $\displaystyle\inf_{\gamma\in\Gamma}J(\gamma)\geq\inf_{\Tilde{\gamma}\in\Tilde{\Gamma}}J(\Tilde{\gamma})$.

\noindent{\bf Step 3}: For the reverse direction, suppose that $\Tilde{\gamma}^{i}_{t}$ is a transition kernel which is realized by some measurable function $D(I^{C}_{t},\Bar{w}):(I^{C}_{t},\Bar{w})\mapsto \theta \in \Theta^{i}$. Now, we will present two different arguments each of which can be used to complete the proof.

First, we have that 

\begin{align*}
&P^{\tilde{\gamma}^{i}_{t}}(du^{i}, dy^{i} | I^{C}_t):=\int_{\Theta^{i}\times\mathbb{U}^{i}\times\mathbb{Y}^{i}} P(du^i,dy^i,df^i | I^{C}_{t})=\int_{\Theta^{i}\times\mathbb{U}^{i}\times\mathbb{Y}^{i}} P(du^{i} | y^{i},f^{i},I^{C}_{t})P(dy^i,df^i | I^{C}_{t})\\
&=\int_{\Theta^{i}\times\mathbb{U}^{i}\times\mathbb{Y}^{i}}P(du^{i} | y^{i},f^{i}) P(dy^{i},df^{i} | I^{C}_{t})=\int_{\Theta^{i}\times\mathbb{U}^{i}\times\mathbb{Y}^{i}}P(du^{i} | y^{i},f^{i}) P(dy^{i} | I^{C}_{t}) P(df^{i} | I^{C}_{t})\\
&=\int_{\Theta^{i}\times\mathbb{U}^{i}\times\mathbb{Y}^{i}}f^{i}(du^{i}|y^{i})P(dy^{i} | I^{C}_{t})\tilde{\gamma}^{i}_{t}(df^{i}|I^{C}_{t})
\end{align*}

Hence, $P^{\Tilde{\gamma}^{i}_{t}}(A |I^{C}_{t})$ is a well defined conditional probability measure. In particular, one can write the measure induced by $\Tilde{\gamma}^{i}_{t}$ on the measurement and action of Agent $i$ as 
\begin{align*}
P^{\Tilde{\gamma}^{i}_{t}}((u^{i},y^{i})\in . |I^{C}_{t})&=\int_{\cdot} P^{\Tilde{\gamma}^{i}_{t}}(du^{i}|y^{i}, I^{C}_{t})P(dy^{i}|I^{C}_{t})=\int_{\cdot} P^{\Tilde{\gamma}^{i}_{t}}(du^{i}|y^{i}, I^{C}_{t})Q^{i}(dy^{i}_{t}|x_{t})P(dx_{t}|I^{C}_{t})
\end{align*}
Then, by letting $\gamma^{i}_{t}(du^{i}|y^{i}, I^{C}_{t})=P^{\Tilde{\gamma}^{i}_{t}}(du^{i}|y^{i}, I^{C}_{t})$, we can see that both $\Tilde{\gamma}^{i}_{t}$ and $\gamma^{i}_{t}$ induce the same measure on $\mathbb{U}^{i}$. Hence, we have shown that for every policy  $\tilde{\gamma}^{i}_{t}\in \tilde{\Gamma}$ there exists a policy $\gamma^{i}_{t}\in \Gamma$ such that equation (\ref{check1}) holds. Thus,
$\displaystyle\inf_{\gamma\in\Gamma}J(\gamma)\leq\inf_{\Tilde{\gamma}\in\Tilde{\Gamma}}J(\Tilde{\gamma})$.

We can also present a different argument which holds given that the cost function $c$ is bounded. To this end, we recall a version of \cite[Theorem 2]{Kohur1972measurabilityoftwovariablefunction} that will be useful. A version of this theorem for complex-valued functions was proved in \cite[Lemma 9.2]{Mackey1952inducedrepresentationsoflocallycompactgroups}.
\begin{theorem} \cite[Theorem 2]{Kohur1972measurabilityoftwovariablefunction} \label{product-meas-thm}
    Let $(\mathbb{X},\mathcal{T})$ be a measurable space and $\mathbb{Y}$ be a Polish metric space, i.e., a complete and separable metric space. Let $f:\mathbb{X}\times\mathbb{Y}\rightarrow \mathbb{A}$ where $\mathbb{A}$ is a Polish space. Suppose $f$ is $\mathcal{T}$-measurable for every $y\in\mathbb{Y}$ and continuous in $y$ for every $x\in \mathbb{X}$. Then, $f$ is $\mathcal{T} \times \mathcal{B}(\mathbb{Y})$ measurable.
\end{theorem}
\begin{remark}
    As noted in \cite{Kohur1972measurabilityoftwovariablefunction} continuity is crucial in order for measurability to hold. See \cite{Deuxexemplesdefonctionsnonmesurables} for counterexamples.
\end{remark}Let $\gamma^{i}_{t}$ be a conditional probability measure which is realized by the function $\Delta(y^{i},I^{C}_{t},\Bar{w})=D(I^{C}_{t},\Bar{w})(y^{i}_{t})$. Because the image of $D$ consists entirely of measurable $f \in \Theta^{i}$, it follows that $\Delta$ is measurable in its first component. By Lusin's theorem, for any arbitrary $\epsilon>0$ there exists a set $K\in \mathcal{B}(\mathbb{I}^{C}_{t}\times \Omega)$, where $\Omega$ is the space for $\Bar{w}$, such that the restriction of $D$ to $K$ is continuous and $P(K^{C})<\epsilon$ where $K^{C}$ denotes the complement of $K$ and $P$ is the marginal probability measure on $\mathcal{B}(\mathbb{I}^{C}_{t}\times \Omega)$. Hence, we have that the restriction of $\Delta$ to the set $K\times \mathbb{Y}^{i}$ is continuous in $I^{C}_{t},\Bar{w}$ and measurable in its first component, it then follows by Theorem \ref{product-meas-thm} that this restriction is measurable. By defining $\Delta$ to be any arbitrary value in $\mathbb{U}^{i}$ outside of $K\times \mathbb{Y}^{i}$ we get that $\Delta$ is measurable. Hence, because $g$ in equation \ref{check1} is bounded , we have shown that for every policy $\Tilde{\gamma}$ one can construct a policy $\gamma$ that achieves an arbitrarily close cost at time $t$ to the one incurred if the policy  $\Tilde{\gamma}$ was used instead. It is easy to see that for any finite horizon cost criteria and for any policy of the form $\tilde{\gamma}$ there is a policy of the form $\gamma$ that achieves an arbitrarily close cost. Because $c$ is bounded, one can approximate, uniformly over all policies, the infinite cost criteria using finite cost criteria. Thus, it follows that $\displaystyle\inf_{\gamma\in\Gamma}J(\gamma)\leq\inf_{\Tilde{\gamma}\in\Tilde{\Gamma}}J(\Tilde{\gamma})$.  \qed

As noted earlier, for the case of static teams, a related equivalence has been established in \cite{Saldi2023Commoninformationapproachstaticteamproblems}. In \cite{Saldi2023Commoninformationapproachstaticteamproblems}, the argument relies on first viewing policies as maps from common and private information to the space of signed measures on the action spaces endowed with the total variation norm. Then, the space of such kernels is identified as a subset of the set of all equivalence classes of $w^{*}$-measurable functions viewed as an appropriate dual space of continuous space valued maps. By contrast, we endow the policies with the Young topology which, in turn, leads to a more straightforward proof and allows for more direct reasoning. 

Note that, here, the cost criteria can be expressed as
\begin{align*}
   &\sum_{t=0}^{\infty} \beta^{t} E\big[c(x_{t},\mathbf{u_{t}})\big]=\sum_{t=0}^{\infty} \beta^{t} E\big[\int P(dx_{t}|I^{C}_{t})(\prod_{i=1}^{N}Q^{i}(dy^{i}_{t}|x_{t}))c(x_{t},\gamma^{1}_{t}(I^{C}_{t},y^{1}_{t}),...,\gamma^{N}_{t}(I^{C}_{t},y^{N}_{t}))\big]\\   
   &= \sum_{t=0}^{\infty} \beta^{t} E\big[\int P(dx_{t}|I^{C}_{t})(\prod_{i=1}^{N}Q^{i}(dy^{i}_{t}|x_{t}))c(x_{t},\Tilde{\gamma}^{1}_{t}(I^{\small C}_{t})(y^{1}_{t}),...,\Tilde{\gamma}^{N}_{t}(I^{\small C}_{t})(y^{N}_{t}))\big]
\end{align*}

\begin{theorem}\label{them2young}
 Define the predictor at time $t$ as $Z_{t}(\cdot):=P(x_{t} \in \cdot |y_{[0,t-1]}^{[1,N]},u_{[0,t-1]}^{[1,N]})$. Then $( Z_{t} ,f_{t})$, where $f_{t}=(f^{1}_{t},...,f^{N}_{t})$, forms a fully observed MDP and is equivalent to the problem given by (\ref{Prob 1}, \ref{obs 1}, \ref{cost}) under an appropriate choice of the cost functional that is to say the choice of a measurable function $\tilde{c}:(Z_{t},f_{t})\mapsto [0,\infty)$. Here, the team policy at time $t$ consists of $\gamma_{t} \in \{(\gamma^{1}_{t},...\gamma^{N}_{t}) | \forall i \in {1,...,N}:  \gamma^{i}_{t}: Z_{t} \mapsto f^{i}_{t}$ and $\gamma^{i}_{t}$ is $\sigma(Z_{t})$-measurable$\}$.  
 \end{theorem}

\begin{remark} While this theorem holds for deterministic functions $(f^{1}_{t},...,f^{N}_{t})$ we allow for randomized maps $(f^{1}_{t},...,f^{N}_{t})$. This is because our analysis later on will require that the action space be compact and hence closed. The latter can not be guaranteed if we consider deterministic functions under Young topology given by Definition \ref{topology for actions onestep}.
\end{remark}

\textbf{Proof of Theorem \ref{them2young}} To prove the theorem, we first find an update equation for the predictor on the state given the common information: $Z_{t}=P(dx_{t}| I_{t}^{C})$ then we use that to find the stochastic transition kernel for  $Z_{t}$. Subsequently, we establish that $( Z_{t} ,f^{1}_{t},...,f^{N}_{t})$ is a controlled Markov Chain. We then rewrite the cost in \eqref{cost} in the terms of $Z_{t}$. We, then, conclude that the resulting problem with $Z_{t} $ as the state is equivalent to problem (\ref{Prob 1}, \ref{obs 1}, \ref{cost}). \\

Step 1: Update equation for $Z_{t}$. \\
Similar to the update equation in \cite{mcdonald2020exponential}, two separate cases need to be considered.

\noindent Suppose $\int_{\mathbb{X}}\prod_{i=1}^{N} h(x_{t},y_{t}^{i})Z_{t}(dx_{t}) \not = 0$, then 

\begin{align*} 
   &Z_{t+1}(\cdot)=P(x_{t+1} \in \cdot|y^{[1,N]}_{[0,t]}, u^{[1,N]}_{[0,t]})=\int \int_{\cdot} P(dx_{t+1}|x_{t}, y^{[1,N]}_{[0,t]}, u^{[1,N]}_{[0,t]})P(dx_{t}|y^{[1,N]}_{[0,t]}, u^{[1,N]}_{[0,t]}) \\
  &=\int \int_{\cdot} \tau(dx_{t+1}|x_{t}, u_{t}^{[1,N]})P(dx_{t}|y^{[1,N]}_{[0,t]}, u^{[1,N]}_{[0,t-1]}) \\
    &= \int \int_{\cdot}  \frac{\tau(dx_{t+1}|x_{t}, u_{t}^{[1,N]})Z_{t}(dx_{t})\prod_{i=1}^{N} h(x_{t},y_{t}^{i})}{\int_{\mathbb{X}}\prod_{i=1}^{N} h(x_{t},y_{t}^{i})Z_{t}(dx_{t})}=:F(Z_{t},\mathbf{u}_{t},y_{t}^{[1,N]})(\cdot)
\end{align*}

where $\tau$ is the transition kernel of the system (\ref{Prob 1}). If, $\int_{\mathbb{X}}\prod_{i=1}^{N} h(x_{t},y_{t}^{i})Z_{t}(dx_{t})  = 0$, then we define $Z_{t+1}(.)\equiv 0$\\

Step 2: Stochastic Kernel. Given the update equation above, it follows that the stochastic kernel for $Z_{t}$ is given by:
\begin{align*}
P(Z_{t+1} \in \cdot |Z_{[0,t]},\mathbf{f}_{[0,t]})&=\int P(Z_{t+1} \in \cdot |y_{t}^{[1,N]},Z_{[0,t]},\mathbf{f}_{[0,t]})P(dy_{t}^{[1,N]}  |Z_{[0,t]},\mathbf{f}_{[0,t]}) \\
&=\int \prod_{i=1}^{N}f^{i}_{t}(du^{i}_{t}|y^{i}_{t})P(Z_{t+1} \in \cdot |y_{t}^{[1,N]},Z_{t},u_{t}^{[1,N]})P(dy_{t}^{[1,N]}  |Z_{[0,t]})\\
&=\int \mathbbm{1}_{ \{ F(Z_t,\mathbf{u}_{t},y_{t}^{[1,N]}) \in \cdot \}} Z_{t}(dx_{t})\prod_{i=1}^{N} Q^{i}(dy_{t}^{i}|x_{t})f^{i}_{t}(du^{i}_{t}|y^{i}_{t})\\
&=:\eta(Z_{t+1} \in \cdot | Z_{t}, f_{t}^{[1,N]}) 
\end{align*}
Step 3: Cost-equivalence\\
Since $\sigma(Z_{t})\subseteq \sigma(I_{t}^{C})$, one can extend the set admissible policies so that $\gamma^{i}_{t}: Z_{t}\times I_{t}^{C} \mapsto f^{i}_{t}$ and is $\sigma(Z_{t}, I_{t}^{C})$ measurable. Consider the following:
\begin{align*}
   \sum_{t=0}^{\infty} \beta^{t} E^{\gamma}[c(x_{t},\mathbf{u_{t}})]&= \sum_{t=0}^{\infty} \beta^{t} E[E[c(x_{t},\mathbf{{u}_{t}})| Z_{t}, I_{t}^{C}]] \\
   &= \sum_{t=0}^{\infty} \beta^{t} E[\int Z_{t} (dx_{t})\prod_{i=1}^{N}f^{i}(du^{i}_{t}|y^{i}_{t})c(x_{t},u^{1}_{t},...,u^{N}_{t})]\\
   &=\sum_{t=0}^{\infty} \beta^{t} E[\Tilde{c}({Z}_{t},f^{1}_{t},...,f^{N}_{t})]
\end{align*}
where \begin{eqnarray}
    \tilde{c}({Z}_{t},f^{1}_{t},...,f^{N}_{t})=\int Z_{t} (dx_{t})\prod_{i=1}^{N}f^{i}(du^{i}_{t}|y^{i}_{t})Q^{i}(dy^{i}_{t}|x_{t})c(x_{t},u^{1}_{t},...,u^{N}_{t}).
\end{eqnarray}

\noindent Note that because $c$ is bounded, it follows by \cite[Theorem 15.13]{AlBo06} or \cite[p. 215]{Bogachev} (see also \cite[Theorem D.2.1]{lecturenotes}) that $\Tilde{c}$ is measurable with respect to $Z_{t}$ and $f_{t}^{1},...,f_{t}^{N}$. Because $(Z_{t},\mathbf{f}_{t})$ is an MDP, the desired result then follows from \cite[Theorem 5.1.1]{lecturenotes} which can be used to show that when the policy maps the state (a Polish space) and another Polish space onto an action space, which is also Polish, one can construct a Borel measurable policy which only relies on the state without any loss of optimality. Hence, it is sufficient to consider policies such that $\gamma^{i}_{t}: Z_{t} \mapsto f^{i}_{t}$ and is $\sigma(Z_{t})$-measurable.  
\qed

\begin{remark}[Belief vs. Predictor] In the context of a centralized POMDP the belief process refers to the conditional probability of the current state given the current information, i.e., $P(x_{t}\in .\mid y_{[0,t]},u_{[0,t-1]})$ whereas the predictor refers to the conditional probability of the next state given current information, i.e., $P(x_{t+1}\in . \mid y_{[0,t]},u_{[0,t-1]})$. While it is standard in the centralized setting, to express POMDPs in terms of an MDP with the belief process as the state  \cite{Yus76,Rhe74}, such a formulation can not be used identically in this setting. This is because from the centralized controller's perspective only the common information (one-step delayed measurements and actions) is available. This will also have implications for the weak-Feller property. \end{remark}

\begin{remark} [$N$-step delayed sharing pattern with $N\geq 2$] Witsenhausen \cite{wit71} conjectured that policies should be separable for the $N$-step delayed information sharing pattern regardless of $N$. \cite{VaraiyaWalrand} showed this to be incorrect by presenting a counterexample that shows that when the information structure is two-step delayed, policies do not exhibit a separation property. For further separation results see (\cite{NayyarMahajanTeneketzis2}, \cite{kur79}). In \cite{NayyarMahajanTeneketzis2}, a more general setup is considered, and it is shown that optimal policies may depend only on the belief in the state as well as the measurements of all agents and their private information conditioned on the common information. In \cite{kur79} a related separation result is given for the two-step delayed information sharing pattern. 
\end{remark}

\subsection{$K$-step periodic information sharing pattern as a centralized MDP} \label{K-step periodic information sharing pattern}
It was shown in \cite{OoiWornell} (see also \cite{YukTAC09}) that optimal policies for the $K$-step periodic information sharing pattern, under finite horizon cost criteria, are of separated form in the case where the measurement spaces, action spaces, and the state space are all finite. Similar to the argument used in \cite{yos75} and \cite{VaraiyaWalrand}, the authors of \cite{OoiWornell} use dynamic programming to show separability of optimal policies which requires a weak-Feller assumption. The latter is too restrictive for problems with standard Borel spaces. In this section, we will use an approach similar to the one we used for the one-step delayed information sharing pattern and the one used in \cite{YukTAC09} for the KSPISP problem with finite spaces. We will generalize this result to the case where all spaces are standard Borel and for the infinite horizon discounted cost criteria. Additionally, we will show that this problem reduces to a centralized problem. \\

We first recall the KSPISP. For problem (\ref{Prob 1},\ref{obs 1}), let $t=qK+r$ where $q$ and $r$ are non-negative integers such that $r\leq K-1$ and $K\in \mathbb{N}$ represent the period between successive information sharing between the Agents. Suppose that each Agent $i$ now has access to 
\begin{align} \label{k-step info}
   \nonumber  I^{i}_{t}&=\{y^{i}_{[qK,qK+r]},u^{i}_{[qK,qK+r-1]},I^{C}_{t}\} \\
     I^{C}_{t}&=\{y^{[1,N]}_{[0,qK-1]},u^{[1,N]}_{[0,qK-1]}\} \\
   \nonumber  I^{P^{i}}_{t}&= \{y^{i}_{[qK,qK+r]},u^{i}_{[qK,qK+r-1]}\} 
\end{align} 
 In order to reduce notation overload, $I^{i}_{t}$ is defined differently from the one-step delayed information sharing pattern case and consists of all the information available to Agent $i$ at time $t$ including common information.
 
Here, each Agent $i$ at time $t$ picks a policy $\gamma^{i}_{t}: I^{i}_{t} \mapsto u^{i}_{t}$ such that $\gamma^{i}_{t}$ is $\mathcal{B}(\mathbb{I}^{i}_{t} )$-measurable. The objective of the agents is to minimize 
\begin{equation} \label{Cost k-step}
   J(\gamma)=\sum_{t=K}^{\infty} \beta^{t} E^{\gamma}[c(x_{t},\mathbf{u_{t}})] 
\end{equation}
In \cite{OoiWornell,YukTAC09}, it was shown that, when the cost criteria has finite horizon, there is no loss in optimality if the agents replace their common information $(y^{[1,N]}_{[0,qK-1]},u^{[1,N]}_{[0,qK-1]})$ with $\pi_{t}=Z_{qK}=P(x_{qK}|y^{[1,N]}_{[0,qK-1]},u^{[1,N]}_{[0,qK-1]})$. In \cite{YukTAC09}, it was shown that when the state space, action space, and all observation spaces are finite, the problem can be reduced to a Markov decision process where the state is replaced by the predictor on the state, given the common information, at $K$-time periods. In what follows, we will extend this result to general spaces and show that for the cost (\ref{cost}) the policies satisfy a separation property. \\

\noindent{\bf Near Optimality of Finite Action Approximation}
The next proposition shows that, under certain assumptions, quantizing the action space so that only policies with finite range are considered leads to a near optimal solution as the quantization size increases.
\begin{proposition} \label{finite action approximation}
    Suppose the state space $\mathbb{X}$ and the action space $\mathbb{U}=\prod_{i=1}^{N}\mathbb{U}^{i}$ are compact. Additionally, assume that the transition kernel $\tau$ is continuous, under total variation, in $x$ and $u$ that is to say for any $(x_{n},\mathbf{u}_{n})\rightarrow (x,\mathbf{u})$: $\|\tau(x.|x_{n},\mathbf{u}_{n})-\tau(x.|x,\mathbf{u})\|_{TV}\rightarrow 0$. Let $\mathbb{U}_{m}=\{\mathbf{u}_{1},...,\mathbf{u}_{N_{m}}\}$ be such that for any $\mathbf{u}\in \mathbb{U}$ there exists $\mathbf{u}_{j}\in \mathbb{U}_{m}$ such that $d(\mathbf{u},\mathbf{u}_{j})\leq \frac{1}{m}$ where $m\in \mathbb{N}$. Let $\Gamma_{m}$ denote the set of admissible team policies with range in $\mathbb{U}_{m}$. Then, for any $\epsilon>0$ there exists $m'$ such that for all $m\geq m'$ we have that $\displaystyle |\inf_{\gamma\in \Gamma}J(\gamma)-\inf_{\gamma\in \Gamma_{m}}J(\gamma)|\leq \epsilon$.
\end{proposition}
\textbf{Proof.} See the Appendix \ref{Proof of Proposition on finite action approximation}.
\begin{remark}
    Note that the proof of Proposition \ref{finite action approximation} holds regardless of the nature of the information available to the agents. Hence, it holds for a variety of information structures including all the ones considered in this paper. For the one-step delayed information sharing pattern problem the same result can be obtained under different conditions through the use of the weak-Feller property of the centralized reduction of the problem (Theorem \ref{theorem5}). In particular, one can apply \cite[Theorem 3.16]{SaLiYuSpringer} in conjunction with the fact that quantized maps are dense in the set of maps endowed with the Young topology (see the proof of Theorem 5.1 in \cite{yuksel2023borkar}), to obtain the result of Proposition \ref{finite action approximation}. 
\end{remark}

\noindent{\bf Topology on the actions for the K-step periodic belief sharing pattern}
In view of Proposition \ref{finite action approximation}, we assume the following.
\begin{assumption} \label{finite action space} For all $i$, $\mathbb{U}^{i}$ is finite. 
\end{assumption}
    While for the one step delayed case finiteness of the action space was not required, we impose this assumption for the $K$-step periodic information structure in order to introduce a topology on the map-valued action space for the equivalent formulation that will lead to compactness of the new action space. In addition, this will prove helpful in establishing the weak-Feller property of the new transition kernel later on. Here, again, we assume without any loss in generality that the measurement channel of each individual Agent $i$ satisfies (\ref{Young top}).

\begin{definition} \label{top k-step}
For $i\in \{1,...,N\}$ and for all $n\in \mathbb{N}$ let $f^{i}_{qK,n}:\mathbb{Y}^{i}\rightarrow \mathcal{P}(\mathbb{U}^{i})$ be a measurable function and let $f^{i}_{qK}:\mathbb{Y}^{i}\rightarrow\mathcal{P}(\mathbb{U}^{i})$ be a measurable function. We say $f^{i}_{qK,n}\rightarrow f^{i}_{qK}$ if and only if $\mathcal{P}(\mathbb{Y}^{i}\times \mathbb{U}^{i})\ni \mu_{n}(.)=\int_{\cdot}\psi(dy^{i})f^{i}_{qK,n}(du^{i}_{qK}|y^{i}_{qK})\rightarrow \mu(.)=\int_{\cdot }\psi(dy^{i}_{qK})f^{i}_{qK}(du^{i}_{qK}|y^{i}_{qK})$ weakly, i.e., for all $g\in\mathbf{C_{b}}(\mathbb{Y}^{i}\times \mathbb{U}^{i})$: \[\int g(y^{i},u^{i})\psi(dy^{i})f^{i}_{qK,n}(du^{i}_{qK}|y^{i}_{qK})\rightarrow \int g(y^{i},u^{i})\psi(dy^{i})f^{i}_{qK}(du^{i}_{qK}|y^{i}_{qK})\]

Additionally, for $r>0$ we say $f^{i}_{qK+r,n}\rightarrow f^{i}_{qK+r}$ if and only if for all $u^{i}_{qK},...,u^{i}_{qK+r-1}$:\\   $\displaystyle \mathcal{P}((\mathbb{Y}^{i})^{r+1}\times \mathbb{U}^{i})\ni \mu_{n}(.)=\int_{\cdot}\prod_{t=qK}^{t=qK+r}\psi(dy^{i}_{t})f^{i}_{qK+r,n}(du^{i}_{qK+r}|y^{i}_{[qK,qK+r]},u^{i}_{[qK,qK+r-1]})\rightarrow \mu \in \mathcal{P}((\mathbb{Y}^{i})^{r+1}\times \mathbb{U}^{i})$ weakly. Where, \\
$\mu(.) =\int_{\cdot}\prod_{t=qK}^{t=qK+r}\psi(dy^{i}_{t})f^{i}_{qK+r}(du^{i}_{qK+r}|y^{i}_{[qK,qK+r]},u^{i}_{[qK,qK+r-1]})$. Define for each $i\in\{1,...,N\}$: 
$a^{i}_{q}=(f^{i}_{qK},...,f^{i}_{(q+1)K-1})$. Then, we say that $a^{i}_{q,n}=(f^{i}_{qK,n},...,f^{i}_{(q+1)K-1,n})\rightarrow a^{i}_{q}=(f^{i}_{qK},...,f^{i}_{(q+1)K-1})$ if and only if for all $r\in \{0,...,K-1\}$: $f^{i}_{qK+r,n}\rightarrow f^{i}_{qK+r}$.
\end{definition}
\begin{theorem} \label{themkstepMDP} The process $(\pi_{q},(a^{1}_{q},...,a^{N}_{q}); q\geq 1)$ forms a controlled Markov decision process. Here, for all $i$, $f^{i}_{qK}: y^{i}_{qK} \mapsto u^{i}_{t}$ and is $\mathcal{B}(\mathbb{Y}^{i})$-measurable and for $r>0$ $f^{i}_{qK+r}: (y^{i}_{[qK,qK+r]},u^{i}_{[qK,qK+r-1]}) \mapsto u^{i}_{t}$ and is $\displaystyle \mathcal{B}(\prod_{l=1}^{r+1}\mathbb{Y}^{i}\times \prod_{l=1}^{r}\mathbb{U}^{i})$-measurable.
\end{theorem}

In the remainder of the paper we will use the notation: $a_{q}=(a^{1}_{q},...,a^{N}_{q})$. We will use $\mathbb{A}$ to refer to the space in which these actions take values. 

\noindent\textbf{Proof.} To prove this theorem we first find an update equation for $\pi_{q}$ then we use that to compute the stochastic kernel and conclude that $(\pi_{q},(a^{1}_{q},...,a^{N}_{q}))$ is a controlled Markov chain.  

Step 1: Update equation for $\pi_{q}$
\begin{equation*}
 \begin{aligned}
      \pi_{q+1}&= Z_{(q+1)K} =F(Z_{(q+1)K-1},u_{(q+1)K-1}^{[1,N]},y_{(q+1)K-1}^{[1,N]})\\
    &=F(F(Z_{(q+1)K-2},u_{(q+1)K-2}^{[1,N]},y_{(q+1)K-2}^{[1,N]}),u_{(q+1)K-1}^{[1,N]},y_{(q+1)K-1}^{[1,N]})\\
    &=:G(\pi_{q},u_{[qK,qK+K-1]}^{[1,N]},y_{[qK,qK+K-1]}^{[1,N]})  
 \end{aligned}   
\end{equation*}
Step 2: Transition kernel
\begin{align}
  \nonumber  &P(\pi_{q+1}\in \cdot |\pi_{[0,q]},a_{[0,q]}^{[1,N]})=\int P(d\pi_{q+1}\in \cdot|y^{[1,N]}_{[qK,qK+K-1]},\pi_{[0,q]},a_{[0,q]}^{[1,N]}, u_{[qK,qK+K-1]}^{[1,N]})\\
  \nonumber   &P(dy^{[1,N]}_{[qK,qK+K-1]}, du_{[qK,qK+K-1]}^{[1,N]}|\pi_{[0,q]},a_{[0,q]}^{[1,N]})\\
  \nonumber   &=\int \mathbbm{1}_{\{G(\pi_{q},u_{[qK,qK+K-1]}^{[1,N]},y_{[qK,qK+K-1]}^{[1,N]}) \in \cdot\}}P(dy^{[1,N]}_{[qK,qK+K-1]}, du_{[qK,qK+K-1]}^{[1,N]}|\pi_{[0,q]},a_{[0,q]}^{[1,N]})\\
     &=: \theta(\pi_{q+1}\in \cdot|\pi_{[0,q]},a_{[0,q]}^{[1,N]}) \label{theta definition}
\end{align}
Where,
\begin{align*}
    &P(dy^{[1,N]}_{[qK,qK+K-1]},d u_{[qK,qK+K-1]}^{[1,N]}|\pi_{[0,q]},a_{[0,q]}^{[1,N]})= \\
    &\int \pi_{q}(dx_{Kq})(\ \prod_{i=1}^{N}Q^{i}(dy^{i}_{Kq}|x_{Kq})f^{i}_{qK}(du^{i}_{qK}|y^{i}_{Kq}))\tau(dx_{Kq+1}|x_{Kq},u^{[1,N]}_{Kq})\\
    &\big( \prod_{i=1}^{N}Q^{i}(dy^{i}_{Kq+1}|x_{Kq+1})f^{i}_{qK+1}(du^{i}_{qK+1}|y^{i}_{[Kq,Kq+1]},u^{i}_{qK})\big) \tau(dx_{Kq+2}|x_{Kq+1}, u^{[1,N]}_{Kq+1})...\\
    &\tau(dx_{qK+K-1}|x_{qK+K-2}, u^{[1,N]}_{qK+K-2})(\ \prod_{i=1}^{N}Q^{i}(dy^{i}_{qK+K-1}|x_{qK+K-1})\\
    &f^{i}_{qK+K-1}(du^{i}_{qK+K-1}|y^{i}_{[Kq,Kq+K-1]},u^{i}_{[qK,qK+K-2]}))
\end{align*}
Hence, $\theta(\pi_{q+1}\in .|\pi_{[0,q]},a_{[0,q]}^{[1,N]})=\theta(\pi_{q+1}\in .|\pi_{q},a_{q}^{[1,N]})$. 
\qed
\begin{theorem} \label{centralized reduction kstep step part two}
 Problem (\ref{Prob 1}),(\ref{obs 1}),(\ref{k-step info}),(\ref{Cost k-step}) is equivalent to one where the team policy at time $q$ is given by $(\gamma_{q}^{1},...,\gamma_{q}^{N})$ such that for all $i$ $\gamma^{i}_{q}:I_{qK}^{C} \mapsto a^{i}_{q}=(f^{i}_{qK},...,f^{i}_{(q+1)K-1})$ and is $\sigma(I_{qK}^{C})$-measurable.
 \end{theorem}

\textbf{Proof.} The proof here will follow similar steps to those of Theorem \ref{centralized reduction one step part one} with $I^{P^{i}}_{t}$ now acting as $y^{i}_{t}$. Let $(\gamma_{t}^{1},...,\gamma_{t}^{N})$ be a team policy as defined in Section 3. 

\noindent Step 1: Write
\begin{align*}
  & \sum_{t=K}^{t=\infty} \beta^{t} E[c(x_{t},\mathbf{u_{t}})]= \sum_{q=1}^{\infty}\sum_{t=qK}^{t=(q+1)K-1} \beta^{t} E[E[c(x_{t},\mathbf{{u}_{t}})| I^{C}_{t}]] \\
   &= \sum_{q=1}^{\infty}\sum_{t=qK}^{t=(q+1)K-1} \beta^{t} E[\int \pi_{q}(dx_{qK})(\prod_{i=1}^{N}Q^{i}(dy^{i}_{qK}|x_{qK})\gamma^{i}_{qK}(du^{i}_{qK}|I^{C}_{qK},y^{i}_{qK}))\\
   &\tau(dx_{qK+1}|x_{qK},\mathbf{u}_{qK})...\tau(dx_{t}|x_{t-1},\mathbf{u}_{t-1})(\prod_{i=1}^{N}Q^{i}(dy^{i}_{t}|x_{t})\gamma^{i}_{t}(du^{i}_{t }|y^{i}_{[qK,qK+r]},u^{i}_{[qK,qK+r-1]},I^{C}_{t}))c(x_{t},\mathbf{u}_{t})]
\end{align*}


Step 2 and 3 on the equivalent formulation under measurable selections proceed in a similar manner to the one-step delayed information sharing pattern case with $I^{P^{i}}$ now acting as a stand-in for $y^{i}$. \\

Step 4: Thus, we get that 
\begin{align*}
  & \sum_{t=K}^{t=\infty} \beta^{t} E[c(x_{t},\mathbf{u_{t}})]=\sum_{q=1}^{\infty}  \beta^{qK} E[\int \pi_{q}(dx_{qK}) \sum_{t=qK}^{t=(q+1)K-1} \beta^{t-qK} (\prod_{i=1}^{N}Q^{i}(dy^{i}_{qK}|x_{qK})\tilde{\gamma}^{i}_{qK}(I^{C}_{qK})(y^{i}_{qK}))  \\
  &\tau(dx_{qK+1}|x_{qK},\mathbf{u}_{qK})... \tau(dx_{t}|x_{t-1},\mathbf{u}_{t-1})(\prod_{i=1}^{N}Q^{i}(dy^{i}_{t}|x_{t})\tilde{\gamma}^{i}_{t}(I^{C}_{qK})(y^{i}_{[qK,qK+r]},u^{i}_{[qK,qK+r-1]}))c(x_{t},\mathbf{u}_{t})]\\ \qquad
   &= \sum_{q=1}^{\infty}(\beta^{K})^{q} E^{\Tilde{\gamma}}[\Tilde{c}(\pi_{q},a^{1}_{q},...,a^{N}_{q})] \quad
\end{align*}


where, 
\begin{align}  \label{New k-step Cost}
   \nonumber &\Tilde{c}(\pi_{q},a^{1}_{q},...,a^{N}_{q})= \\
   \nonumber &\int \pi_{q}(dx_{qK})\sum_{t=qK}^{t=(q+1)K-1} \beta^{t-qK} (\prod_{i=1}^{N}Q^{i}(dy^{i}_{qK}|x_{qK})f^{i}_{qK}(du^{i}_{qK}|y^{i}_{qK})) 
    \tau(dx_{qK+1}|x_{qK},\mathbf{u}_{qK})\\
    &...\tau(dx_{t}|x_{t-1},\mathbf{u}_{t-1})(\prod_{i=1}^{N}Q^{i}(dy^{i}_{t}|x_{t})f^{i}_{t}(du^{i}_{t }|y^{i}_{[qK,qK+r]},u^{i}_{[qK,qK+r-1]}))c(x_{t},\mathbf{u}_{t})
\end{align}

Here, since $c$ is bounded, the measurability of $\Tilde{c}$ with respect to $\pi_{q}$ and $a^{1}_{q},...,a^{N}_{q}$ follows by \cite[Theorem 15.13]{AlBo06} or \cite[p. 215]{Bogachev} (see also \cite[Theorem D.2.1 ]{lecturenotes}). \qed

\begin{corollary} \label{K-periodic belief sharing} The controlled Markov model $(\pi_{q},(a^{1}_{q},...,a^{N}_{q}); q\geq 1)$ defines an equivalent problem to the problem (\ref{Prob 1},\ref{obs 1}, \ref{k-step info}, \ref{Cost k-step}) under an appropriate choice of the cost criteria, i.e, a reformulation of the cost given in (\ref{Cost k-step}). Additionally, one can replace the policies $\{\Tilde{\gamma}^{i}\}$ with policies $\{\gamma^{i}\}$ such that $\gamma^{i}_{q}: \pi_{q} \mapsto a^{i}_{q}$ and $\gamma^{i}_{q}$ is $\sigma(\pi_{q})$-measurable. 
\end{corollary}

\textbf{Proof.} It suffices to consider the following cost:
\begin{equation}
    J=\sum_{q=1}^{\infty} (\beta^{K})^{q} E^{\Tilde{\gamma}}[\Tilde{c} \label{mykperiod}(\pi_{q},a^{1}_{q},...,a^{N}_{q})]
\end{equation}
Since $\sigma(\pi_{q})\subseteq \sigma(I_{qK}^{C})$, one can, without loss of optimality, extend the set admissible policies so that $\gamma^{i}_{q}: \pi_{q}\times I_{qK}^{C} \mapsto a^{i}_{q}$ and is $\mathcal
\sigma(\pi_{q}, I_{qK}^{C})$ measurable. As earlier, the desired result then follows from \cite[Theorem 2]{Blackwell3}. \qed
\begin{remark}
    While the one-step delayed information sharing pattern is a special case of the $K$-step periodic information sharing pattern, the one-step information structure plays a crucial role in establishing properties of the more general $K$-step periodic information sharing pattern. In this section, we saw that the update equation for the K-step predictor can only be obtained through recursively applying the update equation for the one-step predictor. As will be seen in the next section, in order to obtain the weak-Feller regularity for the $K$-step periodic information sharing pattern problem it will be useful to first establish weak-Feller regularity for the one-step delayed information sharing pattern problem.   
\end{remark}
\section{Weak-Feller Property of the Equivalent MDPs, Existence of Optimal Policies, and their Rigorous Approximations}
In this section, we will use the structural results established in Section \ref{equivalent formulation} as well as the topologies on the map-valued actions introduced for the one-step delayed information sharing pattern and the $K$-step periodic information sharing pattern in order to find under which conditions the kernels $\eta$ and $\theta$ are weak-Feller i.e. weakly continuous. The weak-Feller property is important as it enables one to establish the existence of optimal solutions. In this section as well as in the remainder of the paper, we will assume that the state space $\mathbb{X}$ is compact.  

\subsection{Weak-Feller property for the one-step delayed sharing pattern}

\begin{theorem}\label{theorem5}
Suppose the transition kernel $\tau(x_{t+1}\in \cdot |x_{t}=x,\mathbf{u}_{t}=\mathbf{u})$ is weakly continuous in $x$ and $u$ that is to say for any $(x_{n},\mathbf{u}_{n})\rightarrow (x,\mathbf{u})$ and any $g\in \mathbf{C_{b}}(\mathbb{X})$: $\int g(x)\tau(dx|x_{n},\mathbf{u}_{n})\rightarrow \int g(x)\tau(dx|x,\mathbf{u})$. Additionally suppose $h(x,y)$ is continuous in $x$ (for all $y$: $h(.,y): \mathbb{X} \rightarrow \mathbb{R}$ is a continuous function). Then, the MDP $(Z_{t},f_{t})$ introduced in Theorem \ref{them2young} is weak-Feller.
\end{theorem}
\textbf{Proof.} See the Appendix.


\begin{corollary} \label{existence OSDISP}
    Suppose the transition kernel $\tau(x_{t+1}\in \cdot |x_{t}=x,\mathbf{u}_{t}=\mathbf{u})$ is weakly continuous in $x$ and $\mathbf{u}$. Additionally suppose $h(x,y)$ is continuous in $x$. Then, an optimal solution to the centralized reduction of the OSDISP problem (\ref{Prob 1}),(\ref{obs 1}),(\ref{cost}) exists and stationary policies are optimal.
\end{corollary}
\textbf{Proof}. See the Appendix. 
\subsection{Weak-Feller property for the K-step periodic information sharing pattern}
\begin{theorem}\label{theorem6} Suppose the transition kernel $\tau\left(x_{t+1} \in . \mid x_{t}=x, \mathbf{u}_{t}=\mathbf{u}\right)$ is continuous under total variation in $x$ and $\mathbf{u}$, and $h(x, y)$ is continuous in $x$. Then the MDP $\left(\pi q ; \mathbf{a}_{q}\right)$ introduced in Theorem \ref{themkstepMDP} is weak-Feller.
\end{theorem}
\begin{remark}
    Note that since all the marginals on $\mathbf{Y}^{1},...,\mathbf{Y}^{N}$ are fixed, once can without any loss of generality assume that $\mathbf{Y}^{1},...,\mathbf{Y}^{N}$ are compact.  
\end{remark}
\textbf{Proof.}  
Let $g\in C_{b}(\mathcal{P}(\mathbf{X}))$. 
$$
\begin{aligned}
&\left|E\left[g\left(\pi_{q+1}\right) \mid \pi_{q}^{n}, \mathbf{a}_{q, n}\right]-E\left[g\left(\pi_{q+1}\right) \mid \pi_{q}, \mathbf{a}_{q}\right]\right| \\
& \leq V_{1}^{n}+V_{2}^{n}
\end{aligned}
$$
\\
where, $\displaystyle V_{1}^{n}= \bigg| E\left[g\left(\pi_{q+1}\right) \mid \pi_{q}^{n}, \mathbf{a}_{q,n}\right] - E\left[g\left(\pi_{q+1}\right) \mid \pi_{q}^{n}, \mathbf{a}_{q} \right] \bigg|$\\

\noindent and $\displaystyle V_{2}^{n}= \bigg| E\left[g\left(\pi_{q+1}\right) \mid \pi_{q}^{n}, \mathbf{a}_{q}\right] - E\left[g\left(\pi_{q+1}\right) \mid \pi_{q}, \mathbf{a}_{q} \right] \bigg|$.
\\

Next, we will show that $V_{1}^{n} \rightarrow 0$ and $V_{2}^{n} \rightarrow 0$ as $n \rightarrow \infty$.\\

\begin{eqnarray*}
 &&V_{2}^{n}=\left|\int g\left(\pi_{q+1}\right) \theta\left(\pi_{q+1} \mid \pi_{q}^{n}, \mathbf{a}_{q}\right)-\int g\left(\pi_{q+1}\right) \theta\left(\pi_{q+1} \mid \pi_{q}, \mathbf{a}_{q}\right)\right|\\
&& =\bigg| \int g(\pi_{q+1})\chi_{\{G(\pi_{q}^{n},\mathbf{u}_{[qK,qK+K-1]},y_{[qK,qK+K-1]}^{[1,N]}) \}}\\
&&P(dy_{[qK,qK+K-1]}^{[1,N]}, du_{[qK,qK+K-1]}^{[1,N]}\mid\pi_{q}^{n}, \mathbf{a}_{q})-\\
&&\int g(\pi_{q+1})\chi_{\{G(\pi_{q},\mathbf{u}_{[qK,qK+K-1]},y_{[qK,qK+K-1]}^{[1,N]}) \}}\\
&&P(dy_{[qK,qK+K-1]}^{[1,N]}, du_{[qK,qK+K-1]}^{[1,N]}\mid\pi_{q}, \mathbf{a}_{q})  \bigg|\\
&&\leq \bigg| \int g(\pi_{q+1})\chi_{\{G(\pi_{q}^{n},\mathbf{u}_{[qK,qK+K-1]},y_{[qK,qK+K-1]}^{[1,N]}) \}}\\
&&P(dy_{[qK,qK+K-1]}^{[1,N]}, du_{[qK,qK+K-1]}^{[1,N]}\mid\pi_{q}^{n}, \mathbf{a}_{q})-\\
&&\int g(\pi_{q+1})\chi_{\{G(\pi_{q}^{n},\mathbf{u}_{[qK,qK+K-1]},y_{[qK,qK+K-1]}^{[1,N]}) \}}\\
&&P(dy_{[qK,qK+K-1]}^{[1,N]}, du_{[qK,qK+K-1]}^{[1,N]}\mid\pi_{q}, \mathbf{a}_{q})  \bigg|+\\
&&\bigg| \int g(\pi_{q+1})\chi_{\{G(\pi_{q}^{n},\mathbf{u}_{[qK,qK+K-1]},y_{[qK,qK+K-1]}^{[1,N]}) \}}\\
&&P(dy_{[qK,qK+K-1]}^{[1,N]}, du_{[qK,qK+K-1]}^{[1,N]}\mid\pi_{q}, \mathbf{a}_{q})-\\
&&\int g(\pi_{q+1})\chi_{\{G(\pi_{q},\mathbf{u}_{[qK,qK+K-1]},y_{[qK,qK+K-1]}^{[1,N]})\}}\\
&&P(dy_{[qK,qK+K-1]}^{[1,N]}, du_{[qK,qK+K-1]}^{[1,N]}\mid\pi_{q}, \mathbf{a}_{q})\bigg|
\end{eqnarray*}
For the second term in the last inequality, we have that
\begin{eqnarray*}
    && \bigg| \int g(\pi_{q+1})\chi_{\{G(\pi_{q}^{n},\mathbf{u}_{[qK,qK+K-1]},y_{[qK,qK+K-1]}^{[1,N]}) \}}\\
    &&P(dy_{[qK,qK+K-1]}^{[1,N]}, du_{[qK,qK+K-1]}^{[1,N]}\mid\pi_{q}, \mathbf{a}_{q})\\
&&-\int g(\pi_{q+1})\chi_{\{G(\pi_{q},\mathbf{u}_{[qK,qK+K-1]},y_{[qK,qK+K-1]}^{[1,N]}) \}}\\
&&P(dy_{[qK,qK+K-1]}^{[1,N]}, du_{[qK,qK+K-1]}^{[1,N]}\mid\pi_{q}, \mathbf{a}_{q})\bigg|\\
&&=\bigg| \int g(G(\pi_{q}^{n},\mathbf{u}_{[qK,qK+K-1]},y_{[qK,qK+K-1]}^{[1,N]}) )\\
&&P(dy_{[qK,qK+K-1]}^{[1,N]}, du_{[qK,qK+K-1]}^{[1,N]}\mid\pi_{q}, \mathbf{a}_{q})-\\
&&\int g(G(\pi_{q},\mathbf{u}_{[qK,qK+K-1]},y_{[qK,qK+K-1]}^{[1,N]}))\\
&&P(dy_{[qK,qK+K-1]}^{[1,N]}, du_{[qK,qK+K-1]}^{[1,N]}\mid\pi_{q}, \mathbf{a}_{q})\bigg|\rightarrow 0\\
\end{eqnarray*}
The latter follows from the fact that $G$ is obtained through recursive applications of the update function for the one-step predictor: $F$ which was shown to be continuous in the predictor in the proof of the weak-Feller property for the one-step problem. Thus $G$ is continuous in the predictor and the desired result follows by an application of DCT. Hence, in order to show that $V_{2}^{n}\rightarrow 0$, it is sufficient to show that\\
\begin{eqnarray*}
    &&V_{3}^{n}:=\bigg| \int g(\pi_{q+1})\chi_{\{G(\pi_{q}^{n},\mathbf{u}_{[qK,qK+K-1]},y_{[qK,qK+K-1]}^{[1,N]}) \}}\\
    &&P(dy_{[qK,qK+K-1]}^{[1,N]}, du_{[qK,qK+K-1]}^{[1,N]}\mid\pi_{q}^{n}, \mathbf{a}_{q})-\\
    &&\int g(\pi_{q+1})\chi_{\{G(\pi_{q}^{n},\mathbf{u}_{[qK,qK+K-1]},y_{[qK,qK+K-1]}^{[1,N]}) \}}\\
    &&P(dy_{[qK,qK+K-1]}^{[1,N]}, du_{[qK,qK+K-1]}^{[1,N]}\mid\pi_{q}, \mathbf{a}_{q})  \bigg|\rightarrow 0 
\end{eqnarray*}
We have that
\begin{eqnarray*}
&&V_{3}^{n}\leq \|g\|_{\infty}||P(dy_{[qK,qK+K-1]}^{[1,N]}, du_{[qK,qK+K-1]}^{[1,N]}\mid\pi_{q}^{n}, \mathbf{a}_{q})-\\
&&P(dy_{[qK,qK+K-1]}^{[1,N]}, du_{[qK,qK+K-1]}^{[1,N]}\mid\pi_{q}, \mathbf{a}_{q})||_{TV}\\
&& =2 \|g\|_{\infty}\sup _{A} \bigg| \int \pi_{q}^{n}\left(d x_{q K}\right)[\int \left(\prod _ { i = 1 } ^ { N } Q ^ { i } ( d y _ { q K } ^ { i } | x _ { qK }  )f^{i}_{qK}(du^{i}_{qK}|y _ { q K } ^ { i })\right) \\
&&\tau(d x_{q K+1} \mid x_{q K}, u_{q K})...\tau(dx_{ qK+K-1} \mid x_{qK+K-2}, u_{qK+K-2})\\
&&\prod_{i=1}^{N} Q^{i}\left(dy_{ qK+K-1}^{i} \mid x_{qK+K-1}\right)f^{i}_{qK+K-1}(du^{i}_{qK+K-1}|y^{i}_{[qK,qK+K-1]},y^{i}_{[qK,qK+K-2]})] \\
&&- \int \pi_{q}\left(d x_{q K}\right)[\int \left(\prod _ { i = 1 } ^ { N } Q ^ { i } ( d y _ { q K } ^ { i } | x _ {q K }  )f^{i}_{qK}(du^{i}_{qK}|y _ { q K } ^ { i })\right) \tau(d x_{q K+1} \mid x_{q K} u_{q K})...  \\
&&\tau(dx_{ qK+K-1} \mid x_{qK+K-2}, u_{qK+K-2})\prod_{i=1}^{N} Q^{i}\left(dy_{ qK+K-1}^{i} \mid x_{qK+K-1}\right)\\
&&f^{i}_{qK+K-1}(du^{i}_{qK+K-1}|y^{i}_{[qK,qK+K-1]},y^{i}_{[qK,qk+K-2]})]\bigg|
\end{eqnarray*}

Since $\tau$ is continuous in $x$ under total variation and $h(x,y)$ is continuous in $x$, $\{\Tilde{G}_{A}(x)\}_{A}$, where \\
$\Tilde{G}_{A}(x)=\int \left(\prod _ { i = 1 } ^ { N } Q ^ { i } ( d y _ { q K } ^ { i } | x _ { qK } ) f^{i}_{qK}(du^{i}_{qK}|y _ { q K } ^ { i })\right) \tau(d x_{q K+1} \mid x_{q K} u_{q K})...\\ \tau(dx_{ qK+K-1} \mid x_{qK+K-2}, u_{qK+K-2}) \prod_{i=1}^{N} Q^{i}\left(dy_{ qK+K-1}^{i} \mid x_{qK+K-1}\right)$\\
$f^{i}_{qK+K-1}(du^{i}_{qK+K-1}|y^{i}_{[qK,qK+K-1]},y^{i}_{[qK,qK+K-2]})$, satisfies the condition
of lemma \ref{lemma weak feller one step}. It, then, follows that

$$
\begin{aligned}
& \sup _{A}\left|\int \pi_{q}^{n}(dx_{qK})\Tilde{G}_{A}(x)-\int \pi_{q}(dx_{qK})\Tilde{G}_{A}(x)\right|\rightarrow 0
\end{aligned}
$$
as $n \rightarrow 0$. Hence, $V_{3}^{n}\rightarrow 0$. Next, we will show that $V_{1}^{n}\rightarrow 0$. We have that\\
\begin{eqnarray*}
    &&V_{1}^{n}=\bigg| E\left[g\left(\pi_{q+1}\right) \mid \pi_{q}^{n}, \mathbf{a}_{q,n}\right] - E\left[g\left(\pi_{q+1}\right) \mid \pi_{q}^{n}, \mathbf{a}_{q} \right] \bigg |\\
    &&=\bigg| \int g(G(\pi_{q}^{n},\mathbf{u}_{[qK,qK+K-1]},y_{[qK,qK+K-1]}^{[1,N]}) )\\
    &&P(dy_{[qK,qK+K-1]}^{[1,N]}, du_{[qK,qK+K-1]}^{[1,N]}\mid\pi_{q}^{n}, \mathbf{a}_{q,n})-\\
&&\int g(G(\pi_{q}^{n},\mathbf{u}_{[qK,qK+K-1]},y_{[qK,qK+K-1]}^{[1,N]}))\\
&&P(dy_{[qK,qK+K-1]}^{[1,N]}, du_{[qK,qK+K-1]}^{[1,N]}\mid\pi_{q}^{n}, \mathbf{a}_{q})\bigg|
\end{eqnarray*}
It follows from the proof of Theorem \ref{theorem5} that $F$ is continuous and hence $G$ is continuous which entails that $g\circ G$ is continuous. Thus, because $\mathcal{\mathbf{X}}$ is compact, it is sufficient to check that\\
\begin{eqnarray*}
    &&\sup_{\|g\|_{BL}\leq 1}\bigg| \int g(y_{[qK,qK+K-1]}^{[1,N]}, u_{[qK,qK+K-1]}^{[1,N]} )\\
    &&P(dy_{[qK,qK+K-1]}^{[1,N]}, du_{[qK,qK+K-1]}^{[1,N]}\mid\pi_{q}^{n}, \mathbf{a}_{q,n})-\\
    &&\int g(y_{[qK,qK+K-1]}^{[1,N]}, u_{[qK,qK+K-1]}^{[1,N]} )\\
    &&P(dy_{[qK,qK+K-1]}^{[1,N]}, du_{[qK,qK+K-1]}^{[1,N]}\mid\pi_{q}^{n}, \mathbf{a}_{q})\bigg|=\\
    &&\rho(P(dy_{[qK,qK+K-1]}^{[1,N]}, du_{[qK,qK+K-1]}^{[1,N]}\mid\pi_{q}^{n}, \mathbf{a}_{q,n}),\\
    &&P(dy_{[qK,qK+K-1]}^{[1,N]},du_{[qK,qK+K-1]}^{[1,N]}\mid\pi_{q}^{n}, \mathbf{a}_{q}))\rightarrow 0
\end{eqnarray*}
Now,
\begin{eqnarray*}
    &&\sup_{\|g\|_{BL}\leq 1}\bigg| \int g(y_{[qK,qK+K-1]}^{[1,N]}, u_{[qK,qK+K-1]}^{[1,N]} )\\
    &&P(dy_{[qK,qK+K-1]}^{[1,N]}, du_{[qK,qK+K-1]}^{[1,N]}\mid\pi_{q}^{n}, \mathbf{a}_{q,n})-\\
    &&\int g(y_{[qK,qK+K-1]}^{[1,N]}, u_{[qK,qK+K-1]}^{[1,N]} )\\
    &&P(dy_{[qK,qK+K-1]}^{[1,N]}, du_{[qK,qK+K-1]}^{[1,N]}\mid\pi_{q}^{n}, \mathbf{a}_{q})\bigg|\\
    &&=  \sup_{\|g\|_{BL}\leq 1}\bigg|\int g \pi_{q}^{n}(dx_{Kq})(\ \prod_{i=1}^{N}\psi(dy^{i}_{qK})h(x_{qK},y^{i}_{qK})f^{i}_{n,qK}(du^{i}_{Kq}|y_{Kq}^{i}))\\
    &&\tau(dx_{Kq+1}|x_{Kq},u^{[1,N]}_{Kq})(\prod_{i=1}^{N}\psi(dy^{i}_{qK+1})h(x_{qK+1},y^{i}_{qK+1})\\
    &&f^{i}_{n,qK+K-1}(du^{i}_{qK+K-1}|y^{i}_{[qK,qK+K-1]},y^{i}_{[qK,qK+K-2]}))...\\ 
    &&\prod_{i=1}^{N}f^{i}_{n,Kq+K-2}(du^{i}_{Kq+K-2}|y_{[qK,qK+K-2]}^{i},u_{[qK,qK+K-3]}^{i})\\
    &&\tau(dx_{qK+K-1}|x_{qK+K-2}, u^{[1,N]}_{qK+K-2})(\ \prod_{i=1}^{N}\psi(dy^{i}_{qK+K-1})h(x_{qK+K-1},y^{i}_{qK+K-1}))\\
    &&f^{i}_{n,qK+K-1}(du^{i}_{qK+K-1}|y^{i}_{[qK,qK+K-1]},y^{i}_{[qK,qK+K-2]})-\\
    &&\int g\pi_{q}^{n}(dx_{Kq})(\ \prod_{i=1}^{N}\psi(dy^{i}_{qK})h(x_{qK},y^{i}_{qK})f^{i}_{qK}(du^{i}_{Kq}|y_{Kq}^{i}))\tau(dx_{Kq+1}|x_{Kq},u^{[1,N]}_{Kq})\\
    &&(\prod_{i=1}^{N}\psi(dy^{i}_{qK+1})h(x_{qK+1},y^{i}_{qK+1}))...\prod_{i=1}^{N}f^{i}_{qK+K-2}(du^{i}_{Kq+K-2}|y_{[qK,qK+K-2]}^{i},u_{[qK,qK+K-3]}^{i})\\
    &&\tau(dx_{qK+K-1}|x_{qK+K-2}, u^{[1,N]}_{qK+K-2})(\ \prod_{i=1}^{N}\psi(dy^{i}_{qK+K-1})h(x_{qK+K-1},y^{i}_{qK+K-1})\\
    &&f^{i}_{qK+K-1}(du^{i}_{qK+K-1}|y^{i}_{[qK,qK+K-1]},y^{i}_{[qK,qK+K-2]}))\bigg|\\
    && \leq \int  \pi_{q}^{n}(dx_{Kq})\sup_{\|g\|_{BL}\leq 1}\bigg|\int g (\ \prod_{i=1}^{N}\psi(dy^{i}_{qK})h(x_{qK},y^{i}_{qK})f^{i}_{n,qK}(du^{i}_{Kq}|y_{Kq}^{i}))\\
    &&\tau(dx_{Kq+1}|x_{Kq},u^{[1,N]}_{Kq})(\prod_{i=1}^{N}\psi(dy^{i}_{qK+1})h(x_{qK+1},y^{i}_{qK+1})\\
    &&f^{i}_{n,qK+K-1}(du^{i}_{qK+K-1}|y^{i}_{[qK,qK+K-1]},y^{i}_{[qK,qK+K-2]}))...\\ 
    &&\prod_{i=1}^{N}f^{i}_{n,Kq+K-2}(du^{i}_{Kq+K-2}|y_{[qK,qK+K-2]}^{i},u_{[qK,qK+K-3]}^{i})\\
    &&\tau(dx_{qK+K-1}|x_{qK+K-2}, u^{[1,N]}_{qK+K-2})(\ \prod_{i=1}^{N}\psi(dy^{i}_{qK+K-1})h(x_{qK+K-1},y^{i}_{qK+K-1}))\\
    &&f^{i}_{n,qK+K-1}(du^{i}_{qK+K-1}|y^{i}_{[qK,qK+K-1]},y^{i}_{[qK,qK+K-2]})-\int g\\
    &&(\ \prod_{i=1}^{N}\psi(dy^{i}_{qK})h(x_{qK},y^{i}_{qK})f^{i}_{qK}(du^{i}_{Kq}|y_{Kq}^{i}))\tau(dx_{Kq+1}|x_{Kq},u^{[1,N]}_{Kq})(\prod_{i=1}^{N}\psi(dy^{i}_{qK+1})\\
    &&h(x_{qK+1},y^{i}_{qK+1}))...\prod_{i=1}^{N}f^{i}_{qK+K-2}(du^{i}_{Kq+K-2}|y_{[qK,qK+K-2]}^{i},u_{[qK,qK+K-3]}^{i})\\
    &&\tau(dx_{qK+K-1}|x_{qK+K-2}, u^{[1,N]}_{qK+K-2})(\ \prod_{i=1}^{N}\psi(dy^{i}_{qK+K-1})h(x_{qK+K-1},y^{i}_{qK+K-1})\\
    &&f^{i}_{qK+K-1}(du^{i}_{qK+K-1}|y^{i}_{[qK,qK+K-1]},y^{i}_{[qK,qK+K-2]}))\bigg|\\
\end{eqnarray*}
Because the supremum in the last inequality is taken over a compact set, it is sufficient to check that for countably many functions $g_{m}$, we have that\\
\begin{eqnarray*}
    &&\int  \pi_{q}^{n}(dx_{Kq})\bigg|\int g_{m} (\ \prod_{i=1}^{N}\psi(dy^{i}_{qK})h(x_{qK},y^{i}_{qK})f^{i}_{n,qK}(du^{i}_{Kq}|y_{Kq}^{i}))\\
    &&\tau(dx_{Kq+1}|x_{Kq},u^{[1,N]}_{Kq})(\prod_{i=1}^{N}\psi(dy^{i}_{qK+1})h(x_{qK+1},y^{i}_{qK+1})\\
    &&f^{i}_{n,qK+K-1}(du^{i}_{qK+K-1}|y^{i}_{[qK,qK+K-1]},y^{i}_{[qK,qK+K-2]}))...\\ 
    &&\prod_{i=1}^{N}f^{i}_{n,Kq+K-2}(du^{i}_{Kq+K-2}|y_{[qK,qK+K-2]}^{i},u_{[qK,qK+K-3]}^{i})\\
    &&\tau(dx_{qK+K-1}|x_{qK+K-2}, u^{[1,N]}_{qK+K-2})(\ \prod_{i=1}^{N}\psi(dy^{i}_{qK+K-1})h(x_{qK+K-1},y^{i}_{qK+K-1}))\\
    &&f^{i}_{n,qK+K-1}(du^{i}_{qK+K-1}|y^{i}_{[qK,qK+K-1]},y^{i}_{[qK,qK+K-2]})-\int g_{m}\\
    &&(\ \prod_{i=1}^{N}\psi(dy^{i}_{qK})h(x_{qK},y^{i}_{qK})f^{i}_{qK}(du^{i}_{Kq}|y_{Kq}^{i}))\tau(dx_{Kq+1}|x_{Kq},u^{[1,N]}_{Kq})(\prod_{i=1}^{N}\psi(dy^{i}_{qK+1})\\
    &&h(x_{qK+1},y^{i}_{qK+1}))...\prod_{i=1}^{N}f^{i}_{qK+K-2}(du^{i}_{Kq+K-2}|y_{[qK,qK+K-2]}^{i},u_{[qK,qK+K-3]}^{i})\\
    &&\tau(dx_{qK+K-1}|x_{qK+K-2}, u^{[1,N]}_{qK+K-2})(\ \prod_{i=1}^{N}\psi(dy^{i}_{qK+K-1})h(x_{qK+K-1},y^{i}_{qK+K-1})\\
    &&f^{i}_{qK+K-1}(du^{i}_{qK+K-1}|y^{i}_{[qK,qK+K-1]},y^{i}_{[qK,qK+K-2]}))\bigg|\rightarrow 0
\end{eqnarray*}
Define\\
\begin{eqnarray*}
    &&\bigg|\int g_{m} (\ \prod_{i=1}^{N}\psi(dy^{i}_{qK})h(x_{qK},y^{i}_{qK})f^{i}_{n,qK}(du^{i}_{Kq}|y_{Kq}^{i}))\tau(dx_{Kq+1}|x_{Kq},u^{[1,N]}_{Kq})\\
    &&(\prod_{i=1}^{N}\psi(dy^{i}_{qK+1})h(x_{qK+1},y^{i}_{qK+1})\\
    &&f^{i}_{n,qK+K-1}(du^{i}_{qK+K-1}|y^{i}_{[qK,qK+K-1]},y^{i}_{[qK,qK+K-2]}))...\\ 
    &&\prod_{i=1}^{N}f^{i}_{n,Kq+K-2}(du^{i}_{Kq+K-2}|y_{[qK,qK+K-2]}^{i},u_{[qK,qK+K-3]}^{i})\\
    &&\tau(dx_{qK+K-1}|x_{qK+K-2}, u^{[1,N]}_{qK+K-2})(\ \prod_{i=1}^{N}\psi(dy^{i}_{qK+K-1})h(x_{qK+K-1},y^{i}_{qK+K-1}))\\
    &&f^{i}_{n,qK+K-1}(du^{i}_{qK+K-1}|y^{i}_{[qK,qK+K-1]},y^{i}_{[qK,qK+K-2]})-\int g_{m}\\
    &&(\ \prod_{i=1}^{N}\psi(dy^{i}_{qK})h(x_{qK},y^{i}_{qK})f^{i}_{qK}(du^{i}_{Kq}|y_{Kq}^{i}))\tau(dx_{Kq+1}|x_{Kq},u^{[1,N]}_{Kq})(\prod_{i=1}^{N}\psi(dy^{i}_{qK+1})\\
    &&h(x_{qK+1},y^{i}_{qK+1}))...\prod_{i=1}^{N}f^{i}_{qK+K-2}(du^{i}_{Kq+K-2}|y_{[qK,qK+K-2]}^{i},u_{[qK,qK+K-3]}^{i})\\
    &&\tau(dx_{qK+K-1}|x_{qK+K-2}, u^{[1,N]}_{qK+K-2})(\ \prod_{i=1}^{N}\psi(dy^{i}_{qK+K-1})h(x_{qK+K-1},y^{i}_{qK+K-1})\\
    &&f^{i}_{qK+K-1}(du^{i}_{qK+K-1}|y^{i}_{[qK,qK+K-1]},y^{i}_{[qK,qK+K-2]}))\bigg|:=E(x_{qk})
\end{eqnarray*}
It is sufficient to check that $\displaystyle \sup_{x_qK}E(x_{qk})\rightarrow 0$.\\
In what follows, we will assume that $K=3$ for ease of presentation. The more general case can be proved using a similar process. We have that\\
\begin{eqnarray*}
    &&\sup_{x_{qK}}E(x_{qK})\leq \sup_{x_{qK}}M_{1}^{n}+\sup_{x_{qK}}M_{2}^{n}+\sup_{x_{qK}}M^{n}_{3}
\end{eqnarray*}
where, 
\begin{eqnarray*}
    && M_{1}^{n}=  \bigg|\int g_{m}\prod_{i=1}^{N}\psi(dy^{i}_{3q})h(x_{3q},y^{i}_{3q})f^{i}_{n,3q}(du^{i}_{3q}|y_{3q}^{i})\tau(dx_{3q+1}|x_{3q},u^{[1,N]}_{3q})(\prod_{i=1}^{N}\psi(dy^{i}_{3q+1})\\
    &&h(x_{3q+1},y^{i}_{3q+1})f^{i}_{n,3q+1}(du^{i}_{3q+1}|y_{[3q,3q+1]}^{i},u^{i}_{3q}))\tau(dx_{3q+2}|x_{3q+1},u^{[1,N]}_{3q+1})(\prod_{i=1}^{N}\psi(dy^{i}_{3q+2})\\
    &&h(x_{3q+2},y^{i}_{3q+2})f^{i}_{n,3q+2}(du^{i}_{3q+2}|y^{i}_{[3q,3q+2]},u^{i}_{[3q,3q+1]}))-\int g_{m}  \prod_{i=1}^{N}\psi(dy^{i}_{3q})h(x_{3q},y^{i}_{3q})\\
    &&f^{i}_{3q,n}(du^{i}_{3q}|y_{3q}^{i})\tau(dx_{3q+1}|x_{3q},u^{[1,N]}_{3q})(\prod_{i=1}^{N}\psi(dy^{i}_{3q+1})h(x_{3q+1},y^{i}_{3q+1})\\   &&f^{i}_{3q+1,n}(du^{i}_{3q+1}|y_{[3q,3q+1]}^{i},u^{i}_{3q}))\tau(dx_{3q+2}|x_{3q+1},u^{[1,N]}_{3q+1})(\prod_{i=1}^{N}\psi(dy^{i}_{3q+2})\\
    &&h(x_{3q+2},y^{i}_{3q+2})f^{i}_{3q+2}(du^{i}_{3q+2}|y^{i}_{[3q,3q+2]},u^{i}_{[3q,3q+1]}))\bigg|
\end{eqnarray*}
\begin{eqnarray*}
    && M_{2}^{n}= \bigg|\int g_{m} \prod_{i=1}^{N}\psi(dy^{i}_{3q})h(x_{3q},y^{i}_{3q})f^{i}_{n,3q}(du^{i}_{3q}|y_{3q}^{i})\tau(dx_{3q+1}|x_{3q},u^{[1,N]}_{3q})\\
    &&(\prod_{i=1}^{N}\psi(dy^{i}_{3q+1})h(x_{3q+1},y^{i}_{3q+1})f^{i}_{n,3q+1}(du^{i}_{3q+1}|y_{[3q,3q+1]}^{i},u^{i}_{3q}))\\
    &&\tau(dx_{3q+2}|x_{3q+1},u^{[1,N]}_{3q+1})(\prod_{i=1}^{N}\psi(dy^{i}_{3q+2})h(x_{3q+2},y^{i}_{3q+2})f^{i}_{3q+2}(du^{i}_{3q+2}|y^{i}_{[3q,3q+2]},u^{i}_{[3q,3q+1]}))\\
    &&-\int g_{m} \prod_{i=1}^{N}\psi(dy^{i}_{3q})h(x_{3q},y^{i}_{3q})f^{i}_{3q,n}(du^{i}_{3q}|y_{3q}^{i})\tau(dx_{3q+1}|x_{3q},u^{[1,N]}_{3q})\\
    &&(\prod_{i=1}^{N}\psi(dy^{i}_{3q+1})h(x_{3q+1},y^{i}_{3q+1})f^{i}_{3q+1}(du^{i}_{3q+1}|y_{[3q,3q+1]}^{i},u^{i}_{3q}))\tau(dx_{3q+2}|x_{3q+1},u^{[1,N]}_{3q+1})\\
    &&(\prod_{i=1}^{N}\psi(dy^{i}_{3q+2})h(x_{3q+2},y^{i}_{3q+2})f^{i}_{3q+2}(du^{i}_{3q+2}|y^{i}_{[3q,3q+2]},u^{i}_{[3q,3q+1]}))\bigg|   
\end{eqnarray*}

and,

\begin{eqnarray*}
    && M_{3}^{n}=  \bigg|\int g_{m} \prod_{i=1}^{N}\psi(dy^{i}_{3q})h(x_{3q},y^{i}_{3q})f^{i}_{n,3q}(du^{i}_{3q}|y_{3q}^{i})\tau(dx_{3q+1}|x_{3q},u^{[1,N]}_{3q})\\
    &&(\prod_{i=1}^{N}\psi(dy^{i}_{3q+1})h(x_{3q+1},y^{i}_{3q+1})f^{i}_{3q+1}(du^{i}_{3q+1}|y_{[3q,3q+1]}^{i},u^{i}_{3q}))\tau(dx_{3q+2}|x_{3q+1},u^{[1,N]}_{3q+1})\\
    &&(\prod_{i=1}^{N}\psi(dy^{i}_{3q+2})h(x_{3q+2},y^{i}_{3q+2})f^{i}_{3q+2}(du^{i}_{3q+2}|y^{i}_{[3q,3q+2]},u^{i}_{[3q,3q+1]}))-\\
    &&\int g_{m} \prod_{i=1}^{N}\psi(dy^{i}_{3q})h(x_{3q},y^{i}_{3q})f^{i}_{3q}(du^{i}_{3q}|y_{3q}^{i})\tau(dx_{3q+1}|x_{3q},u^{[1,N]}_{3q})\\
    &&(\prod_{i=1}^{N}\psi(dy^{i}_{3q+1})h(x_{3q+1},y^{i}_{3q+1})f^{i}_{3q+1}(du^{i}_{3q+1}|y_{[3q,3q+1]}^{i},u^{i}_{3q}))\tau(dx_{3q+2}|x_{3q+1},u^{[1,N]}_{3q+1})\\
    &&(\prod_{i=1}^{N}\psi(dy^{i}_{3q+2})h(x_{3q+2},y^{i}_{3q+2})f^{i}_{3q+2}(du^{i}_{3q+2}|y^{i}_{[3q,3q+2]},u^{i}_{[3q,3q+1]}))\bigg|   
\end{eqnarray*}
To show that $\sup_{x_{qK}}M_{1}^{n}\rightarrow 0$, it is sufficient to check that for all $\mathbf{u}_{[qk,qk+2]}$: $\sup_{x_{qK}}M_{1}^{n}(\mathbf{u}_{[qk,qk+2]})\rightarrow 0$. Where,\\
\begin{eqnarray*}
    && M_{1}^{n}(\mathbf{u}_{[qk,qk+2]}):=  \bigg|\int g_{m}\prod_{i=1}^{N}\psi(dy^{i}_{3q})h(x_{3q},y^{i}_{3q})f^{i}_{n,3q}(du^{i}_{3q}|y_{3q}^{i})\tau(dx_{3q+1}|x_{3q},u^{[1,N]}_{3q})\\
    &&(\prod_{i=1}^{N}\psi(dy^{i}_{3q+1})h(x_{3q+1},y^{i}_{3q+1})f^{i}_{n,3q+1}(du^{i}_{3q+1}|y_{[3q,3q+1]}^{i},u^{i}_{3q}))\tau(dx_{3q+2}|x_{3q+1},u^{[1,N]}_{3q+1})\\
    &&(\prod_{i=1}^{N}\psi(dy^{i}_{3q+2})h(x_{3q+2},y^{i}_{3q+2})f^{i}_{n,3q+2}(du^{i}_{3q+2}|y^{i}_{[3q,3q+2]},u^{i}_{[3q,3q+1]}))-\\
    &&\int g_{m}  \prod_{i=1}^{N}\psi(dy^{i}_{3q})h(x_{3q},y^{i}_{3q})f^{i}_{3q,n}(du^{i}_{3q}|y_{3q}^{i})\tau(dx_{3q+1}|x_{3q},u^{[1,N]}_{3q})\\
    &&(\prod_{i=1}^{N}\psi(dy^{i}_{3q+1})h(x_{3q+1},y^{i}_{3q+1})f^{i}_{3q+1,n}(du^{i}_{3q+1}|y_{[3q,3q+1]}^{i},u^{i}_{3q}))\tau(dx_{3q+2}|x_{3q+1},u^{[1,N]}_{3q+1})\\
    &&(\prod_{i=1}^{N}\psi(dy^{i}_{3q+2})h(x_{3q+2},y^{i}_{3q+2})f^{i}_{3q+2}(du^{i}_{3q+2}|y^{i}_{[3q,3q+2]},u^{i}_{[3q,3q+1]}))\bigg|\\
\end{eqnarray*}
Note that in the expression above there is no integration with respect to the action spaces. We have that
\begin{eqnarray*} 
   && M_{1}^{n}(\mathbf{u}_{[qk,qk+2]})= \bigg|\int (\prod_{i=1}^{N}\psi(dy^{i}_{3q})\psi(dy^{i}_{3q+1})\psi(dy^{i}_{3q+2}))\\
   &&f^{i}_{n,3q+2}(du^{i}_{3q+2}|y^{i}_{[3q,3q+2]},u^{i}_{[3q,3q+1]})\int g_{m} \prod_{i=1}^{N}h(x_{3q},y^{i}_{3q})f^{i}_{n,3q}(du^{i}_{3q}|y_{3q}^{i})\\
   &&\tau(dx_{3q+1}|x_{3q},u^{[1,N]}_{3q})(\prod_{i=1}^{N}h(x_{3q+1},y^{i}_{3q+1})f^{i}_{n,3q+1}(du^{i}_{3q+1}|y_{[3q,3q+1]}^{i},u^{i}_{3q}))\\  &&\tau(dx_{3q+2}|x_{3q+1},u^{[1,N]}_{3q+1})\prod_{i=1}^{N}h(x_{3q+2},y^{i}_{3q+2})-\int (\prod_{i=1}^{N}\psi(dy^{i}_{3q})\psi(dy^{i}_{3q+1})\psi(dy^{i}_{3q+2}))\\
   &&f^{i}_{3q+2}(du^{i}_{3q+2}|y^{i}_{[3q,3q+2]},u^{i}_{[3q,3q+1]})\int g_{m} \prod_{i=1}^{N}h(x_{3q},y^{i}_{3q})f^{i}_{n,3q}(du^{i}_{3q}|y_{3q}^{i})\\
   &&\tau(dx_{3q+1}|x_{3q},u^{[1,N]}_{3q})(\prod_{i=1}^{N}h(x_{3q+1},y^{i}_{3q+1})f^{i}_{n,3q+1}(du^{i}_{3q+1}|y_{[3q,3q+1]}^{i},u^{i}_{3q}))\\  &&\tau(dx_{3q+2}|x_{3q+1},u^{[1,N]}_{3q+1})\prod_{i=1}^{N}h(x_{3q+2},y^{i}_{3q+2})\bigg|  \\   
\end{eqnarray*}
Let $\epsilon >0$. Since the functions $h, f_{n,3q}, f_{n,3q+1}, f_{n,3q+2}, f_{3q}, f_{3q+1}, f_{3q+2}$ are all measurable in $\mathbf{y}$ there exists a compact set $K$ such that $\psi(K^{C})<\epsilon$ and all those functions are continuous on $K$. Where, here, we abuse notation to denote by $\psi$ the measure induced by the latter on the product space $\prod_{i=1}^{N}\mathbf{Y}^{i}$. Without loss of generality $K$ can be chosen so that continuity holds for all $n$. Thus, we get
\begin{eqnarray*}
     && M_{1}^{n}(\mathbf{u}_{[qk,qk+2]})\leq \bigg|\int_{K} (\prod_{i=1}^{N}\psi(dy^{i}_{3q})\psi(dy^{i}_{3q+1})\psi(dy^{i}_{3q+2}))\\
     &&f^{i}_{n,3q+2}(du^{i}_{3q+2}|y^{i}_{[3q,3q+2]},u^{i}_{[3q,3q+1]})\int g_{m} \prod_{i=1}^{N}h(x_{3q},y^{i}_{3q})f^{i}_{n,3q}(du^{i}_{3q}|y_{3q}^{i})\\
     &&\tau(dx_{3q+1}|x_{3q},u^{[1,N]}_{3q})(\prod_{i=1}^{N}h(x_{3q+1},y^{i}_{3q+1})f^{i}_{n,3q+1}(du^{i}_{3q+1}|y_{[3q,3q+1]}^{i},u^{i}_{3q}))\\  &&\tau(dx_{3q+2}|x_{3q+1},u^{[1,N]}_{3q+1})\prod_{i=1}^{N}h(x_{3q+2},y^{i}_{3q+2})-\int_{K} (\prod_{i=1}^{N}\psi(dy^{i}_{3q})\psi(dy^{i}_{3q+1})\psi(dy^{i}_{3q+2}))\\
     &&f^{i}_{3q+2}(du^{i}_{3q+2}|y^{i}_{[3q,3q+2]},u^{i}_{[3q,3q+1]})\int g_{m} \prod_{i=1}^{N}h(x_{3q},y^{i}_{3q})f^{i}_{n,3q}(du^{i}_{3q}|y_{3q}^{i})\\
     &&\tau(dx_{3q+1}|x_{3q},u^{[1,N]}_{3q})(\prod_{i=1}^{N}h(x_{3q+1},y^{i}_{3q+1})f^{i}_{n,3q+1}(du^{i}_{3q+1}|y_{[3q,3q+1]}^{i},u^{i}_{3q}))\\
     &&\tau(dx_{3q+2}|x_{3q+1},u^{[1,N]}_{3q+1})\prod_{i=1}^{N}h(x_{3q+2},y^{i}_{3q+2})\bigg| +2\epsilon\|g_{m}\|_{\infty}  \\
\end{eqnarray*}
Now, over $K$ for all $n$ the function \\
\begin{eqnarray*}
    &&\int g_{m} \prod_{i=1}^{N}h(x_{3q},y^{i}_{3q})f^{i}_{n,3q}(du^{i}_{3q}|y_{3q}^{i})\tau(dx_{3q+1}|x_{3q},u^{[1,N]}_{3q})(\prod_{i=1}^{N}h(x_{3q+1},y^{i}_{3q+1})\\
    &&f^{i}_{n,3q+1}(du^{i}_{3q+1}|y_{[3q,3q+1]}^{i},u^{i}_{3q}))\tau(dx_{3q+2}|x_{3q+1},u^{[1,N]}_{3q+1})\prod_{i=1}^{N}h(x_{3q+2},y^{i}_{3q+2}):=\hat{g}_{n}
\end{eqnarray*}
is continuous. Now, one can define a metric on probability measures by considering a countable sequence of continuous and bounded functions $\tilde{g}_{n}$ which includes all $\hat{g}_{n}$ and such that $\rho_{1}(\nu_{1},\nu_{2})=\sum_{n=0}^{\infty} 2^{-n}|\int \tilde{g}_{n}\nu_{1}-\int \tilde{g}_{n}\nu_{1}|$. Since $\rho_{1}$ induces the weak topology there exists some constant $C$ such that for all $\nu_{1},\nu_{2}$ $\rho_{1}(\nu_{1},\nu_{2})\leq C\rho(\nu_{1},\nu_{2})$. Hence, we get that for all $x_{qK}$
\begin{eqnarray*}
   && M_{1}^{n}(\mathbf{u}_{[qk,qk+2]})\leq 2\epsilon\|g_{m}\|_{\infty} +C \sup_{\|g\|_{BL}\leq 1} \bigg|\int g (\prod_{i=1}^{N}\psi(dy^{i}_{3q})\psi(dy^{i}_{3q+1})\psi(dy^{i}_{3q+2}))\\
 \nonumber    &&f^{i}_{n,3q+2}(du^{i}_{3q+2}|y^{i}_{[3q,3q+2]},u^{i}_{[3q,3q+1]})-\int g(\prod_{i=1}^{N}\psi(dy^{i}_{3q})\psi(dy^{i}_{3q+1})\psi(dy^{i}_{3q+2}))\\
    &&f^{i}_{3q+2}(du^{i}_{3q+2}|y^{i}_{[3q,3q+2]},u^{i}_{[3q,3q+1]})\bigg|
\end{eqnarray*}
Thus, it is sufficient to check that
\begin{eqnarray}\label{Need to check I} 
   \nonumber  && \sup_{\|g\|_{BL}\leq 1} \bigg|\int g (\prod_{i=1}^{N}\psi(dy^{i}_{3q})\psi(dy^{i}_{3q+1})\psi(dy^{i}_{3q+2}))\\
 \nonumber    &&f^{i}_{n,3q+2}(du^{i}_{3q+2}|y^{i}_{[3q,3q+2]},u^{i}_{[3q,3q+1]})-\int g(\prod_{i=1}^{N}\psi(dy^{i}_{3q})\psi(dy^{i}_{3q+1})\psi(dy^{i}_{3q+2}))\\
    &&f^{i}_{3q+2}(du^{i}_{3q+2}|y^{i}_{[3q,3q+2]},u^{i}_{[3q,3q+1]})\bigg|\rightarrow 0 
\end{eqnarray}
Let
\begin{eqnarray*}
    \mu_{n}(.)&=&\int_{.}  (\prod_{i=1}^{N}\psi(dy^{i}_{3q})\psi(dy^{i}_{3q+1})\psi(dy^{i}_{3q+2}))f^{i}_{n,3q+2}(du^{i}_{3q+2}|y^{i}_{[3q,3q+2]},u^{i}_{[3q,3q+1]})\\
    \mu(.) &=& \int_{.}  (\prod_{i=1}^{N}\psi(dy^{i}_{3q})\psi(dy^{i}_{3q+1})\psi(dy^{i}_{3q+2}))f^{i}_{3q+2}(du^{i}_{3q+2}|y^{i}_{[3q,3q+2]},u^{i}_{[3q,3q+1]})
\end{eqnarray*}
We have that $\rho(\mu_{n},\mu)\rightarrow 0$. Hence the sequence constructed from the restriction of these measures to sigma field of subsets of $K$ also converges weakly. Thus \ref{Need to check I} holds, which in turn implies that $M_{1}^{n}\rightarrow0$. The proof that $M_{2}^{n}\rightarrow 0$ and $M_{3}^{n}\rightarrow 0$ proceeds in a similar manner. 
\qed
\begin{corollary} \label{existence KSDISP}
    Suppose the transition kernel $\tau\left(x_{t+1} \in . \mid x_{t}=x, \mathbf{u}_{t}=\mathbf{u}\right)$ is continuous under total variation in $x$ and $u$, and $h(x, y)$ is continuous in $x$. Then, an optimal solution to the KSPISP problem exists and stationary policies are optimal.
\end{corollary}
\textbf{Proof.} See the Appendix. 
\subsection{Examples}
Here, we will present various examples which showcase the kernel regularity conditions needed for our weak-Feller and existence results. For centralized Markov decision processes various example leading to different regularity conditions are presented in \cite{KSYContQLearning}. In the following, we will build on \cite{kara2020robustness}. Recall (\ref{Prob 1}).

\begin{example} Suppose the realization function $\Phi(x,u^{1},\cdots,u^{N},v)$ is continuous in $(x,u^{1},\cdots,u^{N})$. Then, the kernel $\tau$ is weak-Feller (weakly continuous), i.e., for any sequence $\{(x_{n},u^{1}_{n},\cdots,u^{N}_{n})\}_{n\in \mathbb{N}}$ such that $(x_{n},u^{1}_{n},\cdots,u^{N}_{n})\rightarrow (x,u^{1},\cdots,u^{N})$ as $n\rightarrow \infty$, we have that for any continuous and bounded function $g\in \mathbf{C_{b}}(\mathbb{X})$ $\int g(x_{1})\tau(dx_{1}|x_{n},u^{1}_{n},\cdots,u^{N}_{n})\rightarrow \int g(x_{1})\tau(dx_{1}|x,u^{1},\cdots,u^{N})$ as $n\rightarrow \infty$. To see this, let $g\in \mathbf{C_{b}}(\mathbb{X})$. Then, by an application of the dominated convergence theorem (DCT), we have that
\begin{eqnarray*}
    &&\int g(x_{1})\tau(dx_{1}|x_{n},u^{1}_{n},\cdots,u^{N}_{n})=\int g(\Phi(x_{n},u^{1}_{n},\cdots,u^{N}_{n},v))P(dv)\\
    &&\rightarrow\int g(\Phi(x,u^{1},\cdots,u^{N},v))P(dv)=\int g(x_{1})\tau(dx_{1}|x,u^{1},\cdots,u^{N})
\end{eqnarray*}
Thus, this model satisfies the continuity condition on the kernel as needed for Theorem \ref{theorem5}. 
\end{example}

\begin{example} Suppose the kernel $\tau$ admits a density $\tilde{\Phi}$ with respect to some measure $\xi\in \mathcal{P}(\mathbb{X})$ so that $\tau(x_{1}\in\cdot|x,u^{1},\cdots,u^{N})=\int_{\cdot}\tilde{\Phi}(x_{1},x,u^{1},\cdots,u^{N})\xi(dx_{1})$. If $\tilde{\Phi}$ is continuous in $(x,u^{1},\cdots,u^{N})$, then $\tau$ is continuous in total variation, since for any sequence $\{(x_{n},u^{1}_{n},\cdots,u^{N}_{n})\}_{n\in \mathbb{N}}$ such that $(x_{n},u^{1}_{n},\cdots,u^{N}_{n})\rightarrow (x,u^{1},\cdots,u^{N})$ as $n\rightarrow \infty$, we have that $\displaystyle  \sup_{\|g\|_{\infty}\leq 1}\bigg|\int g(x_{1})\tau(dx_{1}|x_{n},u^{1}_{n},\cdots,u^{N}_{n})-\int g(x_{1})\tau(dx_{1}|x,u^{1},\cdots,u^{N})\bigg|\rightarrow 0$ as $n\rightarrow \infty$. This follows from an application of Scheff\'e Lemma \cite[Theorem 16.12]{Bil95}. 
This then serves as an example for the conditions on the kernels for both Theorems \ref{theorem5} and \ref{theorem6}.
\end{example}

\begin{example} Suppose the realization function $\Phi$ can be written as $\Phi(x,u^{1},\cdots,u^{N},v)=\tilde{\Phi}(x,u^{1},\cdots,u^{N})+v$ where $\tilde{\Phi}$ is continuous in $(x,u^{1},\cdots,u^{N})$ and $v$ admits a continuous density function $\iota$ with respect to some measure $\xi\in \mathcal{P}([0,1])$. Then, $\tau$ is continuous in total variation. This then serves as an example for the conditions on the kernels for both Theorems \ref{theorem5} and \ref{theorem6}. Similarly, one can construct examples for which regularity of the observation channels hold.
\end{example}

\subsection{Finite approximations of the equivalent centralized MDPs}
The centralized reduction and weak-Feller properties presented in this paper allow one to establish the existence of solutions for the one-step delayed information sharing pattern and $K$-step periodic information sharing pattern problems. Additionally, one can use dynamic programming or the discounted cost optimality equation (DCOE) to find solutions to these problems. Moreover, the DCOE also allows one to establish the optimality of stationary policies. However, a challenge arises due to the fact that numerical learning techniques can not be applied to MDPs for which the action and state spaces are not finite. An important method for dealing with such a challenge is through discretizing, also known as quantizing, the state and action spaces in a manner that leads to a finite state and action space MDP which approximates the original MDP arbitrarily well as the size of the quantization increases \cite{SaYuLi15c} (See also \cite{KSYContQLearning}). In this section we will focus on the MDP $(\pi_{q},a_{q})$. We will first describe the procedure for quantizing the action space which we denote $A$. Then, we will describe the procedure for quantizing the state space. A similar reasoning can be applied to the MDP $(Z_{t},f_{t})$ as well. 
\subsubsection{Action space quantization in the equivalent model}
 Here, we describe a quantization procedure which was first proposed in \cite{SaYuLi15c}. Note that under the topology introduced in Definition \ref{top k-step}, we have that $\mathbb{A}$ is compact. Let $n \in \mathbb{N}$ and $\mathbb{A}^{n}=\{a_{1},a_{2},...,a_{n}\}$ be such that for all $a\in\mathbb{A}$ there exists $i\in \{1,...,n\}$ such that $d_{\mathbb{A}}(a,a_{i})<\frac{1}{n}$ where $d_{\mathbb{A}}$ is the metric on the space $\mathbb{A}$ which stems from embedding local policies in the space of probability measures with a fixed marginal and endowing the latter with the bounded Lipschitz metric as in Definition \ref{top k-step}. Consider an MDP, which we denote $MDP_{n}=(\mathcal{P}(\mathbb{X}),\mathbb{A}^{n},\eta,\tilde{c})$ such that the action space $\mathbb{A}$ is replaced with $\mathbb{A}^n$. We denote the optimal cost for $MDP_{n}$ by $J^{*}_{n}$. 
 \begin{assumption} \label{assumption for quantizing actions} The following hold: 
         (i) $\eta$ (as defined by equation \ref{theta definition}) is weak-Feller . (ii) $\tilde{c}$ (as defined by equation \ref{New k-step Cost}) is bounded .  
 \end{assumption}
 Note that Assumption (ii) holds since $c$ is bounded. Whereas Assumption (i) holds under the conditions of Theorem \ref{theorem6}.
\begin{theorem} \cite[Theorem 3.16 ]{SaLiYuSpringer}
Suppose Assumption \ref{assumption for quantizing actions} holds. Then, for any $\pi\in \mathcal{P}(\mathbb{X})$: $|J^{*}_{n}(\pi)-J^{*}(\pi)|\rightarrow 0$ as $n\rightarrow \infty$.  
\end{theorem}
\subsubsection{State space quantization in the equivalent model}
Next, we will describe a procedure for quantizing the (effective) state space. Here, again the compactness of $\mathbb{A}$ is crucial. Let $n \in \mathbb{N}$ and $\mathcal{P}(\mathbb{X})^{n}=\{\pi_{1},\pi_{2},...,\pi_{n}\}$ be such that for all $\pi\in\mathcal{P}(\mathbb{X})$ there exists $i\in \{1,...,n\}$ such that $W_{1}(\pi,\pi_{i})<\frac{1}{n}$. Let $\nu \in \mathcal{P}(\mathcal{P}(\mathbb{X}))$ be a measure such that $\forall i \in\{1,...,n\}$ $\nu(B_{i})>0$ where $B_{i}=\{\pi\in \mathcal{P}(\mathbb{X})| Q(\pi)=\pi_{i}\}$. Here, $Q$ is the quantizer that maps every element $\pi$ to the nearest element in $\mathcal{P}(\mathbb{X})^{n}$ under $W_1$. If ties occur, it is assumed that they are broken so that the function $Q$ is measurable. Consider the following transition kernel: $P(\pi_{j}|\pi_{i},a)=\int_{B_{i}}\int_{\mathcal{P}(\mathbb{X})} \mathbbm{1}_{\{Q(s)=\pi_{j}\}}\eta(ds|w,a)\nu_{i}(dw)$. Here $\nu_{i}$ denotes the normalized weight measure of $\nu$ on the set $B_{i}$ that it is to say $\nu_{i}(E)=\frac{\nu(E)}{\nu(B_{i})}$ for all $E\subseteq B_{i}$ which are measurable. Consider the MDP given by $MDP_{n}=(\mathcal{P}(\mathbb{X})^{n},\mathbb{A},P,\tilde{c})$.  We denote by $J^{*}_{n}$ the expected cost accrued when the optimal policy for $MDP_{n}$ is extended into a policy for the original MDP which is constant over each bin. 

\begin{theorem} \cite[Theorem 4.3]{SaLiYuSpringer} \label{theorem sate space quantization}
    Suppose Assumption \ref{assumption for quantizing actions} holds. Then, $|J^{*}_{n}(\pi)-J^{*}(\pi)|\rightarrow 0$ as $n\rightarrow \infty$.   
\end{theorem}

Quantization of the state and action spaces is a very useful tool that allows one to use numerical learning techniques such as Q-learning and empirical value iteration \cite{Haskel2014empiricalvalueiteration}. Numerical and empirical results will be reported in future studies. 


\section{Conclusion}
In this paper, we have established that both one-step delayed information sharing pattern and $K$-step periodic information sharing pattern problems admit a centralized reduction, their optimal policies are separable, and the kernels of their centralized reductions are weak-Feller. The combination of these results, allow one to arrive at the existence of solutions (Through the DCOE, combined with the compactness of the space of actions under Young topology, as well the continuity and boundedness of $\tilde{c}$), apply analytical tools to find those solutions, and establish the optimality of stationary policies. Moreover, through the use of state and action space quantization, one can use numerical techniques to arrive at solutions with near optimal performance guarantees. Additionally, through approximation with KSPISP problems, one can also find near optimal solutions for completely decentralized information structure problems.  

\section{Appendix}
\subsection{Proof of Proposition \ref{finite action approximation}} \label{Proof of Proposition on finite action approximation}
Here, we will use the following notation for the measurement channel $Q(d\mathbf{y}|x)=\prod_{i=1}^{N}Q^{i}(dy^{i}|x)$. Because the cost functional $c$ is bounded, it is sufficient to show that the result holds for a finite horizon cost criteria $\displaystyle J^{Tf}(\gamma)=\sum_{t=0}^{T_{f}} \beta^{t} E^{\gamma}[c(x_{t},\mathbf{u}_{t})]$. In what follows we assume that $T_{f}=3$. The more general proof proceeds in a similar manner. \\

Step 1: Let $\epsilon >0$. By Stone-Weierstrass theorem, since $\mathbb{X}$ and $\mathbb{U}$ are compact, there exists a Lipschitz function $L(x,u)$ with coefficient $\alpha$ such that $\|c(x,u)-L(x,u)\|_{\infty}\leq \epsilon$. Additionally, because $\tau:\mathbb{X}\times \mathbb{U}\rightarrow \mathcal{P}(\mathbb{X})$ is continuous, under total variation, over a compact set the exists $\delta>0$ such that for any $(x,\mathbf{u}),(x_{1},\mathbf{u}_{1})\in \mathbb{X}\times \mathbb{U}$ such that $d\big((x,\mathbf{u}),(x_{1},\mathbf{u}_{1})\big)\leq \delta$ we have that $\|\tau(.|x_{1},\mathbf{u}_{1})-\tau(.|x,\mathbf{u})\|_{TV}\leq \epsilon$. \\

Step 2: Let $\gamma=(\gamma_{0},\gamma_{1},\gamma_{2})$ be an optimal team policy. Let $\gamma_{m}=(\gamma_{0}^{m},\gamma_{1}^{m},\gamma_{2}^{m})$ be a team policy such that for any time $t$ and realization of $I^{1}_{t},...,I^{N}_{t}$ we have that $d(\gamma_{t}(I^{1}_{t},...,I^{N}_{t}),\gamma^{m}_{t}(I^{1}_{t},...,I^{N}_{t}))\leq \frac{1}{m}$. Then, we get that $|J^{Tf}(\gamma_{0},\gamma_{1},\gamma_{2})-J^{Tf}(\gamma_{0}^{m},\gamma_{1}^{m},\gamma_{2}^{m})|\leq|J^{Tf}(\gamma_{0},\gamma_{1},\gamma_{2})-J^{Tf}(\gamma_{0},\gamma_{1},\gamma_{2}^{m})|+ |J^{Tf}(\gamma_{0},\gamma_{1},\gamma_{2}^{m})-J^{Tf}(\gamma_{0},\gamma_{1}^{m},\gamma_{2}^{m})|+|J^{Tf}(\gamma_{0},\gamma_{1}^{m},\gamma_{2}^{m})-J^{Tf}(\gamma_{0}^{m},\gamma_{1}^{m},\gamma_{2}^{m})|$ \\

Step 3: We have that 
\begin{eqnarray*}
\nonumber    &&|J^{Tf}(\gamma_{0},\gamma_{1},\gamma_{2})-J^{Tf}(\gamma_{0},\gamma_{1},\gamma_{2}^{m})|\leq 2\beta^{2}\epsilon +\beta^{2}\bigg|\int \mu(dx_{0})Q(d\mathbf{y}_{0}|x_{0})\\
\nonumber    &&\tau(dx_{1}|x_{0},\gamma_{0}(I^{1}_{0},...,I^{N}_{0}))Q(d\mathbf{y}_{1}|x_{1})\tau(dx_{2}|x_{1},\gamma_{1}(I^{1}_{1},...,I^{N}_{1})) Q(d\mathbf{y}_{2}|x_{2})\\
\nonumber    &&L(x_{2},\gamma_{2}(I^{1}_{2},...,I^{N}_{2}))-\int \mu(dx_{0})Q(d\mathbf{y}_{0}|x_{0})\tau(dx_{1}|x_{0},\gamma_{0}(I^{1}_{0},...,I^{N}_{0}))Q(d\mathbf{y}_{1}|x_{1})\\
\nonumber    &&\tau(dx_{2}|x_{1},\gamma_{1}(I^{1}_{1},...,I^{N}_{1})) Q(d\mathbf{y}_{2}|x_{2})L(x_{2},\gamma_{2}^{m}(I^{1}_{2},...,I^{N}_{2}))\bigg|\leq \beta^{2}[2\epsilon+ \frac{\alpha}{m}]
\end{eqnarray*}
Next, by applying a similar argument to the one above when $t=1$, we get that for all $m$ such that $\frac{1}{m}\leq \delta$
\begin{eqnarray*}
\nonumber    &&|J^{Tf}(\gamma_{0},\gamma_{1},\gamma_{2}^{m})-J^{Tf}(\gamma_{0},\gamma_{1}^{m},\gamma_{2}^{m})|\leq \beta[2\epsilon+ \frac{\alpha}{m}]\\
&& \quad +\beta^{2}\bigg|\int \mu(dx_{0})Q(d\mathbf{y}_{0}|x_{0})\tau(dx_{1}|x_{0},\gamma_{0}(I^{1}_{0},...,I^{N}_{0}))Q(d\mathbf{y}_{1}|x_{1})\\
\nonumber &&\int \tau(dx_{2}|x_{1},\gamma_{1}(I^{1}_{1},...,I^{N}_{1})) Q(d\mathbf{y}_{2}|x_{2})c(x_{2},\gamma_{2}^{m}(I^{1}_{2},...,I^{N}_{2})) \\
&& \quad -\int \mu(dx_{0})Q(d\mathbf{y}_{0}|x_{0})\tau(dx_{1}|x_{0},\gamma_{0}(I^{1}_{0},...,I^{N}_{0}))Q(d\mathbf{y}_{1}|x_{1})\int\tau(dx_{2}|x_{1},\gamma_{1}^{m}(I^{1}_{1},...,I^{N}_{1}))\\
&&Q(d\mathbf{y}_{2}|x_{2})c(x_{2},\gamma_{2}^{m}(I^{1}_{2},...,I^{N}_{2}))\bigg|\\
\nonumber    &&\leq \beta[2\epsilon+ \frac{\alpha}{m}]+\beta^{2}\|c\|_{\infty}\int \mu(dx_{0})Q(d\mathbf{y}_{0}|x_{0})\tau(dx_{1}|x_{0},\gamma_{0}(I^{1}_{0},...,I^{N}_{0}))Q(d\mathbf{y}_{1}|x_{1})\\
\nonumber    &&\|\tau(.|x_{1},\gamma_{1}(I^{1}_{1},...,I^{N}_{1}))-\tau(.|x_{1},\gamma_{1}^{m}(I^{1}_{1},...,I^{N}_{1}))\|_{TV}\leq\beta[2\epsilon+ \frac{\alpha}{m}]+\beta^{2}\|c\|_{\infty}\epsilon 
\end{eqnarray*}
Similarly, we get that
\begin{eqnarray*}
\nonumber    &&|J^{Tf}(\gamma_{0},\gamma_{1}^{m},\gamma_{2}^{m})-J^{Tf}(\gamma_{0}^{m},\gamma_{1}^{m},\gamma_{2}^{m})|\leq 2\epsilon+\frac{\alpha}{m}+\|c\|_{\infty}[\beta+\beta^{2}]\epsilon
\end{eqnarray*}
Step 4: Hence, we get that 
\begin{eqnarray*}
\nonumber    &&|J^{Tf}(\gamma_{0},\gamma_{1},\gamma_{2})-J^{Tf}(\gamma_{0}^{m},\gamma_{1}^{m},\gamma_{2}^{m})|\leq (1+\beta+\beta^{2})(2\epsilon+\frac{\alpha}{m})\\
 \nonumber   &&+\quad \|c\|_{\infty}[\beta+2\beta^{2}]\epsilon \rightarrow (1+\beta+\beta^{2})(2\epsilon)+\|c\|_{\infty}[\beta+2\beta^{2}]\epsilon \text{ as } m\rightarrow \infty
\end{eqnarray*}
Since $\epsilon$ was arbitrary this completes the proof.\\
\qed

\subsection{Proof of Theorem \ref{theorem5}.} 
Here, we will need the following lemma which is a weaker of version of \cite[Lemma 1]{KSYWeakFellerSysCont}.

\begin{lemma} \label{lemma weak feller one step}Let $\{G_{A}(x)\}_{A}$ be a family of real-valued functions from some standard Borel space $\mathbb{X}$ which are indexed over some set $A$ and uniformly bounded (for all $A$, $||G_{A}(x)||_{\infty}\leq C$) such that for every $\{x^n\}_n$ such that $x^{n}\rightarrow x$ we have $\sup_{A}|G_{A}(x^{n})-G_{A}(x)|\rightarrow 0$, as $n \rightarrow \infty$. Then for any sequence of probability measures $\{Z_{n}\}_{n}$ such that $Z_{n}\rightarrow Z$ weakly, we have that: $\displaystyle \sup_{A}|\int G_{A}(x)Z_{n}(dx)-G_{A}(x)Z(dx)|\rightarrow 0$.
\end{lemma}

Let $g\in\mathbf{C_{b}}(\mathcal{P}(\mathbb{X}))$ be a continuous and bounded function on $\mathcal{P}(\mathbb{X})$. We have
\begin{align*}
   &\bigg{|} E[g(z_{t+1})|Z^{n}_{t},f_{n,t}]-E[g(z_{t+1})|Z_{t},f_{t}]           \bigg{|}  \quad \\ 
   &=  \bigg|\int g(z_{t+1})\mathbbm{1}_{\{F(Z^{n}_{t},\mathbf{u}_{t},y_{t}^{[1,N]}) \in .\}}\prod_{i=1}^{N}f^{i}_{n,t}(du^{i}_{t}|y^{i}_{t})Q^{i}(dy^{i}_{t}|x_{t})Z^{n}_{t}(dx_{t})\\
   &\quad -  \int g(z_{t+1})\mathbbm{1}_{\{F(Z_{t},\mathbf{u}_{t},y_{t}^{[1,N]}) \in .\}}\prod_{i=1}^{N}Q^{i}(dy^{i}_{t}|x_{t})f^{i}_{t}(du^{i}_{t}|y^{i}_{t})Z_{t}(dx_{t})\bigg|\\
   &\leq L_{1}^{n}+L_{2}^{n}
\end{align*}
Where,
\begin{eqnarray*}
  \nonumber && L_{1}^{n}=\bigg|\int g(z_{t+1})\mathbbm{1}_{\{F(Z^{n}_{t},\mathbf{u}_{t},y_{t}^{[1,N]}) \in .\}}\prod_{i=1}^{N}Q^{i}(dy^{i}_{t}|x_{t})f^{i}_{n,t}(du^{i}_{t}|y^{i}_{t})Z^{n}_{t}(dx_{t})- \int g(z_{t+1})\\
  \nonumber &&\mathbbm{1}_{\{F(Z^{n}_{t},\mathbf{u}_{t},y_{t}^{[1,N]}) \in .\}}\prod_{i=1}^{N}Q^{i}(dy^{i}_{t}|x_{t})f^{i}_{n,t}(du^{i}_{t}|y^{i}_{t})Z_{t}(dx_{t})\bigg|
\end{eqnarray*}
and 
\begin{eqnarray*}
 \nonumber  && L_{2}^{n}= \bigg|\int [g(z_{t+1})\mathbbm{1}_{\{F(Z^{n}_{t},\mathbf{u}_{t},y_{t}^{[1,N]}) \in .\}}\prod_{i=1}^{N}f^{i}_{n,t}(du^{i}_{t}|y^{i}_{t})\\
 \nonumber  && \quad -  g(z_{t+1})\mathbbm{1}_{\{F(Z_{t},\mathbf{u}_{t},y_{t}^{[1,N]}) \in .\}}\prod_{i=1}^{N}f^{i}_{t}(du^{i}_{t}|y^{i}_{t})]\prod_{i=1}^{N}Q^{i}(dy^{i}_{t}|x_{t})Z_{t}(dx_{t})\bigg|
\end{eqnarray*}

To complete the proof, we'll show that $L_{1}^{n}\rightarrow 0$ and \\
$L_{2}^{n}\rightarrow 0$ as $n \rightarrow \infty$. \\

Step 1: Showing $L_{1}^{n}\rightarrow 0$. We have:
\begin{eqnarray*}
\nonumber    L_{1}^{n}&\leq& \|g\|_{\infty}\|P(y^{[1,N]}_{t}|Z^{n}_{t})-P(y^{[1,N]}_{t}|Z_{t})\|_{TV}\\
   \nonumber &=&2\|g\|_{\infty}\sup_{A}\bigg|\int_{X}(\int_{A}\prod_{i=1}^{N}Q^{i}(dy^{i}|x_{t}))Z^{n}(dx_{t})-\int_{X}(\int_{A}\prod_{i=1}^{N}Q^{i}(dy^{i}|x_{t}))Z(dx_{t})\bigg|
\end{eqnarray*}

Let $G_{A}(x)=\int_{A}\prod_{i=1}^{N}Q^{i}(dy^{i}|x_{t})$. Since the Radon-Nykodim derivative of the measurement channels $h$ is continuous in $x$, by Scheff\'e's theorem \cite{BillingsleyProbMeasure}, we have that \\
$\displaystyle \sup_{A}|G_{A}(x^n)-G_{A}(x)|\rightarrow 0$. Thus, by Lemma \ref{lemma weak feller one step}, $L_{1}^{n}\rightarrow 0$.\\

Step 2: Showing $L_{2}^{n}\rightarrow 0$. We have that 
\begin{eqnarray*}
 \nonumber  && L_{2}^{n}= \bigg|\int [g(F(Z^{n}_{t},\mathbf{u}_{t},y_{t}^{[1,N]}))\prod_{i=1}^{N}f^{i}_{n,t}(du^{i}_{t}|y^{i}_{t}) -  g(F(Z_{t},\mathbf{u}_{t},y_{t}^{[1,N]}))\prod_{i=1}^{N}f^{i}_{t}(du^{i}_{t}|y^{i}_{t})]\prod_{i=1}^{N}Q^{i}(dy^{i}_{t}|x_{t})Z_{t}(dx_{t})\bigg|
\end{eqnarray*}
 Let $\epsilon>0$ and consider
\begin{eqnarray*}
  \nonumber  &&\mu_{n}(.)=\int_{.}\prod_{i=1}^{N}f^{i}_{n,t}(du^{i}_{t}|y^{i}_{t})\prod_{i=1}^{N}Q^{i}(dy^{i}_{t}|x_{t})Z_{t}(dx_{t})\text{, }   \mu(.)=\int_{.}\prod_{i=1}^{N}f^{i}_{t}(du^{i}_{t}|y^{i}_{t})\prod_{i=1}^{N}Q^{i}(dy^{i}_{t}|x_{t})Z_{t}(dx_{t})
\end{eqnarray*}
Note that $\mu_{n}\rightarrow \mu$ weakly. Hence, there exists a compact set $K_{1}\subseteq \mathbb{X}\times\mathbb{Y}^{1}\times...\times\mathbb{Y}^{N}\times\mathbb{U}^{1}\times...\times\mathbb{U}^{N}$ such that for all $n$ $\mu_{n}(K_{1}^{C})<\epsilon$ and $\mu(K_{1}^{C})<\epsilon$. Thus, we get that 
\begin{eqnarray*}
    &&L_{2}^{n}\leq 2\epsilon \|g\|_{\infty} + \bigg|\int_{K_{1}} [g(F(Z^{n}_{t},\mathbf{u}_{t},y_{t}^{[1,N]}))\prod_{i=1}^{N}f^{i}_{n,t}(du^{i}_{t}|y^{i}_{t}) \\
    && \quad -  g(F(Z_{t},\mathbf{u}_{t},y_{t}^{[1,N]}))\prod_{i=1}^{N}f^{i}_{t}(du^{i}_{t}|y^{i}_{t})]\prod_{i=1}^{N}Q^{i}(dy^{i}_{t}|x_{t})Z_{t}(dx_{t})\bigg|
\end{eqnarray*}
Consider the measure $\displaystyle \nu\in\mathcal{P}(\mathbb{X}\times \prod_{i=1}^{N}\mathbb{Y}^{i})$ given by $\nu(.)=\int_{.} \prod_{i=1}^{N}Q^{i}(dy^{i}_{t}|x_{t})Z_{t}(dx_{t})$. By Lusin's theorem, the measurability of $\displaystyle \prod_{i=1}^{N} h(x,y^{i})$ implies that there exists a closed set $\displaystyle K_{2}\in \mathbb{X}\times \prod_{i=1}^{N}\mathbb{Y}^{i}$ such that the restriction of $\prod_{i=1}^{N} h(x,y^{i})$ to $K_{2}$ is continuous in $\mathbf{y}$ and $\nu(K_{2}^{C})\leq \epsilon$. Hence, by restricting the canonical projection of $K_{1}$ onto $\mathbb{X}\times \prod_{i=1}^{N}\mathbb{Y}^{i}$ to $K_{2}$ we get that there exists a compact set $K$ such that
\begin{eqnarray*}
    &&L_{2}^{n}\leq 4\epsilon\|g\|_{\infty} + \bigg|\int_{K} [g(F(Z^{n}_{t},\mathbf{u}_{t},y_{t}^{[1,N]}))\prod_{i=1}^{N}f^{i}_{n,t}(du^{i}_{t}|y^{i}_{t}) \\
    && \quad-  g(F(Z_{t},\mathbf{u}_{t},y_{t}^{[1,N]}))\prod_{i=1}^{N}f^{i}_{t}(du^{i}_{t}|y^{i}_{t})]\prod_{i=1}^{N}Q^{i}(dy^{i}_{t}|x_{t})Z_{t}(dx_{t})\bigg|
\end{eqnarray*}

Let $L_{3}^{n}=\bigg|\int_{K} [g(F(Z^{n}_{t},\mathbf{u}_{t},y_{t}^{[1,N]}))\prod_{i=1}^{N}f^{i}_{n,t}(du^{i}_{t}|y^{i}_{t})$ \\
$\quad \text{   } -  g(F(Z_{t},\mathbf{u}_{t},y_{t}^{[1,N]}))\prod_{i=1}^{N}f^{i}_{t}(du^{i}_{t}|y^{i}_{t})]\prod_{i=1}^{N}Q^{i}(dy^{i}_{t}|x_{t})Z_{t}(dx_{t})\bigg|$.

Because $\mu_{n}\rightarrow \mu$ weakly, by the generalized DCT theorem \cite[Theorem 3.5]{Lan81} \cite[Theorem 3.5]{serfozo1982convergence} (see also \cite[Theorem D.3.1 (i) ]{lecturenotes}), to show that $L_{3}^{n}\rightarrow 0$ it is sufficient to show that $F(Z^{n}_{t},\mathbf{u}_{t},y_{t}^{[1,N]})\rightarrow F(Z_{t},\mathbf{u}_{t},y_{t}^{[1,N]})$ continuously over $(\mathbf{u}_{t},y_{t}^{[1,N]})$. Let $g_{1}\in \mathbf{C_{b}}(\mathbb{X})$. Consider a sequence such that $(\mathbf{u}^{n}_{t},y_{n,t}^{[1,N]})\rightarrow (\mathbf{u}_{t},y_{t}^{[1,N]})$. Over $K$, we have that $\prod_{i=1}^{N} h(x_{t},y^{i}_{n,t})\rightarrow \prod_{i=1}^{N} h(x_{t},y^{i}_{t})$. Because $h$ is continuous in $x$, the last sequence converges continuously over $x$ and hence by the generalized DCT \cite[Theorem 3.5]{Lan81} \cite[Theorem 3.5]{serfozo1982convergence},\\
$\displaystyle \int Z^{n}_{t}(dx_{t})\prod_{i=1}^{N}h(x_{t},y^{i}_{n,t})\rightarrow \int Z_{t}(dx_{t})\prod_{i=1}^{N}h(x_{t},y^{i}_{t})$. Thus to show that $F(Z^{n}_{t},\mathbf{u}_{t},y_{t}^{[1,N]})\rightarrow F(Z_{t},\mathbf{u}_{t},y_{t}^{[1,N]})$ continuously over $(\mathbf{u}_{t},y_{t}^{[1,N]})$, it is sufficient to show that 
\begin{eqnarray*}
 \nonumber   &&\bigg| \int Z^{n}_{t}(dx_{t})\int g_{1}(x_{t+1})\tau (dx_{t+1}|x_{t},\mathbf{u}_{n})\prod_{i=1}^{N} h(x_{t},y^{i}_{n,t}) \\
 && \quad -  \int Z_{t}(dx_{t})\int g_{1}(x_{t+1})\tau (dx_{t+1}|x_{t},\mathbf{u})\prod_{i=1}^{N} h(x_{t},y^{i}_{t})\bigg|\rightarrow 0
\end{eqnarray*}
Again by the generalized DCT \cite[Theorem 3.5]{Lan81} \cite[Theorem 3.5]{serfozo1982convergence}, it is sufficient to show that for a sequence $x_{n}\rightarrow x$, we have that
\begin{eqnarray*}
  \nonumber  &&\bigg|\int g_{1}(x_{t+1})\tau (dx_{t+1}|x_{n,t},\mathbf{u}_{n})\prod_{i=1}^{N} h(x_{n,t},y^{i}_{n,t})-\int g_{1}(x_{t+1})\tau (dx_{t+1}|x_{t},\mathbf{u})\prod_{i=1}^{N} h(x_{t},y^{i}_{t})\bigg|\rightarrow 0
\end{eqnarray*}
To this end, note that
\begin{eqnarray*}
   \nonumber &&\bigg| \int g_{1}(x_{t+1})\prod_{i=1}^{N} h(x_{n,t},y^{i}_{n,t})\tau (dx_{t+1}|x_{n,t},\mathbf{u}_{n})- \int g_{1}(x_{t+1})\prod_{i=1}^{N} h(x_{t},y^{i}_{t})\tau (dx_{t+1}|x_{t},\mathbf{u})\bigg|\\
   \nonumber &&\leq \bigg| \int g_{1}(x_{t+1})\prod_{i=1}^{N} h(x_{n,t},y^{i}_{n,t})\tau (dx_{t+1}|x_{n,t},\mathbf{u}_{n})-  \int g_{1}(x_{t+1})\prod_{i=1}^{N} h(x_{t},y^{i}_{t})\tau (dx_{t+1}|x_{n,t},\mathbf{u}_{n})\bigg|\\
    \nonumber &&\quad +  \bigg| \int g_{1}(x_{t+1})\prod_{i=1}^{N} h(x_{t},y^{i}_{t})\tau (dx_{t+1}|x_{n,t},\mathbf{u}_{n}) - \int g_{1}(x_{t+1})\prod_{i=1}^{N} h(x_{t},y^{i}_{t})\tau (dx_{t+1}|x_{t},\mathbf{u})\bigg|
\end{eqnarray*}
 Since $K$ is compact, we have that $\displaystyle \prod_{i=1}^{N} h(x_{n,t},y^{i}_{n,t})\rightarrow \prod_{i=1}^{N} h(x_{t},y^{i}_{t})$ uniformly and hence, \\
 $\bigg| \int g_{1}(x_{t+1})\prod_{i=1}^{N} h(x_{n,t},y^{i}_{n,t})\tau (dx_{t+1}|x_{n,t},\mathbf{u}_{n})- \int g_{1}(x_{t+1})\prod_{i=1}^{N} h(x_{t},y^{i}_{t})\tau (dx_{t+1}|x_{n,t},\mathbf{u}_{n})\bigg| \rightarrow 0$. By the continuity of $\tau$ in $(x,u)$, we get that
 \begin{eqnarray*}
    \nonumber &&\bigg| \int g_{1}(x_{t+1})\prod_{i=1}^{N} h(x_{t},y^{i}_{t})\tau (dx_{t+1}|x_{n,t},\mathbf{u}_{n})-\int g_{1}(x_{t+1})\prod_{i=1}^{N} h(x_{t},y^{i}_{t})\tau (dx_{t+1}|x_{t},\mathbf{u})\bigg|\\
     \nonumber &&=\prod_{i=1}^{N} h(x_{t},y^{i}_{t})\bigg| \int g_{1}(x_{t+1})\tau (dx_{t+1}|x_{n,t},\mathbf{u}_{n})-\int g_{1}(x_{t+1})\tau (dx_{t+1}|x_{t},\mathbf{u})\bigg| \rightarrow 0
 \end{eqnarray*}
Thus we have that $L_{3}^{n}\rightarrow 0$. This, in turn, implies that $\displaystyle \limsup_{n\rightarrow \infty}L^{n}_{2}\leq 4\epsilon \|g\|_{\infty}$. Because $\epsilon$ was arbitrary, it follows that $L^{n}_{2}\rightarrow 0$. 

\qed
\subsection{Proof of Corollary \ref{existence OSDISP}}
It is sufficient to show that a solution to the DCOE (discounted cost optimality equation) exists, where the DCOE is given by
\begin{equation}
J^{*}(Z)=\inf_{f^{1},...,f^{N}}\bigg(\tilde{c}(Z,f^{1},...,f^{N})+\beta \int_{\mathcal{P}(\mathbb{X})}J^{*}(Z_{1})\eta(dZ_{1}|Z,f^{1},...,f^{N})\bigg)
\end{equation}
To this end, it is sufficient to check that measurable selection conditions apply  \cite[Theorem 5.5.3]{lecturenotes}. Note that it follows from the compactness of the action spaces $\mathbb{U}^{1},...,\mathbb{U}^{N}$, and the Young topology introduced in Definition \ref{topology for actions onestep} that the actions $(f^{1},...,f^{N})$ take values in a compact set. By Theorem \ref{theorem5}, $\eta$ is weak-Feller. Thus, it is sufficient to check that 

\[\tilde{c}({Z}_{t},f^{1}_{t},...,f^{N}_{t})=
\int Z_{t} (dx_{t})\prod_{i=1}^{N}f^{i}(du^{i}_{t}|y^{i}_{t})Q^{i}(dy^{i}_{t}|x_{t})c(x_{t},u^{1}_{t},...,u^{N}_{t})\] is bounded and continuous in
$(Z_{t},f^{1}_{t},...,f^{N}_{t})$. The boundedness of $\tilde{c}$ follows from the fact that the cost functional $c$ is bounded. Consider a sequence $(Z_{t}^{n},f^{1}_{n,t},...,f^{N}_{n,t})\rightarrow (Z_{t},f^{1}_{t},...,f^{N}_{t})$. Where, here $Z_{t}^{n}\rightarrow Z_{t}$ weakly and $(f^{1}_{n,t},...,f^{N}_{n,t})\rightarrow (f^{1}_{t},...,f^{N}_{t})$ under the topology introduced in Definition \ref{topology for actions onestep}. Let $F:\mathbb{X}\times \mathbb{U}\times \mathbb{Y}\rightarrow \mathbb{R}$ be continuous and bounded. Consider a sequence $x_{n}\rightarrow x$. We get that
\begin{eqnarray*}
    &&\bigg| \int F(x_{n},\mathbf{u},\mathbf{y}) \prod_{i=1}^{N} h(x_{n},y^{i})\psi(dy^{i})f^{i}_{n}(du^{i}|y^{i})  -\int F(x,\mathbf{u},\mathbf{y}) \prod_{i=1}^{N} h(x,y^{i})\psi(dy^{i})f^{i}(du^{i}|y^{i})  \bigg|\\
    &\leq& \bigg| \int F(x_{n},\mathbf{u},\mathbf{y}) \prod_{i=1}^{N} h(x_{n},y^{i})\psi(dy^{i})f^{i}_{n}(du^{i}|y^{i})  - \int F(x_{n},\mathbf{u},\mathbf{y}) \prod_{i=1}^{N} h(x_{n},y^{i})\psi(dy^{i})f^{i}(du^{i}|y^{i})\bigg| \\
    && \quad +\bigg|\int F(x_{n},\mathbf{u},\mathbf{y}) \prod_{i=1}^{N} h(x_{n},y^{i})\psi(dy^{i})f^{i}(du^{i}|y^{i})  -\int F(x,\mathbf{u},\mathbf{y}) \prod_{i=1}^{N} h(x,y^{i})\psi(dy^{i})f^{i}(du^{i}|y^{i})  \bigg|
\end{eqnarray*}
By DCT, we have that \\
$\bigg|\int F(x_{n},\mathbf{u},\mathbf{y}) \prod_{i=1}^{N} h(x_{n},y^{i})\psi(dy^{i})f^{i}(du^{i}|y^{i}) -\int F(x,\mathbf{u},\mathbf{y}) \prod_{i=1}^{N} h(x,y^{i})\psi(dy^{i})f^{i}(du^{i}|y^{i})  \bigg|\rightarrow 0$. Let $\epsilon>0$. There exists some compact set $K_{1}$ such that $\psi(K_{1}^{C})\leq \epsilon$. By Lusin's Theorem, there exists a further compact set $K$ such that the restriction of $h$ to $\mathbb{X}\times K$ is continuous and $\psi(K)\leq 2\epsilon$. Let $S=K^{N}\times \prod_{i=1}^{N}\mathbb{U}^{i}$ Thus, we get
\begin{eqnarray*}
    &&\bigg| \int F(x_{n},\mathbf{u},\mathbf{y}) \prod_{i=1}^{N} h(x_{n},y^{i})\psi(dy^{i})f^{i}_{n}(du^{i}|y^{i}) - \int F(x_{n},\mathbf{u},\mathbf{y}) \prod_{i=1}^{N} h(x_{n},y^{i})\psi(dy^{i})f^{i}(du^{i}|y^{i})\bigg| \\
    &&\leq 4\epsilon\|F\|_{\infty}+\bigg| \int_{S} F(x_{n},\mathbf{u},\mathbf{y}) \prod_{i=1}^{N} h(x_{n},y^{i})\psi(dy^{i})f^{i}_{n}(du^{i}|y^{i}) - \int_{S} F(x_{n},\mathbf{u},\mathbf{y}) \prod_{i=1}^{N} h(x_{n},y^{i})\psi(dy^{i})f^{i}(du^{i}|y^{i})\bigg|
\end{eqnarray*}
By Stone-Weierstrass Theorem, there exists bounded Lipschitz function $G$ such that $\|F.h-G\|_{\infty}\leq\epsilon$. Hence, 
\begin{eqnarray*}  
        && \bigg| \int_{S} F(x_{n},\mathbf{u},\mathbf{y}) \prod_{i=1}^{N} h(x_{n},y^{i})\psi(dy^{i})f^{i}_{n}(du^{i}|y^{i})  - \int_{S} F(x_{n},\mathbf{u},\mathbf{y}) \prod_{i=1}^{N} h(x_{n},y^{i})\psi(dy^{i})f^{i}(du^{i}|y^{i})\bigg|\\
    &&\leq 2\epsilon + \bigg| \int_{S} G\prod_{i=1}^{N} \psi(dy^{i})f^{i}_{n}(du^{i}|y^{i}) - \int_{S} G \prod_{i=1}^{N} \psi(dy^{i})f^{i}(du^{i}|y^{i})\bigg|  \\
    &&\leq 2\epsilon+\|G\|_{BL}\rho(\prod_{i=1}^{N} \psi(dy^{i})f^{i}_{n}(du^{i}|y^{i}),\prod_{i=1}^{N}  \psi(dy^{i})f^{i}(du^{i}|y^{i}))\rightarrow 2\epsilon
\end{eqnarray*}
 Because $\epsilon$ was arbitrary, it follows that
\[\bigg| \int F(x_{n},\mathbf{u},\mathbf{y}) \prod_{i=1}^{N} h(x_{n},y^{i})\psi(dy^{i})f^{i}_{n}(du^{i}|y^{i}) -\int F(x,\mathbf{u},\mathbf{y}) \prod_{i=1}^{N} h(x,y^{i})\psi(dy^{i})f^{i}(du^{i}|y^{i})  \bigg|\rightarrow 0\]
Since, $Z_{n}\rightarrow Z$ weakly, we then get by the generalized DCT \cite[Theorem 3.5]{Lan81} \cite[Theorem 3.5]{serfozo1982convergence} that \[\int F Z_{t}^{n} (dx_{t})\prod_{i=1}^{N}f^{i}_{n,t}(du^{i}_{t}|y^{i}_{t})Q^{i}(dy^{i}_{t}|x_{t})\rightarrow \int F Z_{t} (dx_{t})\prod_{i=1}^{N}f^{i}(du^{i}_{t}|y^{i}_{t})Q^{i}(dy^{i}_{t}|x_{t}).\]

Thus, we have that \[Z_{t}^{n} (dx_{t})\prod_{i=1}^{N}f^{i}_{n,t}(du^{i}_{t}|y^{i}_{t})Q^{i}(dy^{i}_{t}|x_{t})\rightarrow Z_{t} (dx_{t})\prod_{i=1}^{N}f^{i}(du^{i}_{t}|y^{i}_{t})Q^{i}(dy^{i}_{t}|x_{t})\] weakly. Hence, since $c$ is continuous and bounded the desired result follows.
\qed
\subsection{Proof of Corollary \ref{existence KSDISP}}
Similar to the proof of Corollary \ref{existence OSDISP}, it is sufficient to check that 

\begin{eqnarray*}
   \nonumber &&\Tilde{c}(\pi_{q},a^{1}_{q},...,a^{N}_{q})= \\
   \nonumber &&\int \pi_{q}(dx_{qK})\sum_{t=qK}^{t=(q+1)K-1} \beta^{t-qK} (\prod_{i=1}^{N}Q^{i}(dy^{i}_{qK}|x_{qK})f^{i}_{qK}(du^{i}_{qK}|y^{i}_{qK})) 
    \tau(dx_{qK+1}|x_{qK},\mathbf{u}_{qK})\\
    &&...\tau(dx_{t}|x_{t-1},\mathbf{u}_{t-1})(\prod_{i=1}^{N}Q^{i}(dy^{i}_{t}|x_{t-1})f^{i}_{t}(du^{i}_{t }|y^{i}_{[qK,qK+r]},u^{i}_{[qK,qK+r-1]}))c(x_{t},\mathbf{u}_{t})
\end{eqnarray*} 
is continuous. Consider a sequence $(\pi_{q}^{n},a_{q,n})\rightarrow (\pi_{q},a_{q})$. When $t=qk$, the continuity of the $qk$ component of $\tilde{c}$ follows from the proof of Corollary \ref{existence OSDISP}. Let $t\in \{qK+1,...,(q+1)K-1\}$. For all $x_{qK}$ we have that
\begin{eqnarray*}
    &&E_{n}(x_{qK}):=\bigg|\int (\prod_{i=1}^{N}Q^{i}(dy^{i}_{qK}|x_{qK})f^{i}_{qK}(du^{i}_{qK}|y^{i}_{qK})) 
    \tau(dx_{qK+1}|x_{qK},\mathbf{u}_{qK})...\tau(dx_{t}|x_{t-1},\mathbf{u}_{t-1})\\
    &&\quad (\prod_{i=1}^{N}Q^{i}(dy^{i}_{t}|x_{t-1})f^{i}_{t}(du^{i}_{t }|y^{i}_{[qK,qK+r]},u^{i}_{[qK,qK+r-1]}))c(x_{t},\mathbf{u}_{t})\\
    &&\quad -\int (\prod_{i=1}^{N}Q^{i}(dy^{i}_{qK}|x_{qK})f^{i}_{qK,n}(du^{i}_{qK}|y^{i}_{qK})) 
    \tau(dx_{qK+1}|x_{qK},\mathbf{u}_{qK})...\tau(dx_{t}|x_{t-1},\mathbf{u}_{t-1})\\
    &&\quad (\prod_{i=1}^{N}Q^{i}(dy^{i}_{t}|x_{t-1})f^{i}_{n,t}(du^{i}_{t }|y^{i}_{[qK,qK+r]},u^{i}_{[qK,qK+r-1]}))c(x_{t},\mathbf{u}_{t})\bigg|\\
    &&=\bigg|\int \prod_{i=1}^{N}\prod_{t'=qK}^{t}f^{i}_{n,t'}(du^{i}_{t' }|y^{i}_{[qK,qK+r']},u^{i}_{[qK,qK+r'-1]}))\psi(dy^{i}_{t'})\\
    &&\int \big(\prod_{t'=qK}^{t-1}\tau(dx_{t'+1}|x_{t'},\mathbf{u}_{t'}) (\prod_{i=1}^{N} h(x_{t'},y^{i}_{t'}))\big)\prod_{i=1}^{N}h(x_{t},y^{i}_{t}) c(x_{t},\mathbf{u}_{t})\\
    && \quad - \int \prod_{i=1}^{N}\prod_{t'=qK}^{t}f^{i}_{t'}(du^{i}_{t' }|y^{i}_{[qK,qK+r']},u^{i}_{[qK,qK+r'-1]}))\psi(dy^{i}_{t'})\\
    &&\int \big(\prod_{t'=qK}^{t-1}\tau(dx_{t'+1}|x_{t'},\mathbf{u}_{t'}) (\prod_{i=1}^{N} h(x_{t'},y^{i}_{t'}))\big)\prod_{i=1}^{N}h(x_{t},y^{i}_{t}) c(x_{t},\mathbf{u}_{t}) \bigg|
\end{eqnarray*}
By Lusin's Theorem, since $c$ is bounded, the contribution of the set where $h$ may not be continuous can be made arbitrarily small. Thus, one can without loss of generality assume that $h$ is continuous in $(x,u)$. Thus, throughout the rest of the proof we will assume that $h$ is continuous. We have that $\prod_{i=1}^{N}\prod_{t'=qK}^{t}f^{i}_{n,t'}(du^{i}_{t' }|y^{i}_{[qK,qK+r']},u^{i}_{[qK,qK+r'-1]}))\psi(dy^{i}_{t'})\rightarrow$\\
$\prod_{i=1}^{N}\prod_{t'=qK}^{t}f^{i}_{t'}(du^{i}_{t' }|y^{i}_{[qK,qK+r']},u^{i}_{[qK,qK+r'-1]}))\psi(dy^{i}_{t'})$ weakly and \\
$\int \big(\prod_{t'=qK}^{t-1}\tau(dx_{t'+1}|x_{t'},\mathbf{u}_{t'}) (\prod_{i=1}^{N} h(x_{t'},y^{i}_{t'}))\big)\prod_{i=1}^{N}h(x_{t},y^{i}_{t}) c(x_{t},\mathbf{u}_{t})$ is continuous in $(\mathbf{y}_{[qK,t]},\mathbf{u}_{[qK,t]})$ and bounded. It, thus, follows that $E_{n}(x_{qK})\rightarrow 0$ for all $x_{qK}$.
Consider a sequence $x_{qK,n}\rightarrow x_{qK}$. Then, we have that for every $\mathbf{u}_{qK}$
\begin{eqnarray*}
     && \bigg|\int \tau(dx_{qK+1}|x_{qK,n},\mathbf{u}_{qK})\int \prod_{i=1}^{N}\prod_{t'=qK}^{t}f^{i}_{n,t'}(du^{i}_{t' }|y^{i}_{[qK,qK+r']},u^{i}_{[qK,qK+r'-1]}))\psi(dy^{i}_{t'})\\
    &&\int \big(\prod_{t'=qK+1}^{t-1}\tau(dx_{t'+1}|x_{t'},\mathbf{u}_{t'}) (\prod_{i=1}^{N} h(x_{t'},y^{i}_{t'}))\big)\prod_{i=1}^{N}h(x_{t},y^{i}_{t}) c(x_{t},\mathbf{u}_{t})\\
    && \quad - \int \tau(dx_{qK+1}|x_{qK},\mathbf{u}_{qK})\int \prod_{i=1}^{N}\prod_{t'=qK}^{t}f^{i}_{t'}(du^{i}_{t' }|y^{i}_{[qK,qK+r']},u^{i}_{[qK,qK+r'-1]}))\psi(dy^{i}_{t'})\\
    &&\int \big(\prod_{t'=qK+1}^{t-1}\tau(dx_{t'+1}|x_{t'},\mathbf{u}_{t'}) (\prod_{i=1}^{N} h(x_{t'},y^{i}_{t'}))\big)\prod_{i=1}^{N}h(x_{t},y^{i}_{t}) c(x_{t},\mathbf{u}_{t}) \bigg|\\
    &&\leq   \int \tau(dx_{qK+1}|x_{qK,n},\mathbf{u}_{qK})\bigg|\int \prod_{i=1}^{N}\prod_{t'=qK}^{t}f^{i}_{n,t'}(du^{i}_{t' }|y^{i}_{[qK,qK+r']},u^{i}_{[qK,qK+r'-1]}))\psi(dy^{i}_{t'})\\
    &&\int \big(\prod_{t'=qK+1}^{t-1}\tau(dx_{t'+1}|x_{t'},\mathbf{u}_{t'}) (\prod_{i=1}^{N} h(x_{t'},y^{i}_{t'}))\big)\prod_{i=1}^{N}h(x_{t},y^{i}_{t}) c(x_{t},\mathbf{u}_{t})\\
    && \quad - \int \prod_{i=1}^{N}\prod_{t'=qK}^{t}f^{i}_{t'}(du^{i}_{t' }|y^{i}_{[qK,qK+r']},u^{i}_{[qK,qK+r'-1]}))\psi(dy^{i}_{t'})\\
    &&\int \big(\prod_{t'=qK+1}^{t-1}\tau(dx_{t'+1}|x_{t'},\mathbf{u}_{t'}) (\prod_{i=1}^{N} h(x_{t'},y^{i}_{t'}))\big)\prod_{i=1}^{N}h(x_{t},y^{i}_{t}) c(x_{t},\mathbf{u}_{t}) \bigg|\\
    &&\quad + 2\|c\|_{\infty}\|\tau(dx_{qK+1}|x_{qK,n},\mathbf{u}_{qK})-\tau(dx_{qK+1}|x_{qK},\mathbf{u}_{qK})\|_{TV} 
\end{eqnarray*}
Since the right-hand side of the last inequality converges to zero as $n\rightarrow \infty$, the desired result then follows by an application of the generalized DCT \cite[Theorem 3.5]{Lan81}  \cite[Theorem 3.5]{serfozo1982convergence}. 
\qed



\begin{thebibliography}{72}
\ifx \bisbn   \undefined \def \bisbn  #1{ISBN #1}\fi
\ifx \binits  \undefined \def \binits#1{#1}\fi
\ifx \bauthor  \undefined \def \bauthor#1{#1}\fi
\ifx \batitle  \undefined \def \batitle#1{#1}\fi
\ifx \bjtitle  \undefined \def \bjtitle#1{#1}\fi
\ifx \bvolume  \undefined \def \bvolume#1{\textbf{#1}}\fi
\ifx \byear  \undefined \def \byear#1{#1}\fi
\ifx \bissue  \undefined \def \bissue#1{#1}\fi
\ifx \bfpage  \undefined \def \bfpage#1{#1}\fi
\ifx \blpage  \undefined \def \blpage #1{#1}\fi
\ifx \burl  \undefined \def \burl#1{\textsf{#1}}\fi
\ifx \doiurl  \undefined \def \doiurl#1{\url{https://doi.org/#1}}\fi
\ifx \betal  \undefined \def \betal{\textit{et al.}}\fi
\ifx \binstitute  \undefined \def \binstitute#1{#1}\fi
\ifx \binstitutionaled  \undefined \def \binstitutionaled#1{#1}\fi
\ifx \bctitle  \undefined \def \bctitle#1{#1}\fi
\ifx \beditor  \undefined \def \beditor#1{#1}\fi
\ifx \bpublisher  \undefined \def \bpublisher#1{#1}\fi
\ifx \bbtitle  \undefined \def \bbtitle#1{#1}\fi
\ifx \bedition  \undefined \def \bedition#1{#1}\fi
\ifx \bseriesno  \undefined \def \bseriesno#1{#1}\fi
\ifx \blocation  \undefined \def \blocation#1{#1}\fi
\ifx \bsertitle  \undefined \def \bsertitle#1{#1}\fi
\ifx \bsnm \undefined \def \bsnm#1{#1}\fi
\ifx \bsuffix \undefined \def \bsuffix#1{#1}\fi
\ifx \bparticle \undefined \def \bparticle#1{#1}\fi
\ifx \barticle \undefined \def \barticle#1{#1}\fi
\bibcommenthead
\ifx \bconfdate \undefined \def \bconfdate #1{#1}\fi
\ifx \botherref \undefined \def \botherref #1{#1}\fi
\ifx \url \undefined \def \url#1{\textsf{#1}}\fi
\ifx \bchapter \undefined \def \bchapter#1{#1}\fi
\ifx \bbook \undefined \def \bbook#1{#1}\fi
\ifx \bcomment \undefined \def \bcomment#1{#1}\fi
\ifx \oauthor \undefined \def \oauthor#1{#1}\fi
\ifx \citeauthoryear \undefined \def \citeauthoryear#1{#1}\fi
\ifx \endbibitem  \undefined \def \endbibitem {}\fi
\ifx \bconflocation  \undefined \def \bconflocation#1{#1}\fi
\ifx \arxivurl  \undefined \def \arxivurl#1{\textsf{#1}}\fi
\csname PreBibitemsHook\endcsname

\bibitem[\protect\citeauthoryear{D.S.Bernstein and S.Zilberstein}{2022}]{bernstein2022comp-maor}
\begin{botherref}
\oauthor{\bsnm{D.S.Bernstein}, \binits{N.I.} \bsuffix{R.Givan}},
\oauthor{\bsnm{S.Zilberstein}}:
The complexity of decentralized control of markov decision processes.
Mathematics of Operations Research
(2022)
\end{botherref}
\endbibitem

\bibitem[\protect\citeauthoryear{Witsenhausen}{1971}]{wit71}
\begin{barticle}
\bauthor{\bsnm{Witsenhausen}, \binits{H.S.}}:
\batitle{Separation of estimation and control for discrete time systems}.
\bjtitle{Proceedings of the IEEE}
\bvolume{59},
\bfpage{1557}--\blpage{1566}
(\byear{1971})
\end{barticle}
\endbibitem

\bibitem[\protect\citeauthoryear{H.S.Witsenhausen}{1968}]{Witsenhausen1968counterexample}
\begin{botherref}
\oauthor{\bsnm{H.S.Witsenhausen}}:
A counterexample in stochastic optimum control.
SIAM J. Control
(1968)
\end{botherref}
\endbibitem

\bibitem[\protect\citeauthoryear{Y\"uksel and Ba\c{s}ar}{2024}]{YukselBasarBook24}
\begin{bbook}
\bauthor{\bsnm{Y\"uksel}, \binits{S.}},
\bauthor{\bsnm{Ba\c{s}ar}, \binits{T.}}:
\bbtitle{Stochastic Teams, Games, and Control Under Information Constraints}.
\bpublisher{Springer},
\blocation{Cham}
(\byear{2024})
\end{bbook}
\endbibitem

\bibitem[\protect\citeauthoryear{Athans}{1974}]{Athans}
\begin{bchapter}
\bauthor{\bsnm{Athans}, \binits{M.}}:
\bctitle{Survey of decentralized control methods}.
\bpublisher{3rd NBER/FRB Workshop on Stochastic Control},
\blocation{Washington D.C}
(\byear{1974})
\end{bchapter}
\endbibitem

\bibitem[\protect\citeauthoryear{Varaiya and Walrand}{1978}]{VaraiyaWalrand}
\begin{barticle}
\bauthor{\bsnm{Varaiya}, \binits{P.}},
\bauthor{\bsnm{Walrand}, \binits{J.}}:
\batitle{On delayed-sharing patterns}.
\bjtitle{IEEE Transactions on Automatic Control}
\bvolume{23},
\bfpage{443}--\blpage{445}
(\byear{1978})
\end{barticle}
\endbibitem

\bibitem[\protect\citeauthoryear{Witsenhausen}{1971}]{WitsenhausenSeparation}
\begin{barticle}
\bauthor{\bsnm{Witsenhausen}, \binits{H.S.}}:
\batitle{Separation of estimation and control for discrete time systems}.
\bjtitle{Proceedings of the IEEE}
\bvolume{59},
\bfpage{1557}--\blpage{1566}
(\byear{1971})
\end{barticle}
\endbibitem

\bibitem[\protect\citeauthoryear{Yoshikawa}{1975}]{yos75}
\begin{barticle}
\bauthor{\bsnm{Yoshikawa}, \binits{T.}}:
\batitle{Dynamic programming approach to decentralized control problems}.
\bjtitle{IEEE Transactions on Automatic Control}
\bvolume{20},
\bfpage{796}--\blpage{797}
(\byear{1975})
\end{barticle}
\endbibitem

\bibitem[\protect\citeauthoryear{Chong and Athans}{1976}]{ChongAthans}
\begin{barticle}
\bauthor{\bsnm{Chong}, \binits{C.Y.}},
\bauthor{\bsnm{Athans}, \binits{M.}}:
\batitle{On the periodic coordination of linear stochastic systems}.
\bjtitle{Automatica}
\bvolume{12},
\bfpage{321}--\blpage{335}
(\byear{1976})
\end{barticle}
\endbibitem

\bibitem[\protect\citeauthoryear{Aicardi et~al.}{1987}]{AicardiDavoli}
\begin{barticle}
\bauthor{\bsnm{Aicardi}, \binits{M.}},
\bauthor{\bsnm{Davoli}, \binits{F.}},
\bauthor{\bsnm{Minciardi}, \binits{R.}}:
\batitle{Decentralized optimal control of {M}arkov chains with a common past information set}.
\bjtitle{IEEE Transactions on Automatic Control}
\bvolume{32},
\bfpage{1028}--\blpage{1031}
(\byear{1987})
\end{barticle}
\endbibitem

\bibitem[\protect\citeauthoryear{Y\"uksel}{2009}]{YukTAC09}
\begin{barticle}
\bauthor{\bsnm{Y\"uksel}, \binits{S.}}:
\batitle{Stochastic nestedness and the belief sharing information pattern}.
\bjtitle{IEEE Transactions on Automatic Control}
\bvolume{54},
\bfpage{2773}--\blpage{2786}
(\byear{2009})
\end{barticle}
\endbibitem

\bibitem[\protect\citeauthoryear{Lamperski and Lessard}{2015}]{lamperski2015optimal}
\begin{barticle}
\bauthor{\bsnm{Lamperski}, \binits{A.}},
\bauthor{\bsnm{Lessard}, \binits{L.}}:
\batitle{Optimal decentralized state-feedback control with sparsity and delays}.
\bjtitle{Automatica}
\bvolume{58},
\bfpage{143}--\blpage{151}
(\byear{2015})
\end{barticle}
\endbibitem

\bibitem[\protect\citeauthoryear{Nayyar et~al.}{2011}]{NayyarMahajanTeneketzis}
\begin{barticle}
\bauthor{\bsnm{Nayyar}, \binits{A.}},
\bauthor{\bsnm{Mahajan}, \binits{A.}},
\bauthor{\bsnm{Teneketzis}, \binits{D.}}:
\batitle{Optimal control strategies in delayed sharing information structures}.
\bjtitle{IEEE Transactions on Automatic Control}
\bvolume{56},
\bfpage{1606}--\blpage{1620}
(\byear{2011})
\end{barticle}
\endbibitem

\bibitem[\protect\citeauthoryear{Nayyar et~al.}{2013}]{NayyarBookChapter}
\begin{bchapter}
\bauthor{\bsnm{Nayyar}, \binits{A.}},
\bauthor{\bsnm{Mahajan}, \binits{A.}},
\bauthor{\bsnm{Teneketzis}, \binits{D.}}:
\bctitle{The common-information approach to decentralized stochastic control}.
In: \beditor{\bsnm{Como}, \binits{G.}},
\beditor{\bsnm{Bernhardsson}, \binits{B.}},
\beditor{\bsnm{Rantzer}, \binits{A.}} (eds.)
\bbtitle{Information and Control in Networks},
pp. \bfpage{123}--\blpage{156}.
\bpublisher{Springer},
\blocation{Berlin, Heidelberg}
(\byear{2013})
\end{bchapter}
\endbibitem

\bibitem[\protect\citeauthoryear{Nayyar and Teneketzis}{2019}]{nayyar2019common}
\begin{barticle}
\bauthor{\bsnm{Nayyar}, \binits{A.}},
\bauthor{\bsnm{Teneketzis}, \binits{D.}}:
\batitle{Common knowledge and sequential team problems}.
\bjtitle{IEEE Transactions on Automatic Control}
\bvolume{64}(\bissue{12}),
\bfpage{5108}--\blpage{5115}
(\byear{2019})
\end{barticle}
\endbibitem

\bibitem[\protect\citeauthoryear{H.Tavafoghi and D.Teneketzis}{2021}]{tavafoghi2018unified}
\begin{barticle}
\bauthor{\bsnm{H.Tavafoghi}, \binits{Y.O.}},
\bauthor{\bsnm{D.Teneketzis}}:
\batitle{A unified approach to dynamic decision problems with asymmetric information: Non-strategic agents}.
\bjtitle{IEEE Transactions on Automatic Control}
\bvolume{67}(\bissue{3}),
\bfpage{1105}--\blpage{1119}
(\byear{2021})
\end{barticle}
\endbibitem

\bibitem[\protect\citeauthoryear{V.Subramanian and H.Kao}{2022}]{Subramanian-2021-common}
\begin{botherref}
\oauthor{\bsnm{V.Subramanian}},
\oauthor{\bsnm{H.Kao}}:
Common information based approximate state representations in multi-agent reinforcement learning.
Proccedings of the 25\textsuperscript{th} International Conference on Artificial Intelligence and Statistics
(2022)
\end{botherref}
\endbibitem

\bibitem[\protect\citeauthoryear{Witsenhausen}{1973}]{WitsenStandard}
\begin{barticle}
\bauthor{\bsnm{Witsenhausen}, \binits{H.S.}}:
\batitle{A standard form for sequential stochastic control}.
\bjtitle{Mathematical Systems Theory}
\bvolume{7},
\bfpage{5}--\blpage{11}
(\byear{1973})
\end{barticle}
\endbibitem

\bibitem[\protect\citeauthoryear{Y\"uksel}{2020}]{YukselWitsenStandardArXiv}
\begin{barticle}
\bauthor{\bsnm{Y\"uksel}, \binits{S.}}:
\batitle{A universal dynamic program and refined existence results for decentralized stochastic control}.
\bjtitle{SIAM Journal on Control and Optimization}
\bvolume{58}(\bissue{5}),
\bfpage{2711}--\blpage{2739}
(\byear{2020})
\end{barticle}
\endbibitem

\bibitem[\protect\citeauthoryear{J.Dibangoye et~al.}{2016}]{dibangoye2016}
\begin{barticle}
\bauthor{\bsnm{J.Dibangoye}},
\bauthor{\bsnm{C.Amato}},
\bauthor{\bsnm{O.Buffet}},
\bauthor{\bsnm{F.Charpillet}}:
\batitle{Optimally solving dec-pomdps as continuous-state mdps}.
\bjtitle{Journal of Artificial Intelligence Research}
\bvolume{55},
\bfpage{443}--\blpage{497}
(\byear{2016})
\end{barticle}
\endbibitem

\bibitem[\protect\citeauthoryear{Charalambous}{2016}]{charalambous2016decentralized}
\begin{barticle}
\bauthor{\bsnm{Charalambous}, \binits{C.D.}}:
\batitle{Decentralized optimality conditions of stochastic differential decision problems via {G}irsanov's measure transformation}.
\bjtitle{Mathematics of Control, Signals, and Systems}
\bvolume{28}(\bissue{3}),
\bfpage{1}--\blpage{55}
(\byear{2016})
\end{barticle}
\endbibitem

\bibitem[\protect\citeauthoryear{Charalambous and Ahmed}{2017}]{charalambous2017centralizedI}
\begin{barticle}
\bauthor{\bsnm{Charalambous}, \binits{C.D.}},
\bauthor{\bsnm{Ahmed}, \binits{N.U.}}:
\batitle{Centralized versus decentralized optimization of distributed stochastic differential decision systems with different information structures-part i: {A} general theory}.
\bjtitle{IEEE Transactions on Automatic Control}
\bvolume{62}(\bissue{3}),
\bfpage{1194}--\blpage{1209}
(\byear{2017})
\end{barticle}
\endbibitem

\bibitem[\protect\citeauthoryear{Ooi et~al.}{1997}]{OoiWornell}
\begin{barticle}
\bauthor{\bsnm{Ooi}, \binits{J.}},
\bauthor{\bsnm{Verbout}, \binits{S.}},
\bauthor{\bsnm{Ludwig}, \binits{J.}},
\bauthor{\bsnm{Wornell}, \binits{G.}}:
\batitle{A separation theorem for periodic sharing information patterns in decentralized control}.
\bjtitle{IEEE Transactions on Automatic Control}
\bvolume{42},
\bfpage{1546}--\blpage{1550}
(\byear{1997})
\end{barticle}
\endbibitem

\bibitem[\protect\citeauthoryear{Saldi}{2023}]{Saldi2023Commoninformationapproachstaticteamproblems}
\begin{botherref}
\oauthor{\bsnm{Saldi}, \binits{N.}}:
Common information approach for static team problems with polish spaces and existence of optimal policies.
arXiv preprint arXiv:2309.07571
(2023)
\end{botherref}
\endbibitem

\bibitem[\protect\citeauthoryear{Mao et~al.}{2023}]{mao2023decentralized}
\begin{barticle}
\bauthor{\bsnm{Mao}, \binits{W.}},
\bauthor{\bsnm{Zhang}, \binits{K.}},
\bauthor{\bsnm{Yang}, \binits{Z.}},
\bauthor{\bsnm{Ba{\c{s}}ar}, \binits{T.}}:
\batitle{Decentralized learning of finite-memory policies in dec-pomdps}.
\bjtitle{IFAC-PapersOnLine}
\bvolume{56}(\bissue{2}),
\bfpage{2601}--\blpage{2607}
(\byear{2023})
\end{barticle}
\endbibitem

\bibitem[\protect\citeauthoryear{J.K.Gupta et~al.}{2017}]{gupta2017cooperative}
\begin{botherref}
\oauthor{\bsnm{J.K.Gupta}},
\oauthor{\bsnm{M.Egorov}},
\oauthor{\bsnm{M.Kochenderfer}}:
Cooperative multi-agent control using deep reinforcement learning.
Proccedings of the sixteenth International conference on Autonomous Agents and Multiagent Systems,
66--83
(2017)
\end{botherref}
\endbibitem

\bibitem[\protect\citeauthoryear{Hu and Foerster}{2019}]{hu2019simplified}
\begin{botherref}
\oauthor{\bsnm{Hu}, \binits{H.}},
\oauthor{\bsnm{Foerster}, \binits{J.N.}}:
Simplified action decoder for deep multi-agent reinforcement learning.
arXiv preprint arXiv:1912.02288
(2019)
\end{botherref}
\endbibitem

\bibitem[\protect\citeauthoryear{Wu et~al.}{2013}]{wu2013monte}
\begin{botherref}
\oauthor{\bsnm{Wu}, \binits{F.}},
\oauthor{\bsnm{Zilberstein}, \binits{S.}},
\oauthor{\bsnm{Jennings}, \binits{N.R.}}:
Monte-carlo expectation maximization for decentralized pomdps.
Proccedings of the International joint conference on artificial intelligence
(2013)
\end{botherref}
\endbibitem

\bibitem[\protect\citeauthoryear{Liu et~al.}{2015}]{liu2015stick}
\begin{botherref}
\oauthor{\bsnm{Liu}, \binits{M.}},
\oauthor{\bsnm{Amato}, \binits{C.}},
\oauthor{\bsnm{Liao}, \binits{X.}},
\oauthor{\bsnm{Carin}, \binits{L.}},
\oauthor{\bsnm{How}, \binits{J.P.}}:
Stick-breaking policy learning in dec-pomdps.
arXiv preprint arXiv:1505.00274
(2015)
\end{botherref}
\endbibitem

\bibitem[\protect\citeauthoryear{Kraemer and Banerjee}{2016}]{kraemer2016multi}
\begin{barticle}
\bauthor{\bsnm{Kraemer}, \binits{L.}},
\bauthor{\bsnm{Banerjee}, \binits{B.}}:
\batitle{Multi-agent reinforcement learning as a rehearsal for decentralized planning}.
\bjtitle{Neurocomputing}
\bvolume{190},
\bfpage{82}--\blpage{94}
(\byear{2016})
\end{barticle}
\endbibitem

\bibitem[\protect\citeauthoryear{Saldi et~al.}{2018}]{SaLiYuSpringer}
\begin{bbook}
\bauthor{\bsnm{Saldi}, \binits{N.}},
\bauthor{\bsnm{Linder}, \binits{T.}},
\bauthor{\bsnm{Y\"uksel}, \binits{S.}}:
\bbtitle{Finite Approximations in Discrete-Time Stochastic Control: Quantized Models and Asymptotic Optimality}.
\bpublisher{Springer},
\blocation{Cham}
(\byear{2018})
\end{bbook}
\endbibitem

\bibitem[\protect\citeauthoryear{Kara et~al.}{2019}]{KSYWeakFellerSysCont}
\begin{barticle}
\bauthor{\bsnm{Kara}, \binits{A.D.}},
\bauthor{\bsnm{Saldi}, \binits{N.}},
\bauthor{\bsnm{Y\"uksel}, \binits{S.}}:
\batitle{Weak {F}eller property of non-linear filters}.
\bjtitle{Systems \& Control Letters}
\bvolume{134},
\bfpage{104}--\blpage{512}
(\byear{2019})
\end{barticle}
\endbibitem

\bibitem[\protect\citeauthoryear{Feinberg et~al.}{2016}]{FeKaZg14}
\begin{barticle}
\bauthor{\bsnm{Feinberg}, \binits{E.A.}},
\bauthor{\bsnm{Kasyanov}, \binits{P.O.}},
\bauthor{\bsnm{Zgurovsky}, \binits{M.Z.}}:
\batitle{Partially observable total-cost {M}arkov decision process with weakly continuous transition probabilities}.
\bjtitle{Mathematics of Operations Research}
\bvolume{41}(\bissue{2}),
\bfpage{656}--\blpage{681}
(\byear{2016})
\end{barticle}
\endbibitem

\bibitem[\protect\citeauthoryear{Saldi et~al.}{2017}]{SaYuLi15c}
\begin{barticle}
\bauthor{\bsnm{Saldi}, \binits{N.}},
\bauthor{\bsnm{Y{\"u}ksel}, \binits{S.}},
\bauthor{\bsnm{Linder}, \binits{T.}}:
\batitle{On the asymptotic optimality of finite approximations to {M}arkov decision processes with {B}orel spaces}.
\bjtitle{Mathematics of Operations Research}
\bvolume{42}(\bissue{4}),
\bfpage{945}--\blpage{978}
(\byear{2017})
\end{barticle}
\endbibitem

\bibitem[\protect\citeauthoryear{Kara et~al.}{2023}]{KSYContQLearning}
\begin{botherref}
\oauthor{\bsnm{Kara}, \binits{A.D.}},
\oauthor{\bsnm{Saldi}, \binits{N.}},
\oauthor{\bsnm{Y\"uksel}, \binits{S.}}:
Q-learning for {M}{D}{P}s with general spaces: Convergence and near optimality via quantization under weak continuity.
Journal of Machine Learning Research,
1--34
(2023)
\end{botherref}
\endbibitem

\bibitem[\protect\citeauthoryear{Kara and Y\"uksel}{2023}]{kara2021convergence}
\begin{barticle}
\bauthor{\bsnm{Kara}, \binits{A.D.}},
\bauthor{\bsnm{Y\"uksel}, \binits{S.}}:
\batitle{Convergence of finite memory {Q}-learning for {POMDP}s and near optimality of learned policies under filter stability}.
\bjtitle{Mathematics of Operations Research}
\bvolume{48}(\bissue{4}),
\bfpage{2066}--\blpage{2093}
(\byear{2023})
\end{barticle}
\endbibitem

\bibitem[\protect\citeauthoryear{Wu and Verd\'u}{2011}]{WuVer11}
\begin{bchapter}
\bauthor{\bsnm{Wu}, \binits{Y.}},
\bauthor{\bsnm{Verd\'u}, \binits{S.}}:
\bctitle{Witsenhausen's counterexample: a view from optimal transport theory}.
In: \bbtitle{Proceedings of the IEEE Conference on Decision and Control, Florida, USA},
pp. \bfpage{5732}--\blpage{5737}
(\byear{2011})
\end{bchapter}
\endbibitem

\bibitem[\protect\citeauthoryear{Witsenhausen}{1968}]{wit68}
\begin{barticle}
\bauthor{\bsnm{Witsenhausen}, \binits{H.S.}}:
\batitle{A counterexample in stochastic optimal control}.
\bjtitle{SIAM J. Contr.}
\bvolume{6},
\bfpage{131}--\blpage{147}
(\byear{1968})
\end{barticle}
\endbibitem

\bibitem[\protect\citeauthoryear{Gupta et~al.}{2015}]{gupta2014existence}
\begin{barticle}
\bauthor{\bsnm{Gupta}, \binits{A.}},
\bauthor{\bsnm{Y\"uksel}, \binits{S.}},
\bauthor{\bsnm{Ba\c{s}ar}, \binits{T.}},
\bauthor{\bsnm{Langbort}, \binits{C.}}:
\batitle{On the existence of optimal policies for a class of static and sequential dynamic teams}.
\bjtitle{SIAM Journal on Control and Optimization}
\bvolume{53},
\bfpage{1681}--\blpage{1712}
(\byear{2015})
\end{barticle}
\endbibitem

\bibitem[\protect\citeauthoryear{Saldi}{2020}]{SaldiArXiv2017}
\begin{barticle}
\bauthor{\bsnm{Saldi}, \binits{N.}}:
\batitle{A topology for team policies and existence of optimal team policies in stochastic team theory}.
\bjtitle{IEEE Transactions on Automatic Control}
\bvolume{65}(\bissue{1}),
\bfpage{310}--\blpage{317}
(\byear{2020})
\end{barticle}
\endbibitem

\bibitem[\protect\citeauthoryear{Ethier and Kurtz}{2009}]{ethier2009markov}
\begin{bbook}
\bauthor{\bsnm{Ethier}, \binits{S.N.}},
\bauthor{\bsnm{Kurtz}, \binits{T.G.}}:
\bbtitle{Markov Processes: Characterization and Convergence}.
\bsertitle{Wiley Series in Probability and Statistics},
vol. \bseriesno{282}.
\bpublisher{John Wiley \& Sons}, \blocation{Hoboken, NJ}
(\byear{2009})
\end{bbook}
\endbibitem

\bibitem[\protect\citeauthoryear{Hern{\'a}ndez-Lerma and Lasserre}{1996}]{HernandezLermaMCP}
\begin{bbook}
\bauthor{\bsnm{Hern{\'a}ndez-Lerma}, \binits{O.}},
\bauthor{\bsnm{Lasserre}, \binits{J.B.}}:
\bbtitle{{Discrete-Time Markov Control Processes: Basic Optimality Criteria}}.
\bpublisher{Springer},
\blocation{New York}
(\byear{1996})
\end{bbook}
\endbibitem

\bibitem[\protect\citeauthoryear{Borkar}{1993}]{BorkarRealization}
\begin{barticle}
\bauthor{\bsnm{Borkar}, \binits{V.S.}}:
\batitle{White-noise representations in stochastic realization theory}.
\bjtitle{SIAM J. on Control and Optimization}
\bvolume{31},
\bfpage{1093}--\blpage{1102}
(\byear{1993})
\end{barticle}
\endbibitem

\bibitem[\protect\citeauthoryear{Y{\"u}ksel}{2024}]{yuksel2023borkar}
\begin{barticle}
\bauthor{\bsnm{Y{\"u}ksel}, \binits{S.}}:
\batitle{On {B}orkar and {Y}oung relaxed control topologies and continuous dependence of invariant measures on control policy}.
\bjtitle{SIAM Journal on Control and Optimization}
\bvolume{62}(\bissue{4}),
\bfpage{2367}--\blpage{2386}
(\byear{2024})
\end{barticle}
\endbibitem

\bibitem[\protect\citeauthoryear{Sch\"al}{1975}]{Schal}
\begin{barticle}
\bauthor{\bsnm{Sch\"al}, \binits{M.}}:
\batitle{Conditions for optimality in dynamic programming and for the limit of n-stage optimal policies to be optimal}.
\bjtitle{Z. Wahrscheinlichkeitsth}
\bvolume{32},
\bfpage{179}--\blpage{296}
(\byear{1975})
\end{barticle}
\endbibitem

\bibitem[\protect\citeauthoryear{Balder}{2001}]{balder2001}
\begin{barticle}
\bauthor{\bsnm{Balder}, \binits{E.J.}}:
\batitle{On ws-convergence of product measures}.
\bjtitle{Mathematics of Operations Research}
\bvolume{26}(\bissue{3}),
\bfpage{494}--\blpage{518}
(\byear{2001})
\end{barticle}
\endbibitem

\bibitem[\protect\citeauthoryear{Borkar}{1991}]{Bor91}
\begin{barticle}
\bauthor{\bsnm{Borkar}, \binits{V.S.}}:
\batitle{On extremal solutions to stochastic control problems}.
\bjtitle{Appl. Math. Optim.}
\bvolume{24}(\bissue{1}),
\bfpage{317}--\blpage{330}
(\byear{1991})
\end{barticle}
\endbibitem

\bibitem[\protect\citeauthoryear{Dudley}{2002}]{Dud02}
\begin{bbook}
\bauthor{\bsnm{Dudley}, \binits{R.M.}}:
\bbtitle{Real Analysis and Probability},
\bedition{2nd} edn.
\bpublisher{Cambridge University Press},
\blocation{Cambridge}
(\byear{2002})
\end{bbook}
\endbibitem

\bibitem[\protect\citeauthoryear{Dynkin and Yushkevich}{1979}]{dynkin1979controlled1}
\begin{bbook}
\bauthor{\bsnm{Dynkin}, \binits{E.B.}},
\bauthor{\bsnm{Yushkevich}, \binits{A.A.}}:
\bbtitle{{Controlled Markov Processes}}
vol. \bseriesno{235}.
\bpublisher{Springer},
\blocation{Berlin}
(\byear{1979})
\end{bbook}
\endbibitem

\bibitem[\protect\citeauthoryear{Dubins and Freedman}{1964}]{DubinsFreedman}
\begin{barticle}
\bauthor{\bsnm{Dubins}, \binits{L.}},
\bauthor{\bsnm{Freedman}, \binits{D.}}:
\batitle{Measurable sets of measures}.
\bjtitle{Pacific J. Math.}
\bvolume{14},
\bfpage{1211}--\blpage{1222}
(\byear{1964})
\end{barticle}
\endbibitem

\bibitem[\protect\citeauthoryear{Borkar}{1993}]{Bor93}
\begin{barticle}
\bauthor{\bsnm{Borkar}, \binits{V.S.}}:
\batitle{White-noise representations in stochastic realization theory}.
\bjtitle{SIAM J. Control Optim.}
\bvolume{31}(\bissue{5}),
\bfpage{1093}--\blpage{1102}
(\byear{1993})
\end{barticle}
\endbibitem

\bibitem[\protect\citeauthoryear{Rudin}{1987}]{Rud87}
\begin{bbook}
\bauthor{\bsnm{Rudin}, \binits{W.}}:
\bbtitle{Real and Complex Analysis},
\bedition{3rd} edn.
\bpublisher{McGraw-Hill},
\blocation{New York}
(\byear{1987})
\end{bbook}
\endbibitem

\bibitem[\protect\citeauthoryear{K.Gowrisankaran}{1972}]{Kohur1972measurabilityoftwovariablefunction}
\begin{botherref}
\oauthor{\bsnm{K.Gowrisankaran}}:
Measurability of functions in product spaces.
Proceedings of the American Mathematical Society
(1972)
\end{botherref}
\endbibitem

\bibitem[\protect\citeauthoryear{G.W.Mackey}{1952}]{Mackey1952inducedrepresentationsoflocallycompactgroups}
\begin{botherref}
\oauthor{\bsnm{G.W.Mackey}}:
Induced representations of locally compact groups i.
Annals of Mathematics
(1952)
\end{botherref}
\endbibitem

\bibitem[\protect\citeauthoryear{Grande}{1979}]{Deuxexemplesdefonctionsnonmesurables}
\begin{botherref}
\oauthor{\bsnm{Grande}, \binits{Z.}}:
Deux exemples de fonctions non mesurables.
Colloquim Mathematicum
(1979)
\end{botherref}
\endbibitem

\bibitem[\protect\citeauthoryear{Huang, Wang, and Wang}{2022}]{huang2022general}
\begin{barticle}
\bauthor{\bsnm{Huang}, \binits{J.}}, \bauthor{\bsnm{Wang}, \binits{G.}}, \bauthor{\bsnm{Wang}, \binits{W.}}:
\batitle{A general linear quadratic stochastic control and information value}.
\bjtitle{Journal of Mathematical Analysis and Applications}
\bfpage{126486}
(\byear{2022})
\end{barticle}
\endbibitem


\bibitem[\protect\citeauthoryear{McDonald and Y\"uksel}{2020}]{mcdonald2020exponential}
\begin{botherref}
\oauthor{\bsnm{McDonald}, \binits{C.}},
\oauthor{\bsnm{Y\"uksel}, \binits{S.}}:
Exponential filter stability via {D}obrushin's coefficient.
Electronic Communications in Probability
\textbf{25}
(2020)
\end{botherref}
\endbibitem

\bibitem[\protect\citeauthoryear{Aliprantis and Border}{2006}]{AlBo06}
\begin{bbook}
\bauthor{\bsnm{Aliprantis}, \binits{C.D.}},
\bauthor{\bsnm{Border}, \binits{K.C.}}:
\bbtitle{Infinite Dimensional Analysis}.
\bpublisher{Springer},
\blocation{Berlin}
(\byear{2006})
\end{bbook}
\endbibitem

\bibitem[\protect\citeauthoryear{Bogachev}{2007}]{Bogachev}
\begin{bbook}
\bauthor{\bsnm{Bogachev}, \binits{V.I.}}:
\bbtitle{Measure Theory}.
\bpublisher{Springer},
\blocation{Berlin}
(\byear{2007})
\end{bbook}
\endbibitem

\bibitem[\protect\citeauthoryear{Y\"uksel}{2024}]{lecturenotes}
\begin{bbook}
\bauthor{\bsnm{Y\"uksel}, \binits{S.}}:
\bbtitle{Optimization and Control of Stochastic Systems}.
\bpublisher{Queen's University},
\blocation{kingston, ON}
(\byear{2024}).
\bcomment{Lecture notes, available online}.
\burl{https://mast.queensu.ca/{\texttildelow }yuksel/LectureNotesOnStochasticOptControl.pdf}
\end{bbook}
\endbibitem

\bibitem[\protect\citeauthoryear{Yushkevich}{1976}]{Yus76}
\begin{barticle}
\bauthor{\bsnm{Yushkevich}, \binits{A.A.}}:
\batitle{Reduction of a controlled {M}arkov model with incomplete data to a problem with complete information in the case of {B}orel state and control spaces}.
\bjtitle{Theory Prob. Appl.}
\bvolume{21},
\bfpage{153}--\blpage{158}
(\byear{1976})
\end{barticle}
\endbibitem

\bibitem[\protect\citeauthoryear{Rhenius}{1974}]{Rhe74}
\begin{barticle}
\bauthor{\bsnm{Rhenius}, \binits{D.}}:
\batitle{Incomplete information in {M}arkovian decision models}.
\bjtitle{Ann. Statist.}
\bvolume{2},
\bfpage{1327}--\blpage{1334}
(\byear{1974})
\end{barticle}
\endbibitem

\bibitem[\protect\citeauthoryear{Rishel}{1970}]{rishel1970necessary}
\begin{barticle}
\bauthor{\bsnm{Rishel}, \binits{R.}}:
\batitle{Necessary and sufficient dynamic programming conditions for continuous time stochastic optimal control}.
\bjtitle{SIAM Journal on Control}
\bvolume{8},
\bfpage{559}--\blpage{571}
(\byear{1970})
\end{barticle}
\endbibitem


\bibitem[\protect\citeauthoryear{Nayyar et~al.}{2013}]{NayyarMahajanTeneketzis2}
\begin{barticle}
\bauthor{\bsnm{Nayyar}, \binits{A.}},
\bauthor{\bsnm{Mahajan}, \binits{A.}},
\bauthor{\bsnm{Teneketzis}, \binits{D.}}:
\batitle{Decentralized stochastic control with partial history sharing: A common information approach}.
\bjtitle{IEEE Transactions on Automatic Control}
\bvolume{58},
\bfpage{1644}--\blpage{1658}
(\byear{2013})
\end{barticle}
\endbibitem


\bibitem[\protect\citeauthoryear{Pradhan and Yüksel}{2023}]{pradhanyuksel2023DTApprx}
\begin{barticle}
\bauthor{\bsnm{Pradhan}, \binits{S.}}, \bauthor{\bsnm{Yüksel}, \binits{S.}}:
\batitle{Controlled diffusions under full, partial and decentralized information: Existence of optimal policies and discrete-time approximations}.
\bjtitle{arXiv preprint}
\bvolume{arXiv:2311.03254}
(\byear{2023})
\end{barticle}
\endbibitem










\bibitem[\protect\citeauthoryear{Saldi and Yüksel}{2022}]{saldiyukselGeoInfoStructure}
\begin{barticle}
\bauthor{\bsnm{Saldi}, \binits{N.}},
\bauthor{\bsnm{Y\"uksel}, \binits{S.}}:
\batitle{Geometry of Information Structures, Strategic Measures and Associated Control Topologies}.
\bjtitle{Probability Surveys}
\bvolume{19},
\bfpage{450}--\blpage{532}
(\byear{2022})
\end{barticle}
\endbibitem



\bibitem[\protect\citeauthoryear{Caines, Huang, and Malhamé}{2017}]{caines2018peter}
\begin{barticle}
\bauthor{\bsnm{Caines}, \binits{P.}}, \bauthor{\bsnm{Huang}, \binits{M.}}, \bauthor{\bsnm{Malhamé}, \binits{R.}}:
\batitle{Mean Field Games}.
\bjtitle{Handbook of Dynamic Game Theory}
\bfpage{345}--\blpage{372}
(\byear{2017})
\end{barticle}
\endbibitem

\bibitem[\protect\citeauthoryear{Jackson and Lacker}{2025}]{jackson2023approximately}
\begin{barticle}
\bauthor{\bsnm{Jackson}, \binits{J.}}, \bauthor{\bsnm{Lacker}, \binits{D.}}:
\batitle{Approximately optimal distributed stochastic controls beyond the mean field setting}.
\bjtitle{The Annals of Applied Probability}
\bvolume{35},
\bfpage{251}--\blpage{308}
(\byear{2025})
\end{barticle}
\endbibitem























\bibitem[\protect\citeauthoryear{Kurtaran}{1979}]{kur79}
\begin{barticle}
\bauthor{\bsnm{Kurtaran}, \binits{B.}}:
\batitle{Corrections and extensions to decentralized control with delayed sharing information pattern}.
\bjtitle{IEEE Transactions on Automatic Control}
\bvolume{24},
\bfpage{656}--\blpage{657}
(\byear{1979})
\end{barticle}
\endbibitem

\bibitem[\protect\citeauthoryear{Blackwell and Ryll-Nadrzewski}{1963}]{Blackwell3}
\begin{barticle}
\bauthor{\bsnm{Blackwell}, \binits{D.}},
\bauthor{\bsnm{Ryll-Nadrzewski}, \binits{C.}}:
\batitle{Non-existence of everywhere proper conditional distributions}.
\bjtitle{Annals of Mathematical Statistics}
\bvolume{34},
\bfpage{223}--\blpage{225}
(\byear{1963})
\end{barticle}
\endbibitem

\bibitem[\protect\citeauthoryear{Billingsley}{1995}]{BillingsleyProbMeasure}
\begin{bbook}
\bauthor{\bsnm{Billingsley}, \binits{P.}}:
\bbtitle{Probability and Measure}.
\bpublisher{Wiley},
\blocation{New York}
(\byear{1995})
\end{bbook}
\endbibitem

\bibitem[\protect\citeauthoryear{Langen}{1981}]{Lan81}
\begin{barticle}
\bauthor{\bsnm{Langen}, \binits{H.J.}}:
\batitle{Convergence of dynamic programming models}.
\bjtitle{Mathematics of Operations Research}
\bvolume{6}(\bissue{4}),
\bfpage{493}--\blpage{512}
(\byear{1981})
\end{barticle}
\endbibitem

\bibitem[\protect\citeauthoryear{Serfozo}{1982}]{serfozo1982convergence}
\begin{botherref}
\oauthor{\bsnm{Serfozo}, \binits{R.}}:
Convergence of {L}ebesgue integrals with varying measures.
Sankhy{\=a}: The Indian Journal of Statistics, Series A,
380--402
(1982)
\end{botherref}
\endbibitem

\bibitem[\protect\citeauthoryear{Sch{\"a}l}{1975}]{schal1975dynamic}
\begin{barticle}
\bauthor{\bsnm{Sch{\"a}l}, \binits{M.}}:
\batitle{On dynamic programming: compactness of the space of policies}.
\bjtitle{Stochastic Processes and their Applications}
\bvolume{3}(\bissue{4}),
\bfpage{345}--\blpage{364}
(\byear{1975})
\end{barticle}
\endbibitem

\bibitem[\protect\citeauthoryear{O.Mrani-Zentar and S.Y\"uksel}{2024}]{OMZSY-Dec-centralized-reduction-and-weak-Feller}
\begin{botherref}
\oauthor{\bsnm{O.Mrani-Zentar}},
\oauthor{\bsnm{S.Y\"uksel}}:
Decentralized stochastic control in standard borel spaces: Centralized mdp reductions, near optimality of finite window local information, and q-learning.
arXiv:2408.13828
(2024)
\end{botherref}
\endbibitem

\bibitem[\protect\citeauthoryear{Kara and Y\"uksel}{2020}]{kara2020robustness}
\begin{barticle}
\bauthor{\bsnm{Kara}, \binits{A.D.}},
\bauthor{\bsnm{Y\"uksel}, \binits{S.}}:
\batitle{Robustness to incorrect system models in stochastic control}.
\bjtitle{SIAM Journal on Control and Optimization}
\bvolume{58}(\bissue{2}),
\bfpage{1144}--\blpage{1182}
(\byear{2020})
\end{barticle}
\endbibitem

\bibitem[\protect\citeauthoryear{Billingsley}{1995}]{Bil95}
\begin{bbook}
\bauthor{\bsnm{Billingsley}, \binits{P.}}:
\bbtitle{Probability and Measure},
\bedition{3rd} edn.
\bpublisher{Wiley}, \blocation{Hoboken, NJ}
(\byear{1995})
\end{bbook}
\endbibitem

\bibitem[\protect\citeauthoryear{W.B.~Haskell and D.Kalathil}{2014}]{Haskel2014empiricalvalueiteration}
\begin{botherref}
\oauthor{\bsnm{W.B.~Haskell}, \binits{R.J.}},
\oauthor{\bsnm{D.Kalathil}}:
Empirical value iteration for approximate dynamic programming.
American Control Conference (ACC)
(2014)
\end{botherref}
\endbibitem

\bibitem[\protect\citeauthoryear{Jacka and Vidmar}{2015}]{jacka2015informational}
\begin{barticle}
\bauthor{\bsnm{Jacka}, \binits{S.}}, \bauthor{\bsnm{Vidmar}, \binits{M.}}:
\batitle{On the informational structure in optimal dynamic stochastic control}.
\bjtitle{arXiv preprint}
\bvolume{arXiv:1503.02375}
(\byear{2015})
\end{barticle}
\endbibitem

\bibitem[\protect\citeauthoryear{Elliott}{1977}]{elliott1977optimal}
\begin{barticle}
\bauthor{\bsnm{Elliott}, \binits{R. J.}}:
\batitle{The optimal control of a stochastic system}.
\bjtitle{SIAM Journal on Control and Optimization}
\bvolume{15},
\bfpage{756}--\blpage{778}
(\byear{1977})
\end{barticle}
\endbibitem



\end{thebibliography}

\section*{Statements and Declarations}
\subsection*{Funding}
This research was partially supported by the Natural Sciences and Engineering Research Council
of Canada (NSERC).
\subsection*{Competing interests}
We have no relevant financial or non-financial interests to disclose.
\end{document}